\documentclass[12pt,twoside]{report} 

\usepackage{amscd}
\usepackage{makeidx}
\usepackage{nref}
\usepackage{amssymb}
\usepackage{latexsym}
\usepackage{CUP}
\usepackage{multicol}

\newtheorem{theorem}{Theorem}[section]
\newtheorem{lemma}[theorem]{Lemma}
\newtheorem{defn}[theorem]{Definition}
\newtheorem{corollary}[theorem]{Corollary}
\newtheorem{question}[theorem]{Question}
\newtheorem{conjecture}[theorem]{Conjecture}
\newtheorem{proposition}[theorem]{Proposition}

\let\froggy\label
\def\label#1{\marginlabel{Label: #1}\froggy{#1}}

\def\Index#1{#1\index{#1}}

\mathcode`:="203A

\def\em{\sl}

\long\def\galley#1{$^\dagger$\marginpar{\raggedright\tiny\em#1}}
\long\def\marginlabel#1{\marginpar{\raggedright\tiny\em#1}}

\long\def\marginlabel#1{\relax}
\long\def\galley#1{\relax}

{\catcode`\A=13\gdef A{{\mathchar"7241}}}\mathcode`\A="8000

\def\subsection#1{\par\goodbreak\bigskip\noindent{\bfseries#1.}\par\nobreak\medskip}

\def\Lip{\mname{Lip}}
\def\projconst{projection constant}
\def\dom{\mname{\rm dom}}
\def\rp{reduction property}
\def\rps{reduction properties}
\def\crp{complete \rp}
\def\trp{total \rp}
\def\RA{reduction algebra}
\def\CRA{complete \RA}
\def\TRA{total \RA}
\def\tra{\TRA}
\def\cra{\CRA}
\def\RAs{\RA{}s}
\def\TRAs{\TRA{}s}
\def\tras{\TRA{}s}
\def\V{V}
\def\Va{V}
\def\Vb{W}
\def\Vc{U}
\def\U{U}
\def\W{W}
\def\v{\V}
\def\w{\W}
\def\u{\U}
\def\y{Y}
\def\Y{Y}
\def\cohomology #1(#2,#3){{\cal H}^{#1}(#2; #3)}
\def\hm{Hilbertian module}
\def\hAm#1{Hilbertian $#1$-module}
\def\TbT#1,#2,#3,#4;{
  \left[\begin{array}{cc}#1 & #2 \\ #3 & #4\end{array}\right]
}
\def\TC#1{TC(#1)}
\def\mname#1{\mathop{\hbox{\rm #1}}}

\def\choose#1#2{\bigg(\begin{array}{cc}#1\\#2\end{array}\bigg)}
\def\TxT{\TbT}
\def\Tvec#1,#2;{\left[\begin{array}{c}#1\\#2\end{array}\right]}

\def\twine{\rightsquigarrow}

\def\Gel#1{\hat{#1}}

\def\st{such that}

\def\dint#1.d#2{\int^\oplus\kern-3pt#1\,d#2}
\def\bd#1{{#1}^{**}}
\def\inv#1{{#1}^{-1}}
\def\p#1{\Phi_{#1}}
\def\amp#1#2{#1_{(#2)}}
\def\pma#1#2{#1^{(#2)}}
\def\ran{\mathop{\hbox{ran}}}

\def\sy{similarity}
\def\cbrep{cb-\rep}

\def\ideal{\unlhd}
\def\ideal{\subseteq}
\def\isoc{isomorphic}
\def\simty#1#2{{#1}{#2}\inv{#1}}

\newcommand\iso         {isomorphism}

\def\iff                {if and only if}

\def\spam#1{\mbox{\rm span} \{{#1}\}}
\def\cspam#1{\overline{\mbox{\rm span}} \{{#1}\}}
\def\Gr{\mname{Gr}}
\def\gr{\mname{Gr}}
\def\ss{semisimple}

\newcommand\rep         {represent\-ation}

\newcommand\reps        {represent\-ations}
\newcommand\srep        {\hbox{*-}represent\-ation}
\newcommand\sreps       {\srep{}s}

\newcommand\cp          {\trp}

\newcommand\ad          {\mname{ad}}

\newcommand\lat         {\mname{Lat}}
\def\lats#1{{\rm Lat}_{#1}}
\newcommand\alg         {\mname{Alg}}

\newcommand\oa          {operator algebra}
\newcommand\oas         {operator algebras}
\newcommand\aoa         {amenable \oa}

\def\f(#1,#2){\funct{#1}{#2}}
\def\<#1,#2>{\funct{#1}{#2}}
\def\conj#1^#2{#2 #1 #2^{-1}}

\newcommand\BH          {\B(\H)}

\newcommand\sv{V}

\newcommand\ie          {i.e.}

\newcommand\bad         {bounded approximate diagonal}

\newcommand\fd          {finite dimensional}
\newcommand\cf          {c.f.}
\newcommand\id          {\mname{id}}
\newcommand\cs          {\mbox{$\mbox{C*}\!$}}

\newcommand\N           {\mbox{\bf N}}

\renewcommand\H         {{H}}
\newcommand\Kom           {{\cal K}}
\newcommand\K           {{\cal K}}
\newcommand\KomH           {{\cal K}(\H)}

\newcommand\X           {{\cal X}}
\newcommand\hs          {Hilbert space}
\newcommand\vnalg       {von Neumann algebra}
\newcommand\vnas        {\vnalg{}s}

\newcommand\vna {\vnalg}

\newcommand\csalg       {\cs-algebra}
\newcommand\csalgs      {\csalg{}s}

\newcommand\sstrong     {$\sigma$-strong}

\newcommand\ws          {$\hbox{weak}^*$}

\newcommand\sw          {$\sigma$-weak}
\renewcommand\L         {{\cal L}}
\newcommand\B           {{\cal B}}
\newcommand\sa          {self-adjoint}
\newcommand\C           {\hbox{{\bf C}}}
\newcommand\R           {\hbox{{\bf R}}}
\newcommand\KH          {{\cal K}(\H)}
\def\net#1{\{{#1}\}}
\renewcommand{\R}{{\rm\bf I\!R}} 
\renewcommand{\C}{{\rm\bf \,l\!\!\!C}}
\renewcommand{\N}{{\rm\bf I\!N}}
\renewcommand{\R}{{\!\!\rm\ I\!R}} 
\renewcommand{\C}{\!\!{\rm\ \,l\!\!\!C}}
\renewcommand{\N}{\!\!{\rm\ I\!N}}
\def\R{{\mathbb R}}
\def\C{{\mathbb C}}
\def\N{{\mathbb N}}

\def\Zed{{\mathbb Z}}

\newcommand\closure[1]  {\overline{#1}}

\newcommand\cts {continuous}

\newcommand\bded        {bounded}

\newcommand\eps {\epsilon}
\newcommand\ba          {Banach algebra}
\newcommand\bs          {Banach space}
\newcommand\bm          {Banach module}
\newcommand\bas {\ba{}s}
\newcommand\ra          {\rightarrow}
\renewcommand\ker       {\mname{ker}}
\newcommand\im          {\mname{im}}
\def\hul{\hbox{hul }}

\newcommand\implies    {\mbox{ $\Rightarrow$ }}

\newcommand\bai {bounded approximate identity}
\newcommand\brai {bounded right approximate identity}
\newcommand\blai {bounded left approximate identity}

\newcommand\bais {bounded approximate identities}
\newcommand\lbai {\blai}
\newcommand\rbai {\brai}

\newcommand\proof       {{\bf Proof:}\quad{}}
\def\qed{{\nopagebreak\unskip\nobreak\hfil\penalty50\hskip2em\hbox{}
                        \nobreak\hfil$\Box$\parfillskip=0pt
                        \finalhyphendemerits=0 \par\medskip\smallskip}}
\newcommand\ot          {\otimes}
\newcommand\csnorm      {\cs~norm}

\newcommand\M   {{\cal M}}
\newcommand\cb  {completely bounded}

\def\sangle#1,#2;{\angle(#1,#2)}
\newcommand\supof [1]{\sup\left\{#1\right\}}
\newcommand\infof [1]{\inf\left\{#1 \right\}}
\newcommand\norm   [1]{\|#1\|}
\newcommand\bnorm[1] {\left\|#1\right\|}
\newcommand\projtensor[2]{#1\hat\otimes#2}

\newcommand\funct[2]{\left\langle#1,#2\right\rangle}
\newcommand\BxBM[1]{Banach $#1$-bimodule}
\newcommand\DBxBM[1]{dual \BxBM{#1}}

\newcommand\sweak{\sw}

\def\setof#1{\{#1\}}
\def\clos{\expandafter\closure}

\newbox\closurebox
\def\closure#1{\setbox\closurebox=\hbox{$#1$}
  \ifdim\wd\closurebox>40cm
  \box\closurebox\overline{\vphantom{#1}\phantom{x}}\else
    \overline{\box\closurebox}\fi}
\def\pclos#1{\setbox\closurebox=\hbox{$#1$}
  \ifdim\wd\closurebox>40cm
  \left(\box\closurebox\right)\overline{\vphantom{#1}\phantom{x}}\else
    \overline{\box\closurebox}\fi}
\def\swclos#1{\clos{#1}^{\mbox{\kern.5mm$\scriptstyle\sigma w$}}}

\def\ssclos{\expandafter\ssclosure}

\def\ssclosure#1{\setbox\closurebox=\hbox{$\scriptstyle#1$}
  \ifdim\wd\closurebox>40cm
  \box\closurebox\overline{\vphantom{#1}\phantom{x}}\else
    \overline{\box\closurebox}\fi}
\def\sspclos#1{\setbox\closurebox=\hbox{$\scriptstyle#1$}
  \ifdim\wd\closurebox>40cm
  \left(\box\closurebox\right)\overline{\vphantom{#1}\phantom{x}}\else
    \overline{\box\closurebox}\fi}

\def\clos#1{\overline{#1}}

\def\xspacea{,}
\def\xspaceb{.}
\def\xspacec{;}
\def\xspaced{}
\def\xspacee{\ }
\def\xspacef{\space}
\def\xspace#1{{\def\junk{#1}%
    \ifx\junk\xspacea\else%
    \ifx\junk\xspaceb\else%
    \ifx\junk\xspacec\else%
    \ifx\junk\xspaced\else%
    \ifx\junk\xspacef\else%
    \ifx\junk\xspacee\else\space%
    \fi\fi\fi\fi\fi\fi#1}}

\comment{
\def\thebibliography#1{
        \vskip1in{\bf \begin{center} References \end{center}}
 \list
 {[\arabic{enumi}]}{\settowidth\labelwidth{[#1]}\leftmargin\labelwidth
 \advance\leftmargin\labelsep
 \usecounter{enumi}}
 \def\newblock{\hskip .11em plus .33em minus .07em}
 \sloppy\clubpenalty4000\widowpenalty4000
 \sfcode`\.=1000\relax}

}

\def\longpage{\enlargethispage{\baselineskip}}

\def\title#1    {\vskip1in{\begin{center} {\LARGE\bf{}#1}\\ \today\end{center}}\vskip1in}

\makeindex

\begin{document}

\pagestyle{myheadings}
\markboth{October 1997}{October 1997}

\makeatletter
\def\titlepage{
  \pagenumbering{roman}\pagestyle{plain}
    \thispagestyle{empty}
  {\@rulehead{\noindent Operator algebras \\ with a reduction
    property}}
  \vskip 1in
  \begin{center}\Large
  James A. Gifford \\
  \bigskip
  October 1997
  \end{center}
  \vfill
  \begin{center}\begingroup\scshape
    A thesis submitted for the degree of \\ Doctor of Philosophy \\ of the
    Australian National University
    \endgroup
  \end{center}\newpage
  \hbox to 0pt
  {}\newpage
  \hbox to 0pt{}\vfill
  I certify that the content of this thesis is original work and that
  any material derived from other sources has been duly acknowledged.

  \vskip 5.5cm
  \eject
  
  \subsection{Acknowledgements}

  I am deeply grateful to my supervisor Dr Rick Loy, for his constant
  helpfulness and enthusiasm. Also I would like to thank Dr George
  Willis and the Department of Mathematics at the University of
  Newcastle, for providing the opportunity for many stimulating
  discussions during the course of this research.

  \vfill

  \eject
  }
\makeatother


\titlepage

\def\dealwithsectionbreaks{}


\tableofcontents

\newpage
\hbox{}
\newpage

\pagenumbering{arabic}
\pagestyle{headings}


\markboth{}{}


\chapter*{Introduction}
\addcontentsline{toc}{chapter}{\numberline{}Introduction}

A famous open problem in the theory of \csalgs\ is Kadison's
similarity question: given a \csalg\ $A$, is every \rep\ of $A$
similar to a \srep? Kadison \cite{Kadison} conjectured in the 50s that
this is so for every \csalg\ $A$.

A characterisation of the \csalgs\ for which the above is true can
be given in terms of a `semisimplicity' criterion. A \csalg\ $A$ has
the above property  \iff\ whenever $\theta:A\ra\BH$ is a \rep\ of $A$ on a
\hs\ $\H$, every invariant subspace of $\H$ is complemented by an
invariant subspace. We refer to this semisimplicity property as the {\em
  total reduction property}. In these terms Kadison's conjecture
becomes: every \csalg\ has the total reduction property.

It is not necessary that $A$ be a \sa\ \oa\ for the definition of the
\trp. In this thesis we study the class of not-necessarily \sa\ \oas\
with the \trp.
The \trp\ has strong consequences for the structure of an \oa. The
central theme of this work is that algebras with this property are
`like' \csalgs\ to some extent. In fact, we conjecture that every \oa\
with the \trp\ is \isoc\ to a \csalg.

After a summary of the notations and preliminary results we need in
chapter~\ref{Preliminaries}, we introduce in chapter~\ref{Chapter with
  definition} the definition of the \trp\ as well as some related
notions, and exhibit some examples of algebras with and without this
property.  In chapter~\ref{Chapter with properties} we investigate the
properties of algebras with the \trp.  Many of the results obtained
for these algebras will be familiar from \csalg\ theory.
Chapter~\ref{Chapter with characterisations} is concerned with
applying the machinery obtained to the question of whether algebras
with the \trp\ are \isoc\ to \csalgs. In several special cases we are
able to prove that this is so.  However, in general there is an
obstacle to progress in the form of the transitive algebra problem.
This famous open problem asks: if $A\subseteq\BH$ has no proper
invariant subspaces, is $A$ weakly dense in $\BH$? Despite
intensive work over many decades the transitive algebra problem
remains open. A positive answer to the transitive algebra problem
would, in many cases, allow a positive answer to our conjecture about
algebras with the \trp. The special cases in which we can establish
the result mainly derive from the partial solution to the transitive
algebra problem due to Lomonosov.

\chapter{Preliminaries}
\label{Preliminaries}

Presented here is a brief summary of the notations and standard
results we use throughout this work. As a general guide, the results
we quote about \bas\ can be found in \cite{BonsallDuncan}, and the
results about \csalgs\ and \vnas\ can be found in \cite{DixmierC} and
\cite{DixmierV}.

\subsection{Banach algebras and operator algebras}

A \ba\ is a complex
associative algebra $A$ equipped with a complete norm such that the
multiplication $A\times A\ra A$ is norm \cts.
For every \bs\ $X$, the set $\B(X)$ of bounded linear operators on $X$
forms a \ba\ with the natural algebraic structure and the operator
norm. In the case where $X$ is a \hs\ we obtain an algebra of central
importance to this work. We use the letter $\H$ to denote
a generic \hs, and refer to norm-closed subalgebras of $\BH$ as \oas.
The \oa\ $\B(\H)$ possesses, in addition to its \ba\ structure, an
isometric involution $*$ given by the adjoint operation $a\mapsto
a^*$. 

Operator algebras  have been heavily studied in the
literature. The overwhelming emphasis is on self-adjoint \oas---%
that is, \oas\  which are closed under the
involution. These algebras are \csalgs, and they enjoy a very
well-developed theory.
The theory of non-\sa\ \oas\ is less developed, although
a certain amount of work has been done on the so-called nest algebras and
CSL algebras \cite{Davidson}.  In recent years interest has been
growing in non-\sa\ \oas\ as they relate to the field of
quantised functional analysis \cite{Effros-quantize}, \cite{cb-operators}.

Most of the studies of non-\sa\ \oas\ in the literature consider
properties of a particular realisation of the algebra as a subalgebra
of $\BH$.  In this thesis we examine \oas\ from a more
\ba{}ic perspective.
In general \ba\ theory, the norm is chiefly a provider of topological,
rather than isometric, data. For this reason \bas\ are typically
considered equivalent if they are isomorphic via a bicontinuous
isomorphism. 
The issue does not usually arise for \csalgs, since any \iso\ between
\csalgs\ which preserves the involution is automatically isometric,
and consequently there is a unique \cs~norm implied by the
*-algebraic structure of a \csalg.
There is no such uniqueness of norm for \oas\ without an involution.
For a simple example, let $S\in\BH$ be any invertible operator on a \hs\ 
$\H$. The map $a\mapsto \simty S a$ is a bicontinuous algebra
automorphism of $\BH$, which is isometric \iff\ $S$ is a scalar
multiple of a unitary operator. In all other cases the norm
$\norm{a}_{\rm new}=\norm{\simty S a}$ is a distinct but equivalent
\oa\ norm on $\BH$. 

We will generally consider \oas\ as objects in the category of \bas, and
consider isomorphic \oas\ as equivalent.
Our 
terminology will take the \ba{}ic meaning unless otherwise qualified.
For instance, any homomorphism between \bas\ is intended to be
\cts, and any subalgebra of a \ba\ is intended to be closed.

Similarly, when discussing \bs{}s our terminology will take the \bs\ 
meaning. We will refer to a closed subspace of a \bs\ simply as a
subspace, and refer to \cts\ linear maps between \bs{}s simply as maps
or operators. We use the term `submanifold' to refer to
not-necessarily-closed vector subspaces.

\subsection{Representations and modules}

Let $A$ be a \ba, $X$ be a \bs\ and $\theta:A\ra \B(X)$ be a (\ba)
homomorphism. We say that $\theta$ is a \rep\index{\rep} of $A$ on $X$.  
If $A$ contains an identity we do not require that $\theta(1)=1_X$.
The set of operators $\theta(A)$ is an associative algebra with the
algebra operations inherited from $\B(X)$. However, 
even though we are restricting attention to \cts\ $\theta$, the set
$\theta(A)$ need
not be norm-closed and so is not always a Banach subalgebra of
$\B(X)$.  The closure $\clos{\theta(A)}$ is always a
norm-closed subalgebra of $\B(X)$.

Closely related to \reps\ of $A$ are the Banach modules\index{module} of $A$. 
If $X$ is a \bs\ and $m:A\times X\ra X$ is a \cts\ bilinear map with 
$m(a_1a_2,x)=m(a_1,m(a_2,x))$ for all $a_1,a_2\in A$ and $x\in X$ then
we say that $X$ is a left Banach $A$-module. The module action will
almost always be written as $m(a,x)=a\cdot x$ or even $ax$.

If $\theta:A\ra X$ is a \rep, then it is readily verified that
defining $a\cdot x=\theta(a)x$ yields a left Banach $A$-module action
on $X$.  Conversely, if $X$ is a left Banach $A$-module, then the map
$\theta:A\ra \B(X)$ defined by $\theta(a)(x)=a\cdot x$ is a \rep\ of
$A$. These constructions are mutual inverses. This means that the
concepts of \reps\ and of left Banach $A$-modules are completely
equivalent---the distinction is purely notational. The module notation
is usually neater, but if we need to refer to the operator $x\mapsto
a\cdot x$ it is convenient to use the \rep\ notation and simply write
$\theta(a)$.

A right Banach $A$-module is similarly defined as a \bs\ $X$ with
a \cts\ bilinear map $(x,a)\mapsto x\cdot a$ such that $x\cdot(a_1a_2)=(x\cdot
a_1)\cdot a_2$. A  Banach $A$-bimodule is a \bs\ $X$ which possesses
both a left and right Banach module action such that $(ax)b=a(xb)$ for
all $a,b\in A$ and $x\in X$.

In the case that $A\subseteq\B(X)$ is a subalgebra of $\B(X)$ for some
\bs\ $X$ there is a canonical associated \rep, namely the identity
\rep\ $\id:A\ra \B(X)$. This makes $X$ into a left Banach $A$-module
in a natural way. The algebra $\B(X)$ possesses a natural Banach
$A$-bimodule structure, given by the algebra multiplication.

For the sake of brevity we will refer to left Banach $A$-modules
simply as $A$-modules. Furthermore, we will use module terminology to
refer to the appropriate Banach module concepts. For instance, by a
module map between $A$-modules we mean a \cts\ module map, and by
a submodule we mean a closed submodule.

Suppose that $X$ is a module for a \ba\ $A$.  The set of all
submodules of $X$ forms a lattice, denoted \Index{$\lats X A$}, with
the operations $\Va\wedge\Vb=\Va\cap \Vb$ and
$\Va\vee\Vb=\clos{\Va+\Vb}$.  When there is no danger of confusion we
write simply $\lat A$ for $\lats X A$.  The closure operation in the
definition of the join is necessary, since the sum of two submodules
need not be closed.  For an explicit example, suppose $A=\C$ and
consider the scalar action of $A$ on a \hs\ $\H$. The submodules of
$\H$ are simply the closed subspaces, and it is a standard exercise
to construct two closed subspaces of a \hs\ whose sum is non-closed.

For any set $B$ of operators on $X$ we say that a
subspace $\V\subseteq X$ is $B$-invariant if $B\V\subseteq \V$.
The set of all $B$-invariant subspaces is again a lattice with the
above operations, which we will denote by $\lat B$. 
This notation is consistent with the first definition of $\lat A$,
since if $A\subseteq\B(X)$ is a subalgebra of $\B(X)$ then 
the $A$-invariant subspaces are exactly the set of $A$-submodules of $X$. 
It is readily
verified that if
$B_1\subseteq B_2$ then $\lat B_1\supseteq\lat B_2$. Also, if
$B\subseteq\B(X)$ is a set of operators and $A$ is
the closed algebra generated by $B$ then $\lat A=\lat B$.

\comment{ Needs topologies defined Finally, if $B\subseteq\BH$ then
  $\lat B=\lat \clos{B}$, where the closure can be taken with respect
  to any of the weak, \sw, strong and \sstrong\ topologies.  }

Let $X$ be a module for a \ba\ $A$ with associated \rep\ $\theta:A\ra
\B(X)$. If $\theta$ is one-one we say that $\theta$ is a
\Index{faithful} \rep\ of $A$.  When $E\subseteq X$ we write $AE$ for
the algebraic span of the set $\{ax: a\in A, x\in E\}$.  The module
$X$ is said to be cyclic\index{\rep!cyclic} with cyclic vector $x\in
X$ if $\clos{Ax}=X$, and $X$ is said to be
irreducible\index{\rep!irreducible} if every $x\in X$ is cyclic for
$X$; in this case we also say that $A$ acts
transitively\index{transitive algebra} on $X$.  If
$Ax=X$ for some $x\in X$ we say that $X$ is strictly cyclic with
strictly cyclic vector $x$. When $X$ is strictly cyclic for every $x\in
X$ we say that $X$ is strictly irreducible and that $A$ acts strictly
transitively.

From the definitions it follows that $X$ is irreducible \iff\ $\lat A=\{0,X\}$.
Similarly, $X$ is strictly irreducible \iff\ $X$ possesses no proper invariant
(not necessarily closed) submanifolds.

We say that the \rep\ $\theta$
is nondegenerate\index{\rep!nondegenerate} if $\clos{AX}=X$.
If $A$ has identity and $\theta(1)=1_X$ we say that $\theta$
is a unital \rep. A \rep\ of an algebra with identity is unital \iff\ it is
nondegenerate. 

If $\theta:A\ra\B(X)$ is a \rep\ of $A$ and $\V\in\lat A$, then the
induced \rep\ $A\ra\B(\V)$ is called a sub\rep\index{\rep!subrepresentation}
of $\theta$.

When $A$ is a \csalg\ and $\theta:A\ra \BH$ is a \rep\ on a \hs\ $\H$, we say that
$\theta$ is a \srep\ if $\theta(a^*)=\theta(a)^*$ for all $a\in
A$. A \srep\ is automatically contractive, and a faithful \srep\ is
automatically isometric.

\subsection{The \oa\ $\BH$}

The most important \ba\ for this work is the algebra $\BH$ where $\H$
is a \hs.
Apart from the norm topology, there are several topologies on $\BH$
which we will need to use. They are all defined in terms of various
families of seminorms.
\begin{enumerate}
\item The \Index{weak topology} is defined by the family of seminorms
$$
a\mapsto |(a\xi|\eta)|
$$
where $\xi$ and $\eta$ are any vectors in $\H$.
\item The \Index{\sw\ topology} is defined by the family of seminorms
$$
a\mapsto \bigg|\sum_{i=1}^{\infty} (a\xi_i|\eta_i)\bigg|
$$
where $\{\xi_i\}, \{\eta_i\}\subseteq \H$ with $\left(\sum
  \norm{\xi}^2\right)\cdot\left(\sum \norm{\eta_i}^2\right)< \infty$. 
\item The \Index{strong topology} is defined by the family of seminorms
$$
a\mapsto \norm{a\xi}$$
where $\xi\in \H$.
\item
The \sstrong\ topology is defined by the family of seminorms
$$
a\mapsto \left(\sum_{i=1}^\infty \norm{a\xi_i}^2\right)^{1\over 2}
$$
where $\{\xi_i\}\subseteq \H$ and $\sum\norm{\xi_i}^2<\infty$.
\end{enumerate}
The multiplication on $\BH$ is
separately \cts\ with respect to all of the above topologies.
The following diagram summarises the relations between these
topologies (where $x<y$ means `$x$ is finer than $y$').
$$
\begin{array}{ccccc}
\mbox{\rm norm} & < & \mbox{\rm $\sigma$-strong} & < & \mbox{\rm \sw} \\
        && \wedge & & \wedge \\
        && \mbox{\rm strong} & < & \mbox{\rm weak}
\end{array}
$$

It is a standard result that the \bs\ $\BH$ is dual to the \bs\ 
\Index{$\TC \H$} of trace class operators on $\H$, where the duality
is given by $\<f, x>=\mname{tr}(fx)$. This means that $\BH$ is endowed
with a \ws\ topology. In fact, this \ws\ topology coincides with the
\sw\ topology. 
Consequently, bounded subsets of $\BH$ are precompact in the \sw\ topology.
We denote by \Index{$\KomH$} the space of compact
operators on $\H$. It is another standard result that
$\KomH^*=\TC{\H}$, and so we have
$$
\bd {\KomH}=\TC{\H}^*=\BH
$$
for all \hs{}s $\H$.

We will denote the \sw\ closure of a set $E\subseteq\BH$ by $\swclos{E}$.

When $\H$ is finite-dimensional we can identify $\H$
with $\C^n$ for $n=\mname{dim}\H$, whereupon $\BH$ is identified as
$M_n(\C)$.
We denote by $\{e_{ij}: 1\leq i,j\leq n\}$ the usual set of matrix
units in $M_n(\C)$. The standard basis for $\C^n$ is denoted
$\{e_i\}$.
Since we only ever consider the complex field we will write $M_n$ for
$M_n(\C)$.

For an invertible operator $S\in\BH$ we define the map $\ad S:\BH\ra
\BH$ by $\ad S(T) = \simty S T$\index{$\ad S$}. This is an algebra
automorphism, but it does not preserve the involution on $\BH$ unless
$S$ is a scalar multiple of a unitary operator.  For brevity we write
$\ad S(T)=T^S$, and say that $S$ is a similarity on $\H$.  If $S_1$
and $S_2$ are both invertible then $\ad (S_1S_2)=\ad S_1\cdot \ad
S_2$, so $\rm ad$ defines a homomorphism from ${\rm inv}(\BH)$ to
${\rm aut}(\BH)$. If $A\subseteq\BH$ is an \oa, then $\ad S$ carries
$A$ to an isomorphic \oa, which we denote by \Index{$A^S$}. It is
appropriate to consider $A$ and $A^S$ as very much `the same' \oas.
Certainly they are isomorphic; moreover, $S$ carries the submodule
lattice of $A$ to that of $A^S$ in a good\galley{nice} way: $\lat
A^S=\{S\V:\V\in\lat A\}$.  We say that $A$ and $A^S$ are
spatially isomorphic, or simply that $A$ and $A^S$ are similar.  It is
important to note that the automorphism $\ad S$ is a homeomorphism of
$\BH$ with respect to any of the topologies defined above.

{\def\nnorm#1{\norm{#1}_{{\rm new}}} The process of applying a
  similarity can be profitably thought of as a renorming of $\H$ with
  an equivalent \hs\ norm. If $\nnorm{\cdot}$ is an equivalent Hilbert
  norm on $\H$, then the `identity' operator
  $\iota:\H\ra(\H,\nnorm{\cdot})$ is an isomorphism.  Letting
  $\iota=US$ be the polar decomposition of $\iota$\label{Polar
    decomposition}, the positive part $S$ gives
  $\nnorm{\xi}=\norm{S\xi}$.  Considering $S$ as a similarity,
  applying $\ad S$ is much the same as changing to the new norm
  $\nnorm{\cdot}$. Specifically, if $a\in \BH$ and $\nnorm{a}$ is the
  operator norm of $a$ considered as an operator on $(\H,
  \nnorm{\cdot})$, then $\nnorm{a}=\norm{a^S}$.  } The `size' of a
similarity $S$ is measured by the product $\norm{S}\norm{\inv S}$.

For a set $B\subseteq\BH$, the commutant of $B$ is the
set 
$$
B'=\{T\in \BH: Tb=bT \hbox{ for all $b\in B$}\}.
$$
The commutant of any set is a weakly closed unital \oa.  The commutant
of the commutant of $B$ is called the bicommutant, and is denoted by
$B''$.  As a concession to readability, when $A$ is an \oa\ with \rep\ 
$\theta:A\ra \BH$ we will refer to the commutant
$\theta(A)'\subseteq\BH$ as $A'$ if there is no possibility of
confusion.

If $A\subseteq\BH$ is a \sa\ \oa\ with $A=A''$, then $A$ is a
\vna. This class of algebras enjoys an extremely extensive theory,
which is discussed at length in \cite{DixmierV}.

\subsection{Bounded approximate identities}

A left \Index{\bai} for a \ba\ $A$ is a bounded net $\net{e_\nu}\subseteq A$ such
that $\norm{e_\nu a - a}\ra 0$ for all $a\in A$. Right \bais\ are
defined analogously, and a net which is both a left and right \bai\ is
simply called a \bai. Not all \bas\ possess \bais; however, a \csalg\
always possesses a \bai\ with bound 1.

When $A$ possesses a \bai\ and $X$ is an $A$-module, \Index{Cohen's
factorisation theorem} \cite[Chapter VIII, theorem 32.22]{HewittRoss}
states that the set $\{ax: a\in A, x\in X\}$ is closed and coincides
with $\clos{AX}$. 
Representations of algebras with \bais\ need not be nondegenerate;
however, they differ from nondegenerate \reps\ in an essentially 
trivial way:
\begin{lemma}
  \label{bai lemma}
  Suppose $A$ is a \ba\ with \bai\ and $\theta:A\ra\BH$ is a \rep\ of
  $A$ on a \hs\ $\H$. Then there is a projection $e\in
  \swclos{\theta(A)}$ with $e\H=\clos{A\H}$ and $A(1-e)\H=0$.
  If $M$ bounds the norm of the \bai\ then $e$ can be chosen with
  $\norm{e}\leq M\norm{\theta}$.
\end{lemma}
\proof Let $\{e_\nu\}$ be a \bai\ for $A$, with $\norm{e_\nu}\leq M$
for all $\nu$. By dropping to a subnet, we may assume that the net
$\net{\theta(e_\nu)}$ converges \sw{}ly to an operator $e\in\BH$ with
$\norm{e}\leq M\norm{\theta}$.  By the Cohen factorisation theorem we have
$\clos{A\H}=A\H$. Let $\xi\in\H$ be arbitrary. Since $e_\nu\xi\in A\H$
we have $e\xi\in\clos{A\H}=A\H$. Moreover, since $\net{e_\nu}$ is a
\bai, for $a\xi\in A\H$ we have that $ea\xi=\lim e_\nu a\xi=a\xi$, and
so $e$ is a bounded projection from $\H$ onto $A\H$. Since $\lim
ae_\nu=\lim e_\nu a=a$ we have $e\in A'$, so $\ker e=(1-e)\H$ is
$A$-invariant, and in fact $A(1-e)\H=(1-e)A\H=0$.  \qed

\subsection{Idempotent operators}

Idempotent operators on \hs\ are important in this work.
It is immediate that if $p\in\BH$ satisfies $p^2=p$ then $\H$ is a
linear direct sum $\H=p\H\oplus\ker p$.  Conversely, if $\V_1$
and $\V_2$ are (closed) subspaces of $\H$ with $\H=\V_1\oplus\V_2$ then the open
mapping theorem shows that the projection of $\H$ onto $\V_1$ with
kernel $\V_2$ is \cts.  To emphasise the geometric qualities of
idempotent operators we will often refer to them as projections. A
projection is uniquely specified by its range and kernel. If
$\H=\V_1\oplus\V_2$ and $p$ is the projection of $\H$ onto $\V_1$ with
kernel $\V_2$, we say that $p$ is the {\em projection onto $\V_1$ along
$\V_2$} and that $\V_1$ and $\V_2$ are complements of each other.

The norm of a non-zero projection is $\geq 1$.
A projection $p$ is of norm $\leq 1$ \iff\ its kernel is 
orthogonal to its range. 
This happens exactly when $p$ is a \sa\ operator.
More generally, a simple exercise shows that
the norm of a projection $p$ is cosec of the angle made between
its range and its kernel. Thus the norm of a projection measures how
`close' its range and kernel are.

Even though a projection may not have orthogonal range and kernel,
it is possible to apply a similarity to $\H$ to rectify
this. Specifically, let $p\in\BH$ be a projection, and define a new
norm on $\H$ by 
$$
\norm{\xi}_{\rm new}^2=\norm{p\xi}^2+\norm{(1-p)\xi}^2.
$$ 
This norm is equivalent to $\norm{\cdot}$ and satisfies the
parallelogram identity. Thus $(\H, \norm{\cdot}_{\rm new})$ is a \hs\ 
of the same Hilbert dimension as $\H$, and the polar decomposition
argument given on page~\pageref{Polar decomposition} shows that there
is a similarity $S\in\BH$ with $\norm{\xi}_{\rm new}=\norm{S\xi}$.  It
is readily verified that $\simty S p$ is a contractive projection and
hence has orthogonal range and kernel.

More generally, a result of Dixmier on amenable groups allows us to
extend this to certain families of commuting projections.  A simple
proof of the following may be conveniently found in \cite[theorem
8.3]{Paulsen}. The definition of amenability of groups need not
concern us here; it will suffice for our purposes to note that abelian
groups are always amenable when considered as discrete groups
\cite{Paterson}.
\begin{theorem} \label{Group result}
Let $G$ be an amenable group and $\pi:G\ra \H$ a strongly \cts\
bounded \rep. Then there is a similarity $S\in \BH$ such that
$g\mapsto \pi(g)^S$ is a unitary \rep. If $\norm{\pi(g)}\leq K$ for
all $g\in G$ then $S$ can be chosen so that $\norm{S}\norm{\inv S}\leq
K$.
\qed
\end{theorem}
{
\let\oldDelta\Delta\def\Delta{\mathop{\oldDelta}}
Suppose that $P\subseteq\BH$ is a uniformly bounded set of mutually commuting
projections,
closed under the `symmetric difference' operation $p_1 \Delta p_2 = p_1+p_2-2p_1p_2$.  
 There is an associated set $G$ of invertible operators
given by $G=\{1-2p: p\in P\}$. Since
$(1-2p_1)(1-2p_2)=1-2p_1-2p_2+4p_1p_2=1-2(p_1\Delta p_2)$, $G$ is an
abelian group. We now use 
\nref{Group result} to find a similarity $S$ such that 
$(1-2p)^S$ is unitary for all $p\in P$. Since $1-2p$
has order two, this means that $(1-2p)^S$ is \sa, and so $p^S$ is \sa\
for all $p\in P$.
We will use this observation repeatedly, and so record it in the
form of a \nameref{Orthogonalise idempotents}.
}
\begin{lemma}
  \label{Orthogonalise idempotents}
Let $P\subseteq\BH$ be a uniformly bounded set of commuting
idempotents, closed under symmetric differences. Then
there exists a similarity $S\in\BH$ with $p^S$ \sa\ for all $p\in
P$. If $\norm{p}\leq K$ for all $p\in P$ then $S$ can be chosen with
$\norm{S}\norm{\inv S}\leq 1+2K$.
\qed
\end{lemma}

If $p\in\BH$ is an idempotent then the decomposition $\H=p\H
\oplus \ker p$ induces a $2\times 2$
matrix form on the operators in $\BH$ given by
$$
T\mapsto \TxT pTp, pT(1-p), (1-p)Tp, (1-p)T(1-p);.
$$
The operators  $p$ and $1-p$ have the matrix forms
$$
p=\TxT 1, 0, 0, 0; \mbox{\quad and\quad} 1-p=\TxT 0,0,0,1;.
$$ 
Multiplication and addition in $\BH$ correspond to matrix
multiplication and addition.  Note however that if $p$ is not \sa\ 
then the involution on $\BH$ will not have a nice
\rep\ in terms of this matrix form.
Any other projection $q$ of $\H$ onto $p\H$ has the matrix
form
$$
q=\TxT 1, pq(1-p), 0, 0;=\TxT 1, q_{12}, 0, 0;, 
$$
where $pq(1-p)=q_{12}\in\B(\ker p, p\H)$ can be an arbitrary operator. The
kernel of $q$ is the subspace 
$$
\ker q=\{\xi-q_{12}\xi: \xi\in \ker p\},
$$
which is a a graph over $\ker p$ into $p\H$. We write $\gr q_{12}$ for
this subspace.

In the case when $\H$ is an  $A$-module for an \oa\ $A$, 
projections with invariant kernel and range are particularly
important. It is readily checked that these projections are exactly the
projections which belong to $A'$. Equivalently, they are precisely the
projections which are module maps for $\H$. We refer to such
projections as module projections.

Suppose $p\in A'$ is a module projection and write $\V_1=p\H$, 
$\V_2=\ker p$. For any $T\in\BH$ the decomposition $\H=\V_1\oplus\V_2$ 
induces the matrix form mentioned above:
$$
T=\TxT T_{11}, T_{12}, T_{21}, T_{22};.
$$
The components $T_{ij}$ are maps from $\V_j$ to $\V_i$.  It is readily
verified that $T$ lies in
$A'$ \iff\ each $T_{ij}$ is a module map between $\V_j$ and $\V_i$.
In particular, a projection $q$ of $\H$ onto $\V_1$ will be a module
projection \iff\ $q_{12}=pq(1-p)$ is a module map from $\V_2$ to
$\V_1$. 

\comment{
  An arbitrary projection of $\H$ onto $\V_1$ has the
  matrix form
  $$ q=\TxT 1, q_{12}, 0, 0;,$$ where $q_{12}:\V_2\ra \V_1$ is any
  operator. The kernel of $q$ is the subspace $\{\xi-q_{12}\xi:
  \xi\in\V_2\}$, which is a graph over $\V_2$ into $\V_1$, denoted $\gr
  q_{12}$.  The projection $q$ will be a module projection exactly when
  $q_{12}$ is a module map.
  
  Generalising this, if $T_{21}:\V_1\ra\V_2$ is any module map we may
  consider the module complements of $\gr T_{21}$. These must be graphs over
  $\V_2$ into $\gr T_{21}$, and so are of the form 
  $$
  \left\{\Tvec T_{12}\xi, (1+T_{21}T_{12})\xi; : \xi\in \V_2\right\}
  $$
  where $T_{12}:\V_2\ra\V_1$ is any module map. A simple calculation
  shows that the projection onto $\gr T_{21}$ with this kernel has the
  matrix form
  $$
  \TxT 1+T_{12}T_{21}, -T_{12}, T_{21}+T_{21}T_{12}T_{21},
  -T_{21}T_{12};.
  $$
}

\comment{ In general we will use the symbol $\oplus$ to denote the
  algebraic internal direct sum of two closed subspaces of a \hs. In
  this context the notation does not imply that the summands are
  orthogonal. We also will have occasion to use this symbol for the
  external direct sum of \hs{}s; in this case the direct sum norm is
  understood to be the $\ell^2$ norm, making the direct summands
  orthogonal. It should be clear from the context whether the internal
  non-orthogonal direct sum or the external orthogonal direct sum is
  intended. In doubtful cases we use the symbol $\oplus_2$ to
  explicitly refer to the orthogonal direct sum.  }

In 
\label{Page with complemented subspaces
  theorem}%
Hilbert space every subspace is the range of a projection in $\BH$.
This is characteristic of \hs{}s: the complemented subspaces theorem
of Lindenstrauss and Tzafriri \cite{LindenstraussTzafriri}, \cite{Day}
states that if $X$ is a \bs\ and every subspace of $X$ is complemented,
then $X$ is \isoc\ to a \hs.

\subsection{Tensor products, amplifications and direct sums}
\index{tensor product}

If $\H_1$ and $\H_2$ are \hs{}s there is a natural \hs\ tensor
product, denoted by $\H_1\ot_2\H_2$. This is defined by providing an
inner product on the algebraic tensor product $\H_1\ot \H_2$ and
completing  with respect to the resultant norm. The
inner product of two elements $\alpha_1=\sum \xi_1^i\ot \zeta_1^i$ and $\alpha_2=\sum
\xi_2^j\ot \zeta_2^j$ is 
$$ 
(\alpha_1|\alpha_2)=\sum_{i,j}(\xi_1^i|\xi_2^j)(\zeta_1^i|\zeta_2^j).
$$
It can be verified that this is a positive definite sesquilinear form,
and so the completion of $\H_1\ot\H_2$ with respect to the
corresponding norm is a \hs. Considering $\H_2$ as 
$\ell^2(\Lambda)$ for some index set $\Lambda$, the tensor product can be
identified as 
$$
\H_1\ot_2\ell^2(\Lambda)=\sum_{\lambda\in \Lambda}{}^{\oplus_2}\H_1,
$$
where $\sum^{\oplus_2}_{\lambda\in \Lambda}\H_1$ refers to the \hs\ of
functions $f:\Lambda\ra \H_1$ with 
$$
\norm{f}^2=
\sum_{\lambda\in\Lambda}
\norm{f(\lambda)}^2< \infty.
$$
 In particular, $\H_1\ot_2\ell^2(\N)$ is
canonically identified as the space of square-summable $\H_1$-valued
sequences. If $e_i\in \ell^2(\N)$ is the the vector with $1$ in the
$i^{th}$ place and zero elsewhere, then $\xi\ot e_i\in
\H_1\ot_2\ell^2(\N)$ 
corresponds to the sequence with $\xi$ in the
$i^{th}$ place and zero elsewhere. 

If $T_1\in \B(\H_1)$ and $T_2\in \B(\H_2)$, then the operator
$T_1\ot T_2\in \B(\H_1\ot_2\H_2)$ is defined by putting $T_1\ot T_2 (\xi_1\ot
\xi_2)=T_1\xi_1\ot T_2\xi_2$ and extending linearly. It is readily
verified that $\norm{T_1\ot T_2}=\norm{T_1}\norm{T_2}$.  This
definition may be extended by linearity to give an injective map
$\B(\H_1)\ot\B(\H_2)\ra \B(\H_1\ot_2\H_2)$. The image of this map is
strongly dense in $\B(\H_1\ot_2\H_2)$.

The tensor product notation is also useful for defining finite rank
operators on \hs. If $\xi,\zeta\in\H$ are arbitrary vectors, then the
elementary tensor product $\xi\ot\zeta$ is interpreted as the rank-one
linear operator $\xi\ot\zeta:\eta\mapsto(\eta|\zeta)\xi$. Every finite
rank operator $T\in\BH$ may be written in the form $T=\sum_{i\leq
  n}\xi_i\ot\zeta_i$ for some vectors $\xi_i$ and $\zeta_i$ in $\H$.

Suppose that $A$ is a \ba\ and $\H$ is a module for $A$ with
associated \rep\ $\theta:A\ra \BH$. If $\H'$ is any other \hs\ then
there is a natural module action of $A$ on $\H\ot_2\H'$ given by
$a\cdot (\xi\ot \xi')=(a\xi) \ot \xi'$. In other words, the action on
$\H\ot_2\H'$ is given by the \rep\ $\theta\ot 1_{\H'}:a\mapsto
\theta(a)\ot 1_{\H'}$.  If we consider $\H'$ as $\ell^2(\Lambda)$, the
underlying \hs\ is $ \H\ot_2\ell^2(\Lambda)=\sum_{\lambda\in
  \Lambda}^{\oplus_2}\H$ and the module action is given by the
pointwise action; that is, if $f\in\sum_{\lambda\in
  \Lambda}^{\oplus_2}\H$ then $(a\cdot f)(\lambda)=a(f(\lambda))$.
In the case of $\H'=\ell^2(n)$ for $n\in\N$ the tensor product module
is denoted $\pma\H n$ for brevity, and is referred to as the $n$-fold
amplification of $\H$. In terms of matrices, the matrix corresponding
to $a$ is simply
$$
\left[\begin{array}{ccc}
    \theta(a)&&0\\
    &\ddots&\\
    0&&\theta(a)
  \end{array}
  \right].
$$
We will  refer to the \rep\ $\theta\ot1_n:A\ra\B(\pma\H n)$ 
as the $n$-fold amplification of
$\theta$, denoted by $\pma\theta n$. It is useful to extend this
notation  to the countably infinite case: we
denote the module $\H\ot_2\ell^2(\N)$ by $\pma\H\infty$, and the
corresponding \rep\ by $\pma\theta\infty$. 
Furthermore, in the case when $A\subseteq\BH$ is a subalgebra
of $\BH$ we denote the image of $A$ under $n$-fold amplification
as $\pma A n$. That is, $\pma A n=\{a\ot 1_n:a\in A\}$
for $1\leq n\leq\infty$. Such amplification is isometric
on $A$ for all $1\leq n\leq\infty$,
and a homeomorphism with respect to the \sw\ topologies.
Moreover, if $f\in A^*$ is a \sw\ \cts{} functional
then there are vectors $\xi,\zeta\in\pma\H\infty$ with $\<f,a>=(\pma
a\infty \xi|\zeta)$ for all $a\in A$.

It is simple to identify the commutant and bicommutant of an amplified
\rep. Suppose that $A$ is an \ba\ and $\H$ is a $A$-module with
associated \rep\ $\theta:A\ra \BH$. Let
$1\leq n\leq \infty$, and take an operator $T\in\B(\pma\H n)$. 
The
identification $\pma\H n=\sum^{\oplus_2}_{i\leq n}\H$ induces a matrix
form on $T$ and matrix multiplication shows that $T$ lies in
$\pma\theta n(A)'$ \iff\ each of the matrix entries of $T$ lie in
$\theta(A)'$. 
\comment{
  Thus, we can identify $\pma\theta n(A)'$
  with $\theta(A)'\ot M_n(\C)$ (or the strong closure of $\theta(A)'\ot
  \B(\ell^2)$ when $n=\infty$).
}%
In particular, if $T_0\in\theta(A)'$, then $T_0\ot e_{ij}\in\pma\theta
n(A)'$ for every matrix unit $e_{ij}$.

If $S\in\pma\theta n(A)''$, then $S$ must commute with $T_0\ot e_{ij}$
for every $T_0\in \theta(A)'$ and all $i, j\leq n$.
This implies that $S$ is 
of the form $S=S_0\ot 1_n$ for some $S_0\in A''$. 
Since all such operators commute with $\pma\theta n(A)'$, we have
$\pma\theta n(A)''=\pma{(\theta(A)'')} n$. 

It is a very interesting problem to try to understand $\lat \pma A n$
in terms of $\lat A$.  Even the case $n=2$ is not at all trivial.
However, several classes of submodules of $\pma\H 2$ are immediately
apparent.  If $\V_1, \V_2\in\lat A$, then $\V_1\otimes e_1 +
\V_2\otimes e_2 \in\lat\pma A 2$. At the other extreme, if $T\in A'$
then the graph subspace $\Gr T=\{\xi\otimes e_1+ T\xi\otimes e_2 :
\xi\in \H\}$ is in $\lat \pma A 2$. More generally, suppose that
$Q:\dom Q\ra \H$ is a (possibly unbounded) closed linear operator
defined on a submanifold of $\H$ such that $\dom Q$ is $A$-invariant
and $aQ=Qa$ for all $a\in A$. Then $\Gr Q=\{\xi\otimes e_1 +
Q\xi\otimes e_2:\xi\in\dom Q\}$ is a (closed) submodule of $\lat \pma
A 2$. Such linear operators $Q$ are called \Index{graph
  transformations} for $A$.

For brevity of notation, we will write, for instance, $\V_1\oplus\V_2$
for the subspace $\V_1\otimes_2 e_1 + \V_2\otimes_2 e_2$, whenever 
the meaning is clear from the context.

If $A_1$ and $A_2$ are two \bas, there are several standard ways of
norming the tensor product $A_1\ot A_2$. The only norm we will use is
the projective tensor norm, defined by 
$$
\norm{x}_{\rm proj}=\infof{\sum_{i=1}^n \norm{a_i}\norm{b_i}: 
  x=\sum_{i=1}^n a_i\ot b_i}.
$$
It is readily verified that this defines a norm on $A_1\ot A_2$. The
completion of this normed space is the projective tensor product,
written $\projtensor {A_1}{A_2}$. The algebraic tensor product $A_1\ot
A_2$ has a natural associative algebra product defined by the linear
extension of $(a_1\ot a_2)\cdot(b_1\ot b_2)=(a_1 b_1)\ot (a_2
b_2)$. This multiplication extends to a \cts\ multiplication on
$\projtensor{A_1}{A_2}$, giving the projective tensor product the
structure of a \ba.

\label{Page with B(B(H)) dual topology}
If $X$ and $Y$ are \bs{}s, the dual space $({\projtensor
  {X}{Y}})^*$ is isometrically isomorphic with $\B(X,Y^*)$ under
the map $\alpha:({\projtensor {X_1}{X_2}})^*\ra \B(X,Y^*)$ defined by
$\<\alpha(f)(x),y>=\<f,x \ot y>$. This means that the space
$\B(X,Y^*)$ can be considered as a dual space. Under this
identification a net $\net{T_\nu}\subseteq\B(X,Y^*)$ converges \ws\ to
$T$ \iff\ $\<T_\nu(x)-T(x),y>\ra 0$ for all $x\in X$ and $y\in Y$.

There are several different notions of the direct sum of \bas\ which
we will need to use.  If $\{A_\lambda\}$ is a family of \bas, the most
natural direct sum norm for us is the $\ell^\infty$ norm defined by
$\norm{\lambda\mapsto a_\lambda}=\sup_{\lambda}\{\norm{a_\lambda}\}$.  In
general this
gives rise to two possible direct sums, denoted by
\Index{$\sum^{\ell^\infty} A_\lambda$} and \Index{$\sum^{c_0} A_\lambda$}
respectively. The first is the algebra of bounded functions
$\lambda\mapsto a_\lambda$; the second is the algebra of bounded functions
which vanish at infinity. Of course, there is only a distinction
if $\{A_\lambda\}$ is an infinite family.

\subsection{Completely bounded maps}

Suppose $A\subseteq\B(H')$ is an algebra of operators, and
$\theta:A\ra\BH$ is a \rep. Since $A\ot
M_n$ is embedded in the \csalg\ $\B(\H')\ot M_n=\B(\pma{\H'} n)$, we
may equip $A\ot M_n$ with the corresponding operator norm.  
We define the 
map $\amp\theta n:A\ot M_n\ra
\BH\ot M_n$ by linear extension of
$$
\amp\theta n (a\ot e_{ij})=\theta(a)\ot e_{ij}.
$$
This is a bounded \rep\ of $A\ot M_n$ for all $n\in \N$, and gives
rise to an increasing sequence $\{\norm{\amp\theta n}\}$. In the case
that this sequence is bounded we define the completely bounded
norm\index{completely bounded operator}
$$
\norm{\theta}_{\rm cb}=\sup_{n\in \N}\{\norm{\amp\theta n}\},
$$
and say that $\theta$ is a completely bounded map.
If $S\in\B(\H')$ is a similarity then the homomorphism $\ad
S|_A:a\mapsto a^S$ is \cb\
with $\norm{\ad S|_A}_{\rm cb}\leq \norm{S}\norm{\inv S}$.

\subsection{Alg and Lat}

We have already defined the function $E\mapsto\lat E$ which maps a set
of operators in $\BH$ to a lattice of subspaces of $\H$. Complementary
to this is the function $L\mapsto\alg L$\index{$\alg$}, which maps sets
of subspaces of $\H$ to algebras.  If $L$ is a set of subspaces of
$\H$ we define $\alg L=\{a\in \BH: a\V\subseteq\V \hbox{ for all
  $\V\in L$}\}$. It is easy to verify that $\alg L$ is always a weakly
closed algebra, and that $\alg \lat E\supseteq E$ for all
$E\subseteq\BH$. Similarly, $\lat \alg L\supseteq L$ for all $L$.
Furthermore, $\alg\lat\alg L = \alg L$ and $\lat\alg\lat E=\lat E$ for
all $E$ and $L$. When $A=\alg\lat A$ we say that $A$ is a
\Index{reflexive algebra}; similarly when $L=\lat\alg L$ we say that $L$ is a
\Index{reflexive subspace lattice}.
\label{Page with reflexive} 
A weakly closed unital \oa\ need not be reflexive.  For instance, if
$A$ is the (weakly closed) subalgebra of $M_2$ of matrices of the form
$$
\TxT x,y,0,x;  \;\;\;(x, y\in\C)
$$ 
then $\alg\lat A$ consists of the matrices of the form
$$
\TxT x,y,0,z; \;\;\;(x,y,z\in\C)
$$
and so $A$ is not reflexive.
However, for
infinite amplifications even \sw\ closure is sufficient:
\begin{lemma}
  \label{Reflexive lemma}
  Let $A\subseteq\BH$ be a \sw{}ly closed unital \oa. Then $\pma
  A\infty$ is reflexive.
\end{lemma}
\proof Suppose $T\in\alg\lat \pma A\infty$.  Let $\{\xi_i\}_{i\in\N
  }\subseteq\pma\H\infty$ be a square summable sequence and take
$\epsilon>0$.  Since $\pma{(\pma\H\infty)}\infty\cong\pma\H\infty$ we
may consider the vector
$$
\xi=\xi_1\oplus \xi_2\oplus \cdots \in
\pma\H\infty\oplus\pma\H\infty\oplus\cdots\cong\pma\H\infty.
$$
The
subspace $\clos{\pma A\infty \xi}$ is an $\pma A\infty$-invariant
subspace, containing $\xi$ since $A$ is unital. Thus
$T\xi\in\clos{\pma A\infty \xi}$, and so there is $a\in A$ with
$\norm{(\pma a \infty - T)\xi}< \epsilon$. This implies that
$\norm{(\pma a\infty-T)\xi_i}< \epsilon$ for all $i$, and so $T$
lies in the \sstrong\ closure of $\pma A\infty$; since $A$ is \sw{}ly
closed it follows that $T\in \pma A\infty$ and consequently $\pma
A\infty$ is reflexive.  \qed

\subsection{The dual and bidual of an \oa}

We will consider any \bs\ $X$ to be isometrically contained in its
bidual in the usual way. 
Under this embedding $X$ is \ws\ dense in
$\bd X$, and the elements of $X$ can be characterised in $\bd X$ as
those functionals on $X^*$ which are $\sigma(X^*,X)$-\cts.

The duals and biduals of subspaces and quotients of a \bs\ $X$ are
related pleasantly to $X$. If $Y\subseteq X$ is a subspace, then $Y^*$
may be isometrically identified with the quotient $X^*/Y^\perp$, where
$Y^\perp=\{f\in X^*: \<f, Y>=0\}$. The dual of the quotient $(X/Y)^*$
may be isometrically identified with $Y^\perp$. Thus, $\bd Y$ may be
isometrically identified with $Y^{\perp\perp}$. This latter space is
exactly the \ws\ closure of $Y$ considered as a subspace of $\bd X$.
Similarly the bidual of $X/Y$ is isometrically identified with $\bd
X/\bd Y$.

If we have two \bs{}s $X$ and $Z$, and a map $f:X\ra Z$, then we
denote by $f^*$ the dual map $f^*:Z^*\ra X^*$, and by $\bd f$ the
bidual map $\bd f:\bd X\ra \bd Z$. 
A simple Hahn-Banach exercise shows that
the dual map $f^*$ is one-one \iff\ $f$ has norm dense range. Similarly,
$f^*$ has \ws\ dense range \iff\ $f$ is one-one. 
Moreover, $f^*$ is an
isomorphism \iff\ $f$ is an isomorphism.
In general $\bd{(\ker f)}\subseteq\ker \bd f$, but $\ker\bd f$ need
not equal $\bd{(\ker f)}$. 
\label{Page with closed range}
However, in one important case these two subspaces are equal. Suppose
that $f:X\ra Z$ has closed range. We may replace $Z$ with $f(X)$ and
assume that $f$ is onto.  Let $N=\ker f$, and let $\hat f:X/N\ra Z$ be
the induced (\cts) quotient map. Then the functions ${\hat f}^*:Y^*\ra
N^\perp$ and $f^*:Y^*\ra X^*$ differ only in their codomains---as set
functions their values are the same. Since $\hat f$ is an isomorphism,
so is $\bd{\hat f}:\bd X/\bd N\ra \bd Y$. Let $\bd x\in\bd X$ be an
element with $\bd x\not\in\bd N$. Then there is $y^*\in Y^*$ with
$\<\bd{\hat f}(\bd x + \bd N),y^*>\neq 0$. That is, $0\neq\<\bd x+\bd
N,\hat f^*y^*>=\<\bd x,f^*y^*>=\<\bd f\bd x,y^*>$.  Hence $\bd
x\not\in\ker\bd f$, and so $\ker \bd f=\bd{(\ker f)}$ whenever $f$ has
norm-closed range.

In the case when $X$ is a left module over some \ba\ $A$, there is a
natural right module structure on $X^*$ defined by $\<f\cdot a,
x>=\<f,a\cdot x>$.  Similarly, if $X$ is a right Banach $A$-module, we
define an analogous left module structure on $X^*$.  Finally, if $X$
is an $A$-bimodule there is a corresponding bimodule structure on
$X^*$.  It is readily verified that the `Banach module' nature of the
module action is preserved under these constructions.  An $A$-module
is called a \Index{dual module} if it is obtained from one of these
three processes.

When $A$ is a \ba, it is possible to furnish $\bd A$ with a
multiplication to make it into a \ba\ in such a way that $A$ is a
subalgebra of $\bd A$. In general there are several ways of doing
this, which can lead to different algebra structures on $\bd A$. The
most popular approach is via \Index{Arens products}. There are usually two
distinct Arens products, obtained by considering $A$ either as a left
module over itself or a right module over itself.

If we consider $A$ as a left $A$-module, then the above comments show
that there is a natural left $A$-module structure on $\bd A$. The
action $A\times \bd A\ra \bd A$ may be extended to a multiplication
$\bd A\times \bd A\ra \bd A$ by requiring \ws\ continuity in the left
variable, yielding a \ba\ structure on $\bd A$ extending that of $A$.
This multiplication is the left Arens product on $\bd A$.  The same
construction can be carried out, starting with $A$ as a right
$A$-module; this gives another multiplication on $\bd A$, called the
right Arens product. In the case when the two products coincide the
\ba\ $A$ is said to be Arens regular. For Arens regular algebras, the
multiplication on $\bd A$ is separately \ws\ \cts.

It is a well-known fact that every \csalg\ is Arens regular. In fact,
if $B$ is a \csalg\ there is a faithful \srep{} $\pi:B\ra \BH$ for some
$\H$ such that $\bd B$ is isometrically isomorphic to the \sw\ closure
$\swclos{B}$, and such that the Arens product on $\bd B$ coincides
with the algebra product on $\swclos{B}$ inherited from $\BH$.
Moreover, the \ws\ topology on $\bd B$ corresponds to the relative
\sw\ topology on $\swclos{B}$. This \rep\ is known as the
Gelfand-Naimark-Segal (GNS) \rep\index{GNS \rep} of $B$.  
Consequently, if $A\subseteq\B(\H')$ is a (non-\sa)
algebra of operators on some \hs\ $\H'$ then $A$ is Arens regular
and there is an isometric \rep\ $\pi:A\ra \B(\H)$ such that $\bd A$ is
isometrically isomorphic to $\swclos{A}$. Furthermore the \ws\ topology on $\bd A$
becomes the relative \sw\ topology on $\swclos{A}$, and the Arens
product on $\bd A$ coincides with the algebra product on
$\swclos{A}\subseteq\BH$ \cite{EffrosRuan}. The \rep\ $\pi$ can be
taken to be the GNS \rep\ of any \csalg\ containing $A$ isometrically
($\B(\H')$, for instance).
When we speak of the bidual of an \oa\ as an algebra, it is always
this algebra structure we refer to. Since \oas\ are Arens regular,
multiplication in $\bd A$ is separately \ws\ \cts.  Consequently, if
$A$ is an \oa\ and $J\ideal A$ is an (left/right/two-sided) ideal,
then $\bd J$ is a \ws\ closed ideal of the same type.

The dual of the injection $\iota: \TC\H\ra \BH^*$ is a projection
$\iota^*: \bd
\BH\ra\BH$, which is a \ws--\sw\ \cts\ algebra homomorphism. This
means that if $\theta:A\ra\BH$ is any \rep\ of $A$, then
$\bd\theta:\bd A\ra\bd\BH$ is a \ws-\ws\ \cts\ \rep\ of $\bd A$, and
the map
$\overline\theta:\bd A\ra\BH$ defined by $\bar\theta=\iota^*\bd\theta$ is a \ws--\sw\ \cts\ \rep\ of $\bd A$.
This \rep\ extends $\theta$: \ie, if $a\in A$ then
$\overline\theta(a)=\theta(a)$. Since $A$ is \ws\ dense in $\bd A$,
this extension is uniquely defined by the requirement of
continuity. This is a 
universal property of $\bd A$ familiar from the \csalg\ case.

For \label{Page with BH dual module}%
any \hs\ $\H$, the algebra $\BH$ can be considered as a bimodule
over itself in the natural way. Recall that as a \bs\ $\BH$ is dual to
the space $\TC \H$ of trace-class operators on $\H$. In fact, there is
a natural bimodule structure on $\TC \H$ so that $\BH$ becomes a dual
bimodule.
The module structure on $\TC\H$ is obtained by considering the
embedding $\TC\H\subseteq\BH^*$. The separate \ws\ continuity of the
multiplication in $\BH$ implies that $\TC\H$ is a submodule of
$\BH^*$, and so we may consider $\TC\H$ as a $\BH$-bimodule in its own
right. It is then an immediate consequence of the definitions that the
bimodule $\BH$ is dual to $\TC\H$.

The same construction shows that if $\H$ is a module for a \ba\ $A$,
and $\BH$ is given the corresponding bidmodule structure, then in fact
$\BH$ is a dual module for $A$.

\subsection{Amenability}
\index{amenable algebra}
Amenability is a property of \bas\ introduced by Johnson in
\cite{Johnson1}. The definition of amenability is given in terms of
the Hochschild-Johnson cohomology groups for \bas, which we briefly introduce.
\index{cohomology groups}

{\def\L{{\cal L}}
Let $A$ be a \ba\ and $X$ be an $A$-bimodule. The space of $n$-linear
\cts\ maps from $A\times\cdots\times A$ to $X$ is denoted $\L^n(A,X)$. 
For $n=0$ the space $\L^n(A,X)$ is identified as $X$.
The
coboundary operator $\delta^n:\L^{n-1}(A,X)\ra \L^{n}(A,X)$ is defined
by 
\begin{eqnarray*}
(\delta^n \phi)(a_1,\ldots, a_n) & = & 
  a_1 \phi(a_2, \ldots, a_n) + \sum_{i=1}^{n-1} (-1)^i \phi(a_1,
  \ldots, a_j a_{j+1}, \ldots, a_n) \\
  & + & (-1)^n\phi(a_1,\ldots, a_{n-1})a_n.
\end{eqnarray*}
In the case $n=1$ this is interpreted as 
$$
(\delta^1 x)(a)=ax-xa,
$$
where $x\in X$.
It is readily verified that $\delta^{n+1}\delta^n=0$ for any $n$, so
the complex 
$$0\ra \L^0(A,X)\stackrel{\delta^1}\ra
  \L^1(A,X)\stackrel{\delta^2}{\ra} \L^2(A,X)\ra\cdots$$
gives rise to a collection of cohomology groups
$\cohomology^n(A,X)=\ker\delta^{n+1}/\im\delta^{n}$.

The cocycles in $\ker{\delta^2}$ are the maps $\phi:A\ra X$ with
$a\phi(b)-\phi(ab)+\phi(a)b=0$ for all $a, b\in A$; that is, the maps
$\phi$ with $\phi(ab)=a\phi(b)+\phi(a)b$. Such maps are called
\Index{derivations} of $A$ into $X$.  The coboundaries in
$\im{\delta^1}$ are the derivations $\phi:A\ra X$ of the form
$\phi(a)= ax-xa$ for some $x\in X$. These are called \Index{inner
  derivations}.

A \ba\ $A$ is said to be amenable if $\cohomology 1(A,X^*)=0$ for
every dual Banach $A$-bimodule $X^*$; that is, if every (\cts)
derivation into a dual module is inner.  There are many
characterisations of amenable \bas, one of which is in
terms of so-called \bad{}s\index{\bad}. A bounded net $\net{\sum_i
  a_\nu^i\ot b_\nu^i}\subseteq\projtensor A A$ is called a \bad\ for
$A$ if it satisfies
\begin{enumerate}
\item $\net{\sum_i a_\nu^i b_\nu^i}$ is a \bai\ for $A$.
\item $\bnorm{\sum_i ca_\nu^i\ot b_\nu^i- \sum a_\nu^i\ot b_\nu^i c}\ra
  0$ for all $c\in A$.
\end{enumerate}
A \ba\ $A$ is amenable \iff\ it possesses a \bad\ \cite{Johnson1}.

}

\subsection{Abelian \bas}

We will use a few facts from the theory of abelian \bas\ which warrant
mention here. 
An ideal of a \ba\ is said to be modular if corresponding quotient
algebra is unital.
When $A$ is an abelian \ba\ we denote by \Index{$\p A$} the set of maximal
modular ideals of $A$. Each $\omega\in\p A$ has codimension one and
gives rise (via the quotient map) to a non-zero homomorphism of $A$
into $\C$. Conversely, the kernel of such a homomorphism is always a
maximal modular ideal.

Because of this connection with $\C$-valued homomorphisms,
the set $\p A$ can be considered as a subset of the unit ball
$B_1(A^*)$ of
$A^*$, and so inherits the relative \ws\ topology. This is a locally
compact topology on $\p A$ which is compact \iff\ $A$ is unital. The
evaluation map $\Gamma:A\ra \C^{\p A}$ given by $\Gamma
(a)(\omega)=\<a,\omega>=\omega(a)$ is in fact a contractive
homomorphism $\Gamma:A\ra C_0(\p A)$ of $A$ into the space of continuous
functions vanishing at infinity on $\p A$. This map is called the
\Index{Gelfand transform}, and is often denoted by $\Gamma(a)=\hat a$. 
The radical of $A$ is the kernel of $\Gamma$, and so when $\Gamma$ is one-one
$A$ is \ss.

Apart from the relative \ws\ topology on $\p A$, there is also the
\Index{hull-kernel topology} obtained from the following closure
operation.  If $E\subseteq\p A$, define the ideal $\ker E=\{a: \Gel
a(E)=0\}$. If $J\ideal A$ is an ideal, define $\hul J=\{\omega\in\p A:
J\subseteq\omega\}$. The closure operation on $\p A$ which defines the
hull-kernel topology is $\clos{E}=\mathop{\mbox{hul ker}} E$. It is
readily verified that this defines a (generally non-Hausdorff)
topology which is coarser than the relative \ws\ topology. It is not
generally true that the functions obtained from the Gelfand transform
are \cts\ with respect to the hull-kernel topology.  When the
hull-kernel and \ws\ topologies on $\p A$ coincide $A$ is said to be
regular.

Whenever $X$ is a locally compact Hausdorff space the algebra $C_0(X)$
may be equipped with the \Index{uniform norm} $\norm{f}=\supof{|f(x)|: x\in
  X}$.  This makes $C_0(X)$ into a \ba. The involution
$f^*(x)=\overline{f(x)}$ gives $C_0(X)$ the structure of an abelian
\csalg.\index{abelian \csalg} Conversely, if $A$ is an abelian
\csalg\ the Gelfand transform is an isometric *-\iso\ onto $C_0(\p A)$.
Thus up to *-\iso\ the abelian \csalgs\ are exactly the algebras $C_0(X)$ where $X$
is a locally compact Hausdorff space.  

\subsection{\csalgs\ of compact operators}

\label{Page with compact csalgs}\index{compact operators, algebras of}
Self-adjoint subalgebras of $\KH$ exhibit a very transparent
Wedderburn-type structure.  If $A\subseteq\KH$ is a \sa\ algebra of
compact operators acting nondegenerately on $\H$, then there is a
family $\V_\lambda$ of irreducible submodules of $\H$ and an
integer-valued function $\lambda\mapsto n_\lambda$ (the multiplicity
function) such that $\H$ is unitarily equivalent to the orthogonal
direct sum $\sum^{\oplus_2} \V_\lambda\ot_2 \C^{n_\lambda}$ and $A$ is
unitarily equivalent to the \csalg\ $\sum^{c_0} \K(\V_\lambda)\ot
1_{n_\lambda}$. The central projections of $A''$ correspond to
projections onto subsets of the terms $\K(\V_\lambda)\ot
1_{n_\lambda}$, and the irreducible submodules of $\H$ are the spaces
$\V_\lambda\ot_2 \C\xi$, where $\xi\in \C^{n_\lambda}$
\label{Page with C* details}
(see, e.g.,
\cite{ArvesonBook} for details). 

\subsection{Historical problems concerning \oas}

There are three open problems concerning \oas\ which are relevant to
this work.

The first problem was posed by Kadison in \cite{Kadison}, and concerns
\reps\ of \csalgs\ on \hs{}s which do not preserve the involution.
\begin{question}[Similarity problem]\label{similarity question}
Let $A$ be a \csalg, and let $\theta:A\ra \BH$ be a \rep\ of $A$ on a
\hs\ $\H$. Does there exist a similarity $S\in \BH$ such that
$\theta^S:a\mapsto a^S$ is a \srep?
\end{question}

The similarity question has received much attention in the last 20 years. We
will discuss it further in section~\ref{Section with csalgs}. The most
extensive reference for this problem is Pisier's book \cite{PisierSimilarity}.

\medskip

The next two questions can be considered as non-commutative
generalisations of the invariant subspace problem
from operator theory.

\begin{question}[Reductive algebra problem]
  \label{reductive algebra problem}
  Let $A\subseteq\BH$ be a weakly closed \oa, and suppose that
  $\V\in\lat A\implies\V^\perp\in\lat A$. Does this imply that $A$ is \sa?
\end{question}

A special case of the reductive algebra problem occurs when $\lat
A=\{0,\H\}$, whereupon the condition in the above problem is trivially
satisfied. Since the only weakly closed \sa\ transitive subalgebra of
$\BH$ is $\BH$ itself, the question becomes:
\begin{question}[Transitive algebra problem]\label{transitive algebra problem}
  Let $A\subseteq\BH$ be a weakly closed \oa\ with $\lat
  A=\{0,\H\}$. Does this imply that $A=\BH$?
\end{question}

Both the reductive algebra problem and the transitive algebra problem
have received much attention in the past.  The monograph
\cite{RadjaviRosenthal} gives a survey of what was known in the early
70s on these problems.  Several special cases have been proved, the
most satisfactory of which, due to Lomonosov,  effectively
solves both problems as long as $A$ contains `sufficiently many'
compact operators. Lomonosov's result will be discussed in more
detail in chapter~\ref{Chapter with characterisations}.

\comment{


subreps

comments about vector functionals on GNS rep.

*-representations, definition


In abelian section: mention other approaches after direct integral stuff.

A section on non-tras? Eg L^1(G), etc. Ie use tra properties to rule
out algebras. Observe that being a subalgebra of B(X) relates to bai.

}

\chapter{The \rp}
\label{Chapter with rp}
\label{Chapter with definition}

\dealwithsectionbreaks
\section{Introduction}
\label{Section with triangular algebra}

A central philosophy in mathematics is reductionism: in order to understand
a mathematical object, one tries to decompose the
object into simpler parts. In each area of mathematics, there are
several general questions to be asked: What are the `simpler' objects?
In what sense can the objects of study be decomposed into these? To
what extent does this decomposition help our understanding of the subject?

The theory of group \reps\ provides a good example. A
finite-dimensional complex
\rep\ of a finite group can always be decomposed into a direct sum of
irreducible sub\reps. This strong property permits
the classical theory of group \reps\ to proceed.  

The critical point that makes this decomposition useful 
is the fact that the original \rep\
can be pieced together so simply from its constituent sub\reps. This
quality does not carry over to
associative algebra theory. For instance, let $A\subseteq M_2(\C)$
denote the algebra of matrices of the form
$$ 
\TxT *, *, 0, *;.
$$ 
$\C^2$ is naturally a module for $A$ with submodule $\C e_1$.
There are no other proper submodules, so here
the original \rep\ cannot be reconstructed from its sub\reps.
The problem is that the submodule $\C e_1$ is not
complemented in $\C^2$ by another submodule. 

A topological setting where \rep\ theory is invaluable is the theory of
\csalgs. Here the relevant \reps\ are \sreps---\ie\ homomorphisms
$\pi:A\ra\BH$ of a \csalg\ $A$ which preserve the involution. The
image of $A$ under a \srep\ is \sa, which implies that if $\V\in\H$ is a
submodule then so is $\V^\perp$.  Thus, every submodule for a \srep\ is
complemented by another submodule.  This provides the basis for a
powerful spatial decomposition theory. The topological
nature of \csalgs\ complicates matters, but the fact that submodules
lie inside $\H$ in such a good way facilitates much of the
extensive theory of \csalgs.

\comment{ 
  By way of contrast, such a reduction property is not
  obviously maintained for \reps\ of \csalgs\ which do not preserve
  the involution.  Kadison's \Index{similarity question} asks
  precisely whether this property extends to all \hs\ \reps\ of
  \csalgs\ (see section~\ref{Section with csalgs} for more on this).
  }

These considerations imply that it might be interesting to study \bas\ 
which enjoy a similar \rp\ for their modules.  Since we will be
dealing with \bas, it is appropriate to restrict attention
to Banach modules rather than general algebraic modules.

The most straightforward translation of the above property leads to
the class of \bas\ for which all submodules of arbitrary \bm{}s are
complemented by submodules.  Unfortunately there are no \bas\ with
this property. For, suppose that $A$ is any \ba\ and $X$ is a \bs\ 
with non-complemented subspace $Y$. We may equip $X$ with the
zero module action: $a\cdot x=0$ for all $a\in A$, $x\in X$. Under
this action $Y$ is a submodule which is not complemented in $X$ by
any other submodule, since it is not complemented in the underlying \bs.

A less contrived example is afforded by considering the algebra $c_0$
and the module $\ell^\infty$, where the module action is given by
pointwise multiplication. The subspace $c_0\subseteq\ell^\infty$ is a
submodule which is not complemented in $\ell^\infty$ as a submodule,
since again it is not complemented in the underlying \bs.

Because of this functional-analytic complication, in order to define a
non-trivial class of \bas\ some restriction on the modules we require
to be complemented must be introduced.  There are several approaches
to making this restriction.  The most well known of these gives rise
to the class of algebras with global dimension zero \cite{Helemskii}.
A \ba\ $A$ is said to have global dimension zero if whenever
$$
0\ra X\ra Y\ra Z\ra 0
$$
is a short exact sequence of Banach $A$-modules which splits in the
category of \bs{}s then the sequence splits in the category of \bm{}s.
In other words: if a submodule is complemented in the sense of \bs{}s
then it is complemented in the sense of \bm{}s. This natural
class of \bas\  unfortunately appears to be rather small. Although a complete
characterisation is not presently known, it is conjectured in
\cite{Helemskii} that a \ba\ 
with global dimension zero is of the form
$$
\sum_{i\leq n}{}^\oplus M_{d_i}.
$$
In certain special cases this result is known to be true. As a simple
example, it is elementary to show that the result is true for
commutative $A$; that is, a commutative \ba\ of global dimension zero
is simply $\C^n$ for some~$n$. 

To admit a more diverse class of algebras it is appropriate to
introduce greater restrictions on the class of modules under
consideration. The approach we follow is to restrict attention
to the class of Banach modules for which every subspace is
complemented as a \bs. This certainly removes the topological
obstruction to module complementation. The complemented subspaces
theorem\label{complemented subspaces theorem} 
mentioned on page~\pageref{Page with complemented subspaces theorem}
shows that these Banach modules are exactly those whose underlying
\bs\ is \isoc\ to a \hs.
A similar restriction is found to be efficacious in the classical
theory of group \reps\ when moving from finite groups to compact
groups.

\dealwithsectionbreaks
\section{Definition of the \rp}

We start with a formal definition of an \oa. The definition we give is
phrased to emphasise that we are less concerned with the isometric
structure of an \oa\ than with its \iso\ class as a \ba.
\tracingmacros=1
\begin{defn}
  A \ba\ $A$ is called an \Index{\oa} if there is a 
  \hs\ $\H$ and a \ba\ isomorphism of $A$
  onto a closed subalgebra of $\BH$.
\end{defn}
\tracingmacros=0
It is often convenient to consider `concrete' \oas\ which are given as
subalgebras of $\BH$ for some \hs\ $\H$.  When we say `consider an
\oa\ $A\subseteq\BH$' we always intend the concrete \oa\ 
structure inherited from $\BH$.

Despite the non-isometric emphasis, we will define properties of
concrete \oas\ which depend on the particular realisation as a
subalgebra of $\BH$. For economy of language, we will not notationally
distinguish these from true topological algebra invariants; in the
same spirit we will refer to the norm of a \bs\ which we are really
only interested in up to \iso. It will
always be clear from the context what is happening, so no harm is
done.

\begin{defn}
  Let $A$ be an \oa, and let $X$ be an $A$-module.  If $X$
  is isomorphic to a \hs\ then $X$ is called a \Index{Hilbertian $A$-module}.
\end{defn}
The point of this definition is that again we do not need to consider
the isometric structure of a \hm, but merely its \bs\ structure.
If $X$ is a Hilbertian $A$-module and $X\cong \H$, then
$\H$ is a Hilbertian $A$-module with the induced module action. 
When $A\subseteq\BH$ is a subalgebra of $\BH$ for some \hs\ $\H$, 
there is a natural Hilbertian $A$-module structure on $\H$, which we
will use without further discussion.

\begin{defn}\label{Defn of CP}
\label{Defn of RP}
  Let $A$ be an \oa\ and suppose $\H$ is a \hAm A. $\H$ is said to have the
    \Index{\rp} if for every closed submodule $\Va\subseteq\H$ there is
  another closed submodule $\Vb\subseteq\H$ with $\H=\Va\oplus \Vb$. 
\end{defn}
\begin{defn}
  Let $A\subseteq\BH$ be a concrete \oa. If $\H$ has the \rp\ we say that
  $A$ is a \RA.
\end{defn}

\longpage

The \rp\ is an \oa{}ic version of an old concept from algebra. It can
be interpreted as a kind of semisimplicity (as will be shown in
section~\ref{Section with finite dimension}). The \rp\ has been
discussed in several papers in the literature. An early reference is
Kadison's paper \cite{Kadison}, where the \rp\ is mentioned in
connection with the similarity problem\index{similarity question}.
Later, Fong reintroduced the property, generalised to apply to
algebras of operators on arbitrary \bs{}s, under the name of `complete
reducibility' \cite{Fong}.  There it is shown that if $A\subseteq\BH$
is a \RA\ with the additional property that all its graph
transformations are bounded, then $A$ is of the form
$$ 
A\cong\B(\H_1)\oplus\B(\H_2)\oplus\cdots\oplus\B(\H_k),
$$ 
where $\{\H_i\}$ is a finite collection of \hs{}s.
The strength of this conclusion is a reflection of the strength of the
assumption about the graph transformations, not of the \rp. Indeed, the
prototype \RA\ is any \sa\ \oa---this class of \oas\ is rather wider
than the class Fong investigates.

If $A\subseteq\BH$ is an \RA, it is clear from the definition that
$\lat A$ is a complemented lattice. In fact, the \rp\ implies that
if $\V\in\lat A$ then there is $\V'\in\lat A$ with
$\V\wedge\V'=\{0\}$, $\V\vee \V'=\H$ and
$\V+\V'$ closed. We say that $\lat A$ is topologically complemented
in this case. This is a stronger requirement than lattice-theoretic
complementation. For instance, if $\V_1$ and $\V_2$ are closed
subspaces of $\H$ with $\V_1\cap \V_2=\{0\}$ and $\V_1+\V_2$ dense but 
not closed,
and $A=\alg\{\V_1, \V_2\}$, then simple calculations with rank one
operators show that $A$ is transitive on both $\V_1$ and $\V_2$, and
$\lat A=\{0, \V_1, \V_2, \H\}$. This lattice is complemented but
not topologically complemented, and $A$ does not have the \rp. 
Since we are only interested in topological complementation of
submodules, when we say that two modules are complements we always mean in
the topological sense, not the lattice sense.

\comment{ There is a special case of the \rp\ which has received
  attention in the literature.  A (weakly closed) \oa\ $A\subseteq\BH$
  with the property that $\V\in\lat A$ implies $\V^\perp\in\lat A$ is
  said to be reductive\index{reductive algebra}
  \cite{RadjaviRosenthal}.  All \vnas\ have this property.  The study
  of reductive algebras is extensive and will be discussed further in
  section\galley{Careful}~\ref{Section with reductive algebras}. The
  main open question in this area is the `reductive algebra problem':
  are all reductive algebras \sa?  The requirement that the algebra be
  weakly closed is important here, since it is easy to find
  non-weakly-closed non-\sa\ reductive algebras (there is an example
  on page~\pageref{Page with RP example}).  }

A degenerate case of the \rp\ occurs when $\lat A=\{0, \H\}$. 
Recall from \nref{transitive algebra problem} that it is an open
question whether such an algebra is weakly dense in $\BH$. 
To avoid
necessarily including this rather intractable class of \oas\ in our
study, it is appropriate to strengthen the definition of the \rp\ 
somewhat.

\longpage

\begin{defn}\label{defn of crp}
  Let $A$ be an \oa\ and $\H$ be a \hAm A. We say that $\H$ has the
  \Index{\crp} if the module $\pma \H\infty$ has the \rp. 
  When $A\subseteq\BH$ and $\H$ has the \crp, we say that $A$ is a \Index{\cra}.
\end{defn}
\begin{defn}\label{defn of trp}
  Let $A$ be an \oa. We say that $A$ has the \Index{\trp} if every \hAm A has
  the \rp. For brevity, we will also say that $A$ is a \Index{\tra}.
\end{defn}

\longpage

From the definitions, it is immediate that when $A$ is a \tra\ 
every \hAm A has the \crp\ as well as the \rp.  
We have the following elementary but useful result, which appears in
\cite{Fong}. 
\begin{lemma}
  \label{Submodules have rp}
  Let $A$ be an \oa\ and $\H$ a \hAm A. If $\H$ has the \rp\ then
  every submodule of $\H$ has the \rp.
\end{lemma}
\proof
%
%
Suppose that $A$ is an \oa\ and $\H$ is a \hAm A with the \rp.  Let
$\V\subseteq\H$ be a submodule of $\H$. If $\Vb\subseteq\V$ is a
submodule of $\V$, then $\Vb$ is a submodule of $\H$ and so there is a
submodule $\Vc$ such that $\H=\Vb\oplus\Vc$. The intersection $\Va\cap
\Vc$ is again a closed submodule of $\Va$, and
$\Vb\oplus(\Va\cap\Vc)=\Va$.\qed

Since $\H$ appears as a submodule of $\pma \H\infty$, it follows that
if $\H$ has the \crp\ then it also has the \rp.

The \trp\ for $A$ is the only one of the three properties defined
above which is obviously a \ba{}ic invariant.  All three \rps\ are
stable under spatial isomorphisms of \oas\ (\ie\ under similarities).
Also, every \sa\ \oa\ has the \crp.

The original motivation for studying the \rp\ and its kin came from a
paper of Willis \cite{Willis-operator-normal}, 
\label{Page with Willis}%
in which singly-generated amenable \oas\ were studied.  There it is
shown that if a compact operator $T\in\Kom(\H)$ generates an amenable
\oa\ $A_T$, then $T$ is similar to a normal compact operator.  Since
the spectral projections of a normal compact operator are contained in
the generated \oa, this means that $A_T$ is similar to a \sa\ algebra.
On the other hand, it is known that every abelian \sa\ \oa\ is
amenable \cite{Johnson1}, so Willis' result can be rephrased as: for
$T\in\KH$, $A_T$ is amenable \iff\ $A_T$ is similar to a \csalg.

Away from compact operators, Scheinberg \cite{Helemskii-survey} has a
pretty theorem which can be viewed in a similar light to that of
Willis.  Scheinberg's result says that any amenable uniform algebra is
isomorphic to $C_0(X)$ for some locally compact Hausdorff space $X$.
This can be interpreted as a result about \oas, since uniform algebras
are always \oas.

These results raise the question: what can be said about other classes
of amenable \oas? In fact, an examination of Scheinberg's result shows
that the proof still works if amenability is replaced by the \trp.
This will be discussed in detail in section~\ref{Section abelian}.
Moreover, although it is not immediately obvious, it is also possible
to replace amenability with the \crp\ in Willis' result. Since in
section~\ref{Cohomological section} it is shown that amenability
implies the \trp\ and hence the \crp, the question becomes: what can
be said about \oas\ with either the \crp\ or the \trp?

In section~\ref{section with compact case} we generalise Willis'
result considerably
by showing that if $A\subseteq\KH$ is an arbitrary
algebra of compact operators with the \crp, then $A$ is similar to a
\csalg.

\comment{ 
  These positive results imply that amenability is an
  unnecessarily strong condition for the `spatial' results obtained by
  Willis and Scheinberg. On the other hand, Rosenoer's (and Fong's)
  need for additional hypotheses imply that the bare \rp\ is generally
  too weak on its own to produce good results.  This point of view
  will be borne out in what is to come. The \crp\ and \trp\ seem to be
  good compromises between these extremes.  
  }

Related to this is the work of Rosenoer \cite{Rosenoer1}
\cite{Rosenoer2}, in which \RA{}s which commute with compact operators
are studied. Rosenoer's work stems from \cite{Fong}, and is set
in the context of algebras of operators on a general \bs.
The results obtained there imply that if $T\in\Kom(X)$ is a compact
operator on a \bs\ $X$ such that $A_T$ has the \rp, then under various
additional hypotheses (which amount to something like a
`multiplicity-free' situation) then $T$ is a spectral
operator of scalar type (\ie\ the \bs\ analogue of a normal operator).
When this result is restricted to \hs\ operators, it is 
related to but neither stronger nor weaker than Willis' theorem.
In \cite{RosenoerAnnoying}, Rosenoer extends his previous work
to yield related results, which we will discuss in section~\ref{section with compact case}.

The general thesis of this work is that \oas\ which possess one of the
three \rps\ defined above are `like' \csalgs, in some sense which will
be clarified as we proceed. To provide a framework on which to base
the discussion, we propose the following conjecture as a non-\sa\
analogue of the reductive algebra problem (\nref{reductive algebra problem}).
\begin{conjecture}\label{cra conjecture}
  If $A\subseteq\BH$ is a \sw{}ly closed \cra, then $A$
  is similar to a \sa\ \oa.
\end{conjecture}

The most obvious approach to proving \nref{cra conjecture} is to
firstly try to `orthogonalise' $\lat A$: that is, if $A\subseteq\BH$
is a \cra, try to show that there is a similarity $S$ such that $\lat
A^S$ is {\em orthogonally} complemented. Such a similarity would make
$A^S$ a reductive algebra. If this can be achieved the second step is
to show that $A^S$ is \sa.  In the light of the reductive
algebra problem, it might seem as if the second stage of this program
will present substantial difficulties.  In fact, we will show in \nref{RP
  -> reductive is desirable} that if the first step in this program
can always be completed then \nref{cra conjecture} is true.

\comment{
The reason for requiring \sw\ closure in \nref{cra conjecture} rather
than weak closure is technical:
for \cra{}s the \sw\ topology is more natural, since the
operation of taking
the \sw\ closure commutes with taking the infinite amplification.
In fact we could consider only weakly \cra{}s just as effectively, since
in we show later that
the weak and \sw\ closures of \cra{}s coincide.
}

An example from \cite{Kadison} demonstrates the need for the closure
hypothesis in both the reductive algebra problem and \nref{cra
  conjecture}. Before presenting this example, we make a simple
observation about the stability of the \crp\ under \sw\ closure.
\begin{lemma}\label{cra and weak closure}
  Let $A\subseteq\BH$. Then $A$ has the \crp\ \iff\ the \sw\ closure
  of $A$ has the \crp.
\end{lemma}
\proof
It is readily verified from the definitions that the algebra
$\pma{(\swclos{A})}\infty$ is the weak closure of $\pma A\infty$ in
$\pma\BH\infty$. However, the invariant subspace lattice of an algebra
and its weak closure are the same. 
\qed
For 
the promised example, let $\H=\C^2\ot_2 \ell^2(\N)$ and put
$$
T=\TxT 0, 1, 0, 0;\ot 1_{\ell^2(\N)}\in\BH.
$$ Let $A$\label{Page with RP example}
be the norm-closed algebra generated by $T$ and
$\Kom(\H)$. The \sw\ closure of $A$ is
$\B(\H)$, which, being \sa, is a \cra. This implies that $A$ is a
\cra. However, $A$ cannot be
\isoc\ to a \sa\ \oa\ since the quotient of $A$ by the two-sided
ideal $\K(\H)$ is nilpotent.

Note that in this example $A$ does not have the \trp.  The \rep\ 
$\theta:A/\Kom(\H)\ra M_2(\C)$ given by
$$
[T]\mapsto \TxT 0, 1, 0, 0;
$$ 
induces a \rep\ of $A$. However, $\C^2$ does not have the \rp\ under
this \rep\ since $\C e_1$ is a submodule with no module complement.
We will see later that $\BH$ has the \trp\ (\nref{BH has trp}), and so
this shows that unlike the \crp\ the \trp\ can distinguish a concrete
\oa\ from its \sw\ closure.
This suggests that we should expect a stronger conclusion in \nref{cra
  conjecture} if we replace the \crp\ with the \trp. In fact, all
evidence is consistent with the following conjecture.
\begin{conjecture}\label{tra conjecture}
  If $A$ is a \tra, then $A$ is isomorphic to a \csalg.
\end{conjecture}
The attraction with this reformulation is that it is independent of
any particular \rep\ of the \ba\ $A$.

In the case of a concrete \tra\ $A\subseteq\BH$, the question
arises: if $A$ is \isoc\ to a \csalg, is $A$ necessarily similar to a
\csalg?  The answer is yes, as will be shown in section~\ref{Section about
  similarity problem}.

The converse to \nref{tra conjecture} is also an open question; in
section~\ref{Section with csalgs} it is shown that the converse question is
equivalent to Kadison's similarity problem. Thus if
both the similarity problem and \nref{tra conjecture} have positive
answers, then an \oa\ is \isoc\ to a \csalg\ \iff\ it has the \trp.

\smallskip

The structure of a \csalg\ is reflected in both the algebraic
structure of the underlying *-algebra and the isometric character of
the norm.  It is easy to reconstruct the (unique) \csnorm\ of a
\csalg\ from the *-algebraic structure alone: since the spectrum of an
element of a \ba\ is defined algebraically, the spectral radius of
every element of is definable without reference to the norm. Then the
fact that $\norm{x}^2=\norm{x^*x}=r(x^*x)$ for elements of a \csalg\ 
shows that the norm can be recovered.

On the other hand, it is well-known that if $A$ is a unital \csalg,
then the involution can be reconstructed using only the isometric \ba\ 
data \cite{BonsallDuncan}. Specifically, the states of $A$
may be identified as those norm one functionals $f\in A^*_1$ with
$\<f, 1>=\norm{f}=1$. Then the \sa\ elements of $A$ are identified
as those elements whose image under the states of $A$ is 
real. Specifying the \sa\ elements of $A$ clearly specifies the
involution.

This procedure for defining the `states' and `\sa\ elements' of a
unital \ba\ $A$ can be carried out whether or not $A$ is a \csalg.  In
the context of a general \ba\ the `\sa\ elements' so obtained are
called hermitean elements.  The Vidav-Palmer theorem
\cite{BonsallDuncan} characterises the \csalgs\ amongst unital \bas\ 
in the sense that if the set of hermitean elements spans $A$, then the
algebra is isometrically isomorphic to a \csalg.

In the light of these results, it is natural to wonder if it is
possible to characterise the class of \csalgs\ in terms of their
non-isometric \ba\ data alone. There is an interesting result by
Gardner \cite{Gardner} which is relevant here. Gardner shows that if
$A_1$ and $A_2$ are \csalgs\ which are \isoc\ as \bas, then they are
\isoc\ as \csalgs.  This means that if $A$ is a \ba, then either $A$
is not \isoc\ to any \csalg\ or it is \isoc\ to precisely one \csalg\ 
(up to isometric \iso). Consequently, in some sense the \csalg\
structure of a \csalg\ must be encoded in the \ba\ data.

This \label{Page with non-uniqueness}%
does not mean that the involution of a \csalg\ can be recovered
from knowledge of the \ba\ structure alone.  The algebra $M_2(\C)$
illustrates the difficulty.  Up to isometric \iso\ there is only one
\csalg\ \isoc\ to $M_2(\C)$; however, there are many distinct
involutions under which $M_2(\C)$ admits a \csnorm.  For example, if
$S\in M_2(\C)$ is any similarity then the involution $x\mapsto
\inv S\inv S{}^*x^*S^* S$ is a \csalg\ involution for $M_2(\C)$ under
the norm $\norm{x}_{\rm new}=\norm{\simty S x}$.

There is an analogy here with the \Index{complemented subspaces
  theorem}. A straightforward polar decomposition argument shows if
two \hs{}s\ are \isoc\ as \bs{}s then they are isometrically \isoc,
and so any \hs\ structure on a \bs\ $X$ is determined by the
topological data alone. The complemented subspaces theorem provides a
\bs\ criterion to detect this structure (the existence of complements
to arbitrary subspaces) and a method of `constructing' a \hs\ norm on
$X$ equivalent to the original norm. Of course, this norm is not
uniquely defined, and the `construction' of the \hs\ norm involves
making certain arbitrary choices (see \cite{Day} for details).

Since Gardner's theorem shows that the existence of a
\csalg\ structure is encoded in the \ba\ structure, it is interesting
to try to find ways of detecting this structure.  \Nref{tra
  conjecture} (and its converse, the similarity problem) are an
attempt to do this.

\smallskip

Suppose that $A\subseteq\BH$ is an \oa. In order that $A$ be similar
to a \sa\ algebra it is clearly necessary that there be a constant $K\geq
1$ such that every $\V\in\lat A$ is the range of a module projection
$p\in A'$ with $\norm{p}\leq K$.  The \rp\ implies that every
submodule of $\H$ is the range of a module projection;
however, we will see in section~\ref{rp does not imply crp} an example
of an \RA\ for which no such constant $K$ exists.  Happily, the \crp\ 
does ensure the existence of such a bound.
\begin{proposition}
\label{Existence of \projconst{}s for crp}
Let $A$ be an \oa, and $\H$ a \hAm A with the \crp.
There exists $K\geq 1$ so that
for any submodule $\Va\subseteq \H$ there is a
module projection $p:\H\ra\Va$ with $\norm{p}\leq K$.
\end{proposition}
\proof For every submodule $\Va\subseteq\H$, let $K(\Va)$ denote the
infimum of the norms of the module projections onto $\Va$.  We have
$K(\Va)<\infty$ for all $\Va$.  Suppose that there is a sequence
$\{\V_i\}$ of submodules with $K(\V_i)\ra \infty$.  We may consider
$\V_i$ as embedded into the $i^{th}$ component of $\pma\H\infty$ by
$\xi\mapsto \xi\ot e_i$. With this embedding let $V=\sum^\oplus
\V_i\subseteq \pma \H\infty$. Then $V$ is a submodule of $\pma
\H\infty$, and since $\H$ has the \crp\ there is a complementing
submodule $\U\subseteq\pma \H\infty$. Let $p:\pma\H\infty\ra \V$ be
the corresponding module projection. Let $\H_i$ denote the copy of
$\H$ appearing in the $i^{th}$ coordinate position and let
$q_i:\pma\H\infty\ra\H_i$ denote the $i^{th}$ coordinate projection.
Then $p_i=q_i p|_{\H_i}:\H_i\ra \V_i$ is a module projection onto
$\V_i$ with $\norm{p_i}\leq\norm{p}$. However, by assumption
$\norm{p_i}\geq K(V_i)\ra \infty$. This contradiction establishes the
result.  \qed
\begin{defn}
  \label{\projconst}
  Let $A$ be an \oa\ and $\H$ a \hm\ for $A$ with the \crp. The
  smallest $K>0$ such that every submodule of $\H$ is the range of a
  module projection $p$ with $\norm{p}\leq K$ is called the
  \projconst\ of $\H$. When $A\subseteq\BH$ is a concrete \cra, the
  \Index{\projconst} of $\H$ will also be referred to as the
  \projconst\ of $A$.
\end{defn}

For \tra{}s the idea of \nref{Existence of \projconst{}s for
  crp} can be extended by treating more than one \rep\ at once.
\begin{proposition}
\label{Bound on \projconst{}}
  Let $A$ be an \oa\ with the \trp. 
  There is an increasing function $K:\R^+\ra \R^+$
  such that  if $\theta:A\ra \BH$ is a
  \rep\ of $A$ and $\sv\subseteq \H$ is a submodule
  then there is a module projection $p:\H\ra\sv$ with $\norm{p}\leq
  K(\norm{\theta})$. 
\end{proposition}
\proof 
Take $C>0$. 
Suppose that there is a sequence $\{\theta_i:A\ra \B(\H_i)\}$
of \reps\ with $\norm{\theta_i}\leq C$ and a sequence
$\{\V_i\subseteq\H_i\}$ of submodules such that $K(\V_i)\ra \infty$.
Consider the direct sum \rep\ $\theta:A\ra \B(\sum^\oplus\H_i)$ given
by $\theta(a)(\xi_i)=(\theta_i(a)\xi_i)$. Then $\norm{\theta}\leq C$,
and since $A$ is a \tra\ the module $\H=\sum^\oplus\H_i$ has the \crp,
and there is a projection $p\in \theta(A)'$ onto
$\V=\sum^\oplus\V_i$. As before, if we denote by $q_i$ the projection
from $\H$ onto $\H_i$, then $p_i=q_ip|_{\H_i}$ is a projection in
$\theta_i(A)'$ onto $\V_i$, with $\norm{p_i}\leq \norm{p}$ for all
$i$. This contradiction implies the existence of the function $K$.
\qed

It is worth observing that, notwithstanding these results, it is
unusual for there to be a uniform bound on the set of {\em all\/}
module projections. In fact, such a bound would imply that each
submodule was uniquely complemented in $\lat A$, since if $p$ and $q$
are two module projections with the same range, then $p+\lambda(q-p)$
is also a module projection for any $\lambda\in\C$. In the case of a
\sa\ \oa\ the phenomenon of unique complementation in $\lat A$
corresponds to the `multiplicity-free' situation \cite{ArvesonBook}.

\medskip

In defining the \rp\ we have chosen to restrict attention to \hm{}s
rather than general Banach modules in order to avoid the issue of
non-complemented subspaces of \bs{}s. There is another, more general
approach, introduced by Fong in \cite{Fong}, where a `completely
reducible algebra' is defined to be a closed subalgebra $A\subseteq\B(X)$ for
some \bs\ $X$ such that $\lat_X A$ is topologically complemented (\ie, for all
$\V\in\lat_X A$ there is $\W\in\lat_X A$ with $X=\V\oplus\W$).  Many
of our results can be applied without great modification to such
algebras. This setting has not been adopted here for several reasons.
One reason is that the greater generality obtained from the \bs\ 
setting does not seem to shed any new light on the subject. Much of
the motivation for studying the \rp\ and its kin arise from
connections with established areas and problems in non-\sa\ operator
theory, all of which are found in the \hs\ setting. Secondly, while
the \rp\ is easy to generalise to a \bs\ setting, the 
discussion of the previous section\galley{Careful} illustrates that
that \trp\ is not. Since the \trp\ is the most significant of the
three reduction properties discussed, there is nothing to gain by admitting Fong's
completely reducible algebras here.

\dealwithsectionbreaks
\section{A cohomological definition of the \trp}
\label{Cohomological section}

The definition of the \trp\ can be recast into a 
cohomological setting. The cohomological definition is less
illuminating to work with, but has the advantage that it displays the
connection between the \trp\ and other notions already in the
literature.

If $\theta:A\ra \BH$ is a \rep\ of an \oa\ $A$, the
space $\BH$ becomes an $A$-bimodule in the natural way via the algebra
structure of $\BH$. Then we may speak
of derivations from $A$ into $\BH$, and of the cohomology group
$\cohomology 1(A, \BH)$.

\goodbreak
 
\begin{theorem}
  \label{Cohomological definition}
  An \oa\ $A$ has the \trp\ \iff\ $\cohomology 1(A,\BH)=0$ for every
  \rep\ $\theta:A\ra \BH$.
\end{theorem}
\proof
Suppose that $A$ has the \trp\ and that $\theta:A\ra \BH$ is a \rep. Let
$\delta:A\ra \BH$ be a derivation with respect to $\theta$.

Consider the map $\hat \theta:A\ra \B(\H\oplus\H)$ given by 
$$
\hat\theta:a\mapsto \TxT \theta(a), \delta(a), 0, \theta(a);.
$$
This is a \rep\ of $A$, making $\H\oplus\H$ an $A$-module.
The subspace 
$\H\oplus 0$ is a submodule of $\H\oplus\H$, and since $A$ has the \trp\
there exists a complementary submodule $\V$. This subspace must be a
closed graph over $0\oplus \H$ and hence is of the form $\V=\{T\eta\oplus
\eta:\eta\in \H\}=\gr T$ for some $T\in \BH$.

Applying the matrix $\hat\theta(a)$ to $T\eta\oplus\eta$ gives
$$
\TxT \theta(a), \delta(a), 0, \theta(a); \Tvec T\eta,\eta;=
\Tvec (\theta(a)T+\delta(a))\eta, \theta(a)\eta; =
\Tvec T\theta(a)\eta, \theta(a)\eta;
$$ by the invariance of $\V$. Thus $\delta(a)=T\theta(a)-\theta(a)T$
for all $a\in A$, showing that $\delta$ is inner and $\cohomology
1(A,\BH)=0$.

Note that $T$ can be chosen with $\norm{T}\leq
K(\norm{\theta}+\norm{\delta})$, where $K$ is the function of
\nref{Bound on \projconst{}}, since the  module projection $p_T$ associated
with the complement $\gr T$  is 
$$
p_T=\TxT 1, -T, 0, 0;,$$
and 
by suitable choice of $p_T$
we have $\norm{T}\leq \norm{p_T}\leq K(\norm{\hat \theta})\leq
K(\norm{\theta}+\norm{\delta})$. 

Conversely, suppose $\cohomology 1(A,\BH)=0$ for all modules $\BH$ induced by a
\rep\ of $A$ into $\BH$. Suppose that $\theta:A\ra \BH$ is a \rep\ of
$A$ and $\V\subseteq\H$
is a submodule of $\H$.
The orthogonal decomposition $\H=\V\oplus \V^\perp$
gives the matrix form
$$
\theta(a)=\TxT a_{11}, a_{12}, 0, a_{22}; 
$$
for elements of $A$. 
Define
$\hat\theta:A\ra\BH$ and
$\delta:A\ra \BH$ by 
$$
\hat\theta(a)= \TxT a_{11},0 ,0 , a_{22}; \hbox{\quad and\quad}
\delta(a)=\TxT 0,a_{12},0,0;.
$$

If $a, b\in A$ we have
$$
\TxT a_{11}, a_{12}, 0, a_{22}; \TxT b_{11}, b_{12}, 0, b_{22}; =
\TxT a_{11}b_{11}, a_{11}b_{12}+a_{12}b_{22}, 0, a_{22}b_{22}; $$
which shows that  $\hat\theta$ is a \rep\ of $A$ and that
\begin{eqnarray*}
  \delta(ab) & = & \TxT 0, a_{11}b_{12}+a_{12}b_{22}, 0, 0; \\
  & = & \TxT a_{11}, 0, 0, a_{22}; \TxT 0, b_{12}, 0, 0; +
  \TxT 0, a_{12}, 0, 0; \TxT b_{11}, 0, 0, b_{22}; \\
  & = & \hat\theta(a)\delta(b) + \delta(a)\hat\theta(b).
\end{eqnarray*}
Thus $\delta$ is a derivation with respect to the \rep\ 
$\hat\theta$.  
Now considering $\B(H)$ as an $A$-bimodule via $\hat\theta$, the fact
that $\cohomology 1(A,\BH)=0$ implies that 
$\delta(a)=\hat\theta(a)T-T\hat\theta(a)$ for some $T\in\BH$. Writing
$$
T=\TxT T_{11}, T_{12}, T_{21}, T_{22};
$$
and expanding the identity $\delta(a)=\hat\theta(a)T-T\hat\theta(a)$
in components gives $a_{12}=a_{11}T_{12}-T_{12}a_{22}$.  Thus the
subspace
$$
\left\{\Tvec -T_{12}\eta, \eta; : \eta\in\V^\perp\right\}$$ is a
$\theta(A)$-invariant complement to $\V$, and hence $A$ is a \tra.
\qed

This characterisation of \tra{}s should be compared to the definition
of amenability of \bas.  Recall that a \ba\ $A$ is amenable if
$\cohomology 1(A,X^*)=0$ for all \DBxBM As $X^*$.  From the discussion
on page~\pageref{Page with BH dual module} the bimodules $\BH$ which
arise above are in fact dual bimodules.  Consequently, if $A$ is an
amenable \oa\ and $\theta:A\ra\BH$ is a \rep\ then $\cohomology
1(A,\BH)=0$.  This gives the following \nameref{Amenable implies trp}.
\begin{proposition}
  \label{Amenable implies trp}
Let $A$ be an amenable\index{amenable algebra} \oa. Then $A$ has the \trp.\qed
\end{proposition}

A corollary of this is that every abelian \csalg\ has the \trp, since
abelian \csalgs\ are amenable \cite{Johnson1}.

There is a more direct proof that amenability implies the \trp,
which is interesting in its own right. Recall that a \ba\ $A$ is
amenable \iff\ it possesses a \bad. For any nondegenerate \rep\
$\theta:A\ra \BH$ of an amenable \oa\ $A$ the \bad\ allows us to
define a map from invariant subspaces of $\H$ to module projections,
as the next \nameref{Quasiexpectation} shows.

\begin{proposition} \label{Quasiexpectation}
Let $A$ be an amenable \oa, and let $\theta:A\ra \BH$ be a
nondegenerate \rep.
Then there is an $A'$-bimodule
projection $\Xi: \BH\ra A'$, which has the following property: if
$\V\subseteq \H$ is a submodule of $\H$, and $p:\H\ra \V$ is any projection
onto $\V$, then $\Xi(p)$ is  a module projection onto $\V$.
\end{proposition}
\proof Suppose that $\net{\sum_i a^i_\nu\otimes b^i_\nu}$ is a bounded
approximate diagonal for $A$ with
$\sum_i\norm{a^i_\nu}\norm{b^i_\nu}\le M$ for all $\nu$.  The product
net $\net{\sum a_\nu^ib_\nu^i}$ forms a \bai\ for $A$. From
nondegeneracy and \Index{Cohen's factorisation theorem} it follows that
$\H=\clos{A\H}=\{a\xi: a\in A, \xi\in\H\}$, so that every element $\xi\in\H$ may be
factorised as $\xi=a\zeta$ for some $a\in A$ and $\zeta\in\H$.
Consequently $\lim\sum a_\nu^ib_\nu^i\xi=\xi$ for all $\xi\in\H$.  We
consider the family of maps $ \Xi_\nu: \BH \ra \BH $ defined by
$$
\Xi_\nu(T) = \sum_i a^i_\nu T b^i_\nu.$$
Then $\norm{\Xi_\nu} \le M$ for all $\nu$. 

Since $\BH$ is dual to $\TC\H$, $\B(\BH)$ may be identified with the
dual space
$(\projtensor{\BH}{\TC\H})^*$.
Closed balls in $\B(\BH)$ are compact in the resulting
\ws\ topology, and so the net $\net{\Xi_\nu}$ has a cluster point $\Xi$
with $\norm{\Xi}\le M$.  Taking a subnet we may suppose $\Xi_\nu \ra
\Xi$.  We observe that $\Xi$ enjoys several good properties.  
If $c\in A$,
$T\in \BH$, and $\xi, \zeta \in \H$, then
\begin{eqnarray*}
(c\Xi(T)\xi|\zeta) & = & \lim_\nu \sum_i( ca_\nu^iTb_\nu^i\xi|\zeta) \\
        &=& \lim_\nu \<\sum ca_\nu^i\ot b_\nu^i\,, \,
              a\ot b\mapsto (aTb \xi|\zeta)> \\
        &=& \lim_\nu \<\sum a_\nu^i\ot b_\nu^i c\,, \,
              a\ot b\mapsto (aTb \xi|\zeta)> \\
        &=& \lim_\nu \sum_i( a_\nu^iTb_\nu^i c\xi|\zeta) \\
        &=& (\Xi(T)c \xi|\zeta),
\end{eqnarray*}
so $\Xi(T)\in A'$. Also, if $T\in A'$ we have
\begin{eqnarray*}
(\Xi(T)\xi|\zeta) &=& \lim_\nu \sum_i( a_\nu^i T b_\nu^i\xi|\zeta) \\
        & = & \lim_\nu \sum_i( T a_\nu^i b_\nu^i \xi|\zeta) \\
        &=& (T\xi|\zeta).
\end{eqnarray*}
Finally, if $r, s\in A'$ and $T\in \BH$ we have
\begin{eqnarray*}
(\Xi(rTs)\xi|\zeta) &=& \lim_\nu \sum_i( a^i_\nu rTs b^i_\nu \xi | \zeta) \\
        &=& \lim_\nu \sum_i( ra^i_\nu T b^i_\nu s \xi | \zeta) \\
        &=& (r\Xi(T)s\xi|\zeta).
\end{eqnarray*}
In short, $\Xi$ is an  $A'$-bimodule projection from $\BH$ onto
$A'$.

To verify the module projection property, suppose that $\V\subseteq \H$ is a
submodule, let $p:\H\ra \V$ be any projection onto $\V$, and consider
$\Xi(p)$. Let $\xi\in \H$ be an arbitrary vector. Then
$\Xi(p)\xi=\lim_\nu \sum_i a^i_\nu p b^i_\nu \xi$, where the limit is
in the \ws\ topology on $\H$. Since $pb^i_\nu\xi\in\V$, we have
$a^i_\nu pb^i_\nu\xi\in\V$, so $\sum_i a^i_\nu pb^i_\nu\xi\in\V$ and
finally $\Xi(p)\xi\in\V$ since $\V$ is \ws\ closed in $\H$.  Suppose
now that $\xi\in \V$. Then $\Xi(p)\xi=\lim_\nu \sum_i a^i_\nu p
b^i_\nu \xi = \lim_\nu \sum_i a^i_\nu b^i_\nu \xi = \xi$ since
$b^i_\nu\xi\in \V$. That is, $\Xi(p)$ is a projection onto $\V$.  
It is a module projection since it belongs to $A'$.
\qed
\comment{
  Although the above result works only for nondegenerate \reps, the
  existence of a \bai\ in $A$ allows us to reduce to the nondegenerate
  case, as the next result shows.
\begin{lemma}\label{Nondegeneracy basic lemma}
  Let $A\subseteq\BH$ be an \oa. If $\H$ has the \rp\ then there is a
  submodule $\K\subseteq\H$ with $\clos{A\H}\oplus\K=\H$ and
  $A\K=\{0\}$.  If $A$ has a \bai\ and $\clos{A\H}$ has the \rp, then
  $\H$ has the \rp.
\end{lemma}
\proof
Suppose that $\H$ has the \rp.
Let $\K$ be any module complement to $\clos{A\H}$. Then 
$A\K\subseteq \K\cap \clos{A\H}=\{0\}$.

Suppose now that $A$ has a \bai\ $\{e_\nu\}$ and $\clos{A\H}$ has the
\rp.  It is standard that $\{e_\nu\}$ clusters in the \sw\ topology to
a module projection $e:\H\ra \clos{A\H}$ (see for instance,
\cite{Johnson1}). Then $A(1-e)\H=\{0\}$. Moreover, suppose
$\V\subseteq\H$ is $A$-invariant. Since $e\in\swclos{A}$ we have
$e\V\subseteq\V$ and $(1-e)\V\subseteq\V$. On the other hand, since
$e$ is a module projection, $e\V$ and $(1-e)\V$ are $A$-invariant
subspaces of $\V$ with direct sum $\V$. Moreover, $A(1-e)\V=\{0\}$,
and $e\V\subseteq \clos{A\H}$. Since $\clos{A\H}$ has the \rp, there
is a submodule $\W$ of $\clos{A\H}$ with $e\V\oplus\W=\clos{A\H}$.
Then the submodule $\W\oplus((1-e)\H \ominus (1-e)\V)$ is a module
complement to $\V$, demonstrating that $\H$ has the \rp.  \qed }

This gives an improvement of \nref{Amenable implies trp}:
\begin{theorem}
\label{AOAs have the CP}
Let $A$ be an \aoa\index{amenable algebra}. Then $A$ has the \cp.
Moreover, if $\theta:A\ra \BH$ is a \rep\ of $A$ then the \projconst{}
of $\H$ is bounded by $(M\norm{\theta})^2$, where $M$ bounds the norm of
the \bad\ of $A$.\galley{Check this estimate sometime}
\end{theorem}
\proof Suppose that $A$ is an \aoa, with \bad\ $\net{\sum_i a_\nu^i\ot
  b_\nu^i}$. Let $e$ be the projection from \nref{bai lemma}.  Suppose
that $\V\subseteq\H$ is $A$-invariant, and hence invariant for the
\sw\ closure of $A$. The projection $e$ is in the \sw\ closure of $A$,
so $e\V\subseteq\V$ and $(1-e)\V\subseteq\V$.  This shows that
$\V=e\V\oplus(1-e)\V$ is an $A$-module decomposition of $\V$.
Restricting the action of $\theta$ to $e\H=A\H=\clos{A\H}$ yields a
nondegenerate \rep, and $e\V$ is an invariant subspace of $e\H$. The
amenability of $A$ allows us to use \nref{Quasiexpectation} to
construct the map $\Xi$ for the \hm\ $e\H$.  Let $q$ be the orthogonal
projection from $e\H$ onto $e\V$.  Then $\Xi(q)$ is a projection,
commuting with $A$ and with range $e\V$. Let $\w$ be the orthogonal
complement to $(1-e)\V$ in $(1-e)\H$. The sum $\ker \Xi(q)\oplus \w$
is clearly closed, and is a module complement to $\V$.  We have
$\norm{\Xi}\leq M\norm{\theta}$ and $\norm{e}\leq M\norm{\theta}$, so
the projection onto $\V$ along $\ker \Xi(q)\oplus \W$ has norm at most
$(M\norm{\theta})^2$.  \qed

Derivations are very closely related to homomorphisms, and in
\cite{Gourdeau} a characterisation of amenability in terms of
homomorphisms is given. 
Other connections between cohomology and homomorphisms are given in
\cite{JohnsonPerturbation} and \cite{RaeburnPerturbation}.
It is not surprising that the cohomological
definition of the \tra\ can also be stated in terms of properties of
homomorphisms. If $\phi, \psi:A\ra \BH$ are two \reps\ of an \oa\ $A$,
we say that $\phi$ and $\psi$ are similar if there is a similarity
$S\in\BH$ with $\psi=\simty S \phi$. The following result
characterises the \trp\ by the fact that sufficiently close
homomorphisms are always similar.
\begin{lemma}\label{Fred's result}
  Let $A$ be an \oa. Then the following two properties are equivalent:
  \begin{enumerate}
  \item $A$ has the \trp.
  \item For every \rep\ $\phi:A\ra \BH$ there is $\epsilon>0$ such
    that the following holds. 

    If $\psi:A\ra \BH$ is any \rep\ of $A$ with
    $\norm{\phi-\psi}<\epsilon$ then there exists an invertible $S\in
    \BH$ with $\phi^S=\psi$. If $\V\subseteq\H$ is a common
    invariant subspace for $\phi$ and $\psi$ and
    $\phi|_\V=\psi|_\V$ then $S$ can be chosen with $S|_\V=\id_\V$.
  \end{enumerate}
\end{lemma}
\proof
Suppose that $A$ is a \tra, and that $\phi, \psi:A\ra\BH$ are two
\reps. 
Observe from the proof of \nref{Cohomological definition} that
if $\theta:A\ra \BH$ is a \rep\ of $A$ and $\delta:A\ra \BH$ is a derivation
with respect to $\theta$, then there is $T\in\BH$ with
$\delta(a)=\theta(a)T-T\theta(a)$ and $\norm{T}\leq
K(\norm{\theta}+\norm{\delta})$.

Writing $\lambda=\norm{\phi-\psi}$, we define the \rep\ 
$\theta=\phi\oplus\psi: A\ra \B(\pma\H 2)$, and a derivation $\delta: A\ra
\B(\pma\H 2)$ by
$$
\delta(a)=\TxT 0, \inv\lambda(\phi-\psi)(a), 
    0, 0;.$$
Then $\norm{\delta}\leq 1$, and by the \trp\ there is
$T\in\B(\pma\H 2)$ with $\delta(a)=\theta(a)T-T\theta(a)$ and
$\norm{T}\leq K(\norm{\phi\oplus\psi}+1)$. The
operator $T$ has the matrix form 
$$
T=\TxT T_{11}, T_{12}, T_{21}, T_{22};, 
$$
and when the identity $\delta(a)=\theta(a)T-T\theta(a)$ is expanded
into components we obtain the identity $(\phi-\psi)(a)=\phi(a)\lambda
T_{12}-\lambda T_{12}\psi(a)$. Consequently $ \phi(a)(1-\lambda
T_{12})=(1-\lambda T_{12})\psi(a)$, and since $\norm{\lambda
  T_{12}}\leq \lambda \norm{T}$ the operator $1-\lambda T_{12}$ is
invertible if $\lambda< \inv{[K(\norm{\phi\oplus\psi}+1)]}$. This
shows that, for instance,  a suitable choice for $\epsilon$ would be $\epsilon<
\inv{[K(\norm{\phi}+2)]}$.

We now consider the situation where $\V\subseteq\H$ is an invariant
subspace and $\phi|_\V=\psi|_\V$. By the \trp\ there is a
$\phi(A)$-invariant subspace $\W_\phi$ with $\H=\V\oplus\W_\phi$ and
similarly there is a $\psi(A)$-invariant subspace $\W_\psi$ with
$\H=\V\oplus\W_\psi$. The homomorphisms $\phi$ and $\psi$ induce
quotient homomorphisms $\hat\phi:A\ra \B(\H/V)$ and $\hat\psi:A\ra
\B(\H/V)$, with $\norm{\hat \phi-\hat \psi}\leq\norm{\phi- \psi}$.
The $\phi$-induced modules $\H/\V$ and $\W_\phi$ are similar, and the
$\psi$-induced modules $\H/\V$ and $\W_\psi$ are similar.  Thus, by
the above argument applied to the quotient \reps, if
$\norm{\phi-\psi}<\inv{[K(\norm{\phi\oplus\psi}+1]}$ then the
$A$-modules $\W_\phi$ and $\W_\psi$ are similar.  This similarity can
be extended (uniquely) to a similarity of $\H$ which fixes $\V$,
demonstrating that (i)$\implies$(ii).

For the direction (ii)$\implies$(i), 
suppose that $A$ satisfies condition (ii).
Let $\theta:A\ra \BH$ be a \rep\ and
$\delta:A\ra \BH$ a derivation with respect to $\theta$. For $\lambda\in
\C$ we define
\rep{}s $\phi$ and $\psi$ by 
$$
\phi=\TxT \theta, 0, 0, \theta; \hbox{\quad and \quad}
\psi_\lambda=\TxT \theta, \lambda \delta, 0, \theta;.
$$
By assumption, for sufficiently small $\lambda$ there is $S\in\BH$
with $\phi^S=\psi_\lambda$. Moreover, since $\H\oplus0$ is a common
invariant subspace for $\phi$ and $\psi_\lambda$, we may assume that
$S$ fixes $\H\oplus0$. 
Writing 
$$
T=\inv S = \TxT 1, T_{12}, 0, T_{22};,
$$
the matrix component expansion of the equation $\simty S
\phi=\phi^S=\psi$ gives $\delta(a)=\theta(a) T_{12}-T_{12}\theta(a)$
for all $a\in A$. By \nref{Cohomological definition} $A$ is a \tra.
\qed

If $\H$ denotes a fixed \hs, the set of \reps\ from $A$ into $\BH$ is
a topological space with the norm topology. The above result shows
that for \tras\ the similarity orbits in this space are both open and
closed.

\dealwithsectionbreaks
\section{\csalgs\ and the \rp}
\label{Section about similarity problem}
\label{Section with csalgs}

We have observed that if $A$ is a \csalg\ and $\pi:A\ra \BH$
is a \srep\ then $\H$ has the \crp\ as an $A$-module. For non-* \reps\
the situation is not so simple. 
We say that a \csalg\ $A$ has the {\em similarity
property} if every \rep\ of $A$ is similar to a \srep.
With this definition the similarity question (\nref{similarity
  question}) becomes: do all \csalgs\ have the similarity property?
It is elementary to prove that any \csalg\ with the similarity property
also has the \trp.
\begin{lemma}
  \label{Similarity property implies trp}
Let $A$ be a \csalg\ with the similarity property.
Then $A$ has the \trp.
\end{lemma}
\proof
Suppose that $\theta:A\ra\BH$ is a
\rep, and $\V\subseteq\H$ is a submodule of $\H$. There is a
similarity $S\in\BH$ such that $\theta^S$ is a \srep. Then $S\V$
is invariant under $\theta^S(A)$, so $(S\V)^\perp$ is also
invariant. Reversing the similarity, $\inv S((S\V)^\perp)$ is
invariant under $\theta(A)$, and complements $\V$.
\qed

\comment{
An intriguing open question in the theory of \csalgs\ is the {\em
  similarity question} \cite{PisierSimilarity} \cite{Kadison}:
\begin{question}[Similarity Question]\label{Similarity
    question}\index{similarity question}
Does every \csalg\ have the similarity property?
\end{question}

This question has attracted much attention in the past, and although
it remains open in general, the general feeling seems to be that the
answer is `yes'.
}

There are several partial results for the similarity question, most of
which can be found in Pisier's book \cite{PisierSimilarity}. 
Many of these find restrictions on either the algebra or the module
which force a \rep\ to be similar to a \srep. 
The most prominent result of the latter type
 is Haagerup's result for finitely generated 
modules \cite{HaagerupCyclic}, \cite{PisierSimilarity}:
\begin{theorem}\label{Haagerup}
  Let $A$ be a \csalg\ and $\theta:A\ra \BH$ a \rep\ of $A$. If there
  is a finite set $\{\xi_i\}$ of vectors in $\H$ such that the
  submodule generated by $\{\xi_i\}$ is $\H$, then $\theta$ is 
  similar to a \srep.\qed
\end{theorem}

In general, it is not known whether the size of the similarity
obtained in this result can be bounded in terms of $\norm{\theta}$
alone. Indeed, solving the similarity problem is equivalent to
establishing such a bound \cite{PisierSimilarity}.  However, when $A$
has the \trp, the \projconst\ function obtained in
\nref{Bound on \projconst{}} can be used to bound the norms of the
similarities.

The following proof is adapted
from the proof of the \Index{complemented subspaces theorem} which appears in
\cite{Day}. 

\begin{lemma}\label{Uniform bound on similarities}
  Let $A$ be a \csalg\ with the \trp, and let $\theta:A\ra \BH$ be a
  \rep\ which is similar to a \srep. If $K$ is the projection constant
  function of \nref{Bound on \projconst{}} then there is a similarity
  $S$ so that $\theta^S$ is a \srep\ and $\norm{S}\norm{\inv S}\leq
  128 K(\norm{\theta})^2$.
\end{lemma}
\proof
Let us say that a \hm\ $\H'$ is a $*$-module for $A$ if the
corresponding \rep\ is a \srep. Define
$$
\alpha=\infof{\norm{S}\norm{\inv S}: \hbox{$S$ is a module \iso\ from $\H$
     onto a $*$-module $\H'$}}
$$ 
By assumption $\alpha<\infty$ and we may find a $*$-module $\H'$ and a
contractive module \iso\ $S:\H\ra \H'$ with $\norm{\inv S}\leq
2\alpha$, say.  Since $\H'$ is a $*$-module the corresponding \rep\ 
$\theta':A\ra\B(\H')$ is contractive, and the \rep\ 
$\theta\oplus\theta':A\ra\BH\oplus\B(\H')$ has
$\norm{\theta\oplus\theta'}=\norm{\theta}$.  Thus $\H\oplus\H'$ has
the \rp\ with \projconst\ $M\leq K(\norm{\theta})$ and consequently
for any $\mu\in\R^+$ there is a module projection $p$ from
$\H\oplus\H'$ onto $\gr \mu S$ with $\norm{p}\leq M$.  Writing $p$ in
matrix components reveals that $p$ has the form
$$
p=\TxT 1+R \mu S, -R, \mu S(1+R \mu S), -\mu S R;
$$
for some module map $R:\H'\ra \H$. The fact that $\norm{p}\leq M$
implies that $\norm{\mu S(1+R \mu S)}\leq M$ and $\norm{R}\leq M$.

We consider now the module map $T:\H\ra \H'\oplus\H'$ given by 
$$
T\xi={1\over 2}S \xi \oplus {1\over 2M} (\mu S(1+R \mu S)\xi).
$$
Since $S$ is bounded below, $T$ is also bounded below and is
a contractive module isomorphism onto some closed submodule of
$\H'\oplus\H'$. Since $\H'$ is a $*$-module for $A$, $\H'\oplus\H'$ is
also a $*$-module and any submodule of $\H'\oplus\H'$ is a $*$-module. By the
definition of $\alpha$, this means that there is $\xi_0\in\H$ with
$\norm{\xi_0}=1$ and $\norm{T\xi_0}\leq 2\inv \alpha$. 

Suppose that $\norm{S \xi_0}\leq \inv{(2M\mu)}$. From the second
term in the definition of $T$ we see that $\norm{T\xi_0}\geq
\mu/8M\alpha$. Since we know that $2\inv \alpha\geq \norm{T\xi_0}$ this
is impossible if we choose $\mu>16M$. 

Thus if we choose $\mu=16M+\epsilon$, we must have $\norm{S\xi_0}>
\inv{(2M\mu)}=\inv{(32M^2+2M\epsilon)}$. But the first term in the
definition of $T$ then gives us the inequality
$$
2\inv \alpha \geq \inv 2\norm{S\xi_0}>\inv{(64M^2 + 4M\epsilon)},
$$
and hence $\alpha\leq 128M^2$.
\qed

This allows us to pass from finitely generated modules to general
modules, proving that the \trp\ is equivalent to the similarity
property for \csalgs.

\begin{proposition}\label{ultrafilter argument}
  Let $A$ be a \csalg\ with the \trp. Then every \rep\ of $A$ is
  similar to a \srep.
\end{proposition}
\proof Let $\theta:A\ra \BH$ be a \rep\ of $A$.  From Haagerup's
result we know that every finitely generated sub\rep\ of $\theta$ is
similar to a \srep. To pass to the full \rep\ consider the set
$\Lambda$ of finite subsets of $\H$, and let ${\cal F}$
denote an ultrafilter on $\Lambda$ which contains the filterbase of sets
of the form $\{\lambda: \lambda\supseteq \lambda_0\}$. For each
$\lambda=\{\xi_i\}\in \Lambda$ let $\H_\lambda$ be the
submodule of $\H$ generated by $\{\xi_i\}$. Then 
by \nref{Uniform bound on similarities}
there are
contractive similarities $S_\lambda\in\B(\H_\lambda)$ with
$\norm{\inv{S_\lambda}}\leq 128[K(\norm{\theta})]^2$ such that applying
$S_\lambda$ makes $\theta|_{\H_\lambda}$ a \srep. In particular, if
$\xi,\eta\in \lambda\in \Lambda$, then $(S_\lambda
\theta(a)\xi|S_\lambda \eta)=(S_\lambda \xi| S_\lambda
\theta(a^*)\eta)$ for all $a\in A$. We define a new inner product on
$\H$ by
$$
(\xi|\eta)_{\rm new}=\lim_{\cal F} (S_\lambda \xi|S_\lambda \eta).
$$
Routine verification shows that this defines an inner product norm on
$\H$ equivalent to the original norm, and $\theta$ becomes a \srep\ 
with respect to this inner product.  \qed

\begin{corollary}\label{Similarity property equivalent to trp}
  The \trp\ is equivalent to the similarity property for \csalgs.
  \qed
\end{corollary}

\begin{corollary}
  \label{Amenable implies sp}
  Any amenable \csalg\ has the similarity property. 
\end{corollary}
\proof
\Nref{Amenable implies trp} shows that an amenable \oa\ has the \trp. 
\qed

This result is known (see, for example, \cite{PisierSimilarity})
but the proofs in the literature use some relatively sophisticated
machinery concerning the biduals of amenable \csalgs.

In \cite{HaagerupCyclic} Haagerup uses \nref{Haagerup} to show that
any \csalg\ with no tracial states has the similarity property. This
can be used to exhibit examples of \tras\ which are not amenable.
\begin{corollary}\label{BH has trp}
  If $A$ is a \csalg\ with no tracial states, then $A$ has the \trp.
  For any infinite dimensional \hs\ $\H$ the \csalgs\ $\BH$ and
  $\sum^{\ell^\infty}_{i\in\N} \BH$
  both have the \trp.
\end{corollary}
\proof 
Since $\H$ is infinite dimensional there are partial isometries $u,
v\in\BH$ with $uu^*+vv^*=1$. If $\tau\in\BH^*$ is tracial then
$\tau(1)=\tau(uu^*+vv^*)=\tau(u^*u+v^*v)=\tau(2)$, so $\tau(1)=0$. But
if $\tau$ is a state on $\BH$ then $\tau(1)=1$. Thus $\BH$ has no
tracial states. The algebra $\sum^{\ell^\infty}\BH$ contains $\BH$ via
the constant functions; this is a unital embedding. Since any tracial
state on $\sum^{\ell^\infty}\BH$ would restrict to a tracial state on
$\BH$, it follows that $\sum^{\ell^\infty}\BH$ has no tracial states.
\qed

For infinite dimensional $\H$, the algebra $\BH$ does not have the
approximation property and so is not amenable \cite{LauLoyWillis}. 

\comment{ 
A problem related to the similarity question refers to derivations
rather than homomorphisms \cite{PisierSimilarity}.

  Let $A$ be a \csalg. We say that $A$ has the derivation property if,
  for $\cohomology 1 (A,{\BH})=0$ for any \srep\ $\pi:A\ra \BH$.
\begin{question}[Derivation Question]\label{Derivation question}
  Does every \csalg\ have the derivation property?
\end{question}
For general \csalgs\ the derivation question is also open.
\comment{
Suppose that $A$ is a \csalg\ with the \trp. Let $\pi:A\ra \BH$ be a
\srep\ of $A$, and $\delta:A\ra\BH$ a derivation relative to
$\pi$. Then by 
\nref{Cohomological definition}, $\delta$ is inner.
Thus, for \csalgs, the \trp\ implies the derivation property.
}
\Nref{Cohomological definition} shows that the \trp\ implies the
derivation property for \csalgs.

A less obvious argument due to Kirchberg completes the
circle of implications for \csalgs:
\begin{theorem}\label{Similarity property equivalent to trp}
  Any \csalg\ with the derivation property also has the
  similarity property. The similarity property, the derivation
  property and the \trp\ are all equivalent for \csalgs. 
\end{theorem}
\proof
See Theorem 7.21 of \cite{PisierSimilarity} for a proof of the first
part. The second part follows from this and the above results.
\qed
}

If $A\subseteq\BH$ is a \tra, the statement of \nref{tra conjecture}
superficially suggests that $A$ is only algebraically \isoc\ to a
\csalg. In fact it would be equivalent to demand spatial \iso:
\begin{corollary}\label{result about conjectures}
  Let $A\subseteq\BH$ be a \tra\ which is isomorphic to a \csalg. Then
  $A$ is similar to a \csalg.
\end{corollary}
\proof Suppose $\theta:A\ra B$ is an isomorphism from $A$ to a \csalg\ 
$B$. Since the \trp\ is stable under isomorphism, $B$ is a \csalg\ 
with the \trp. The inverse \iso\ $\inv\theta:B\ra A\subseteq\BH$ is a
\rep\ of $B$; consequently there is $S\in \BH$ with
$(\inv\theta)^S:B\ra \BH$ a \srep. The image of $(\inv\theta)^S$ is
$A^S$, and the image of a \csalg\ under a \srep\ is \sa.  \qed

\Nref{Similarity property equivalent to trp}
means that it would be of great interest to know which
\csalgs\ have the \trp. 
It also means that the \trp\ `extends' the
similarity property from the class of \csalgs\ to the class of general
\oas. 

There are other ways of extending the similarity property to non-\sa\ 
\oas. A remarkable result of Haagerup shows
that the similarity property is closely related to complete
boundedness of \reps.
\begin{theorem}\label{cb reps are similar to sreps}
  Let $A$ be a \csalg. A \rep\ $\theta:A\ra\BH$ is
  similar to a \srep\ \iff\ $\theta$ is \cb.
\end{theorem}
\proof
See, for instance, \cite[Theorem 8.1]{Paulsen} for a proof.
\qed

This suggests the introduction of another property of (possibly
non-\sa) \oas:
\begin{defn}
  Let $A$ be an \oa. We say that $A$ has the \Index{\cbrep\ property} if
  every \rep\ of $A$ is \cb.
\end{defn}

The \cbrep\ property is
explicitly introduced in \cite{PisierDegree} as an extension of the
similarity property to non-\sa\ \oas.  
Although the \cbrep\ property and the \trp\ coincide for \sa\ \oas,
the two are not the same 
for non-\sa\ \oas. We shall see in the next
section that there are finite-dimensional algebras with the \cbrep\
property but not the \trp. On the other hand, if \nref{tra conjecture}
is true, then every \tra\ has the \cbrep\ property.

\dealwithsectionbreaks
\section{Examples of \RA{}s}

In this section we examine some examples of algebras with the \rp\ or
one of the stronger variants. The first example demonstrates the \rp\
as a kind of semisimplicity.

\subsection{Finite-dimensional algebras}
\label{Section with finite dimension}\index{Finite-dimensional algebras}%

\begin{proposition}\label{Finite-dimensional plus RA implies ss}
Let $A\subseteq M_n$ be a \RA. Then $A$ is semisimple and 
isomorphic to a direct sum of full matrix algebras.
That is, there is a finite sequence of integers $\{k_i\}_{i=1}^m$ such
that $A$ is isomorphic to $\sum^\oplus M_{k_i}$
\end{proposition}
\proof Since $A\subseteq M_n$, we may consider $\C^n$ as an
$A$-module.  The finite-dimensionality of $\C^n$ and the \rp\ imply
that $\C^n$ may be decomposed into a finite direct sum of irreducible
submodules. The corresponding set of sub\reps\ give a faithful family of
irreducible \reps\ of $A$, so the intersection of the corresponding
primitive ideals is $\{0\}$. Thus $A$ is a semisimple algebra.
An application of the
Wedderburn structure theorem for semisimple associative algebras and the
Gelfand--Mazur\ theorem shows that any
finite-dimensional semisimple (complex) \ba\ is isomorphic to a finite direct
sum of full matrix algebras \cite{Palmer}. 
\qed

Clearly such an algebra can be endowed with a \csalg\
structure. Conversely, since \csalgs\ are semisimple, all
finite-dimensional \csalgs\ are of this form. 

The \rep\ theory for these algebras is well-known. However, it is
illuminating to revise the details here, since the finite-dimensional
case provides the motivation for dealing with the case of compact
operators in section~\ref{section with compact case}.
It also gives a direct demonstration that such algebras enjoy the \trp.

Suppose then that $A=\sum^\oplus_{i\leq n} M_{k_i}$ is a
finite-dimensional \csalg. $A$ possesses $n$ minimal central
projections, $\{p_i\}$ say, and identity $1=\sum p_i$. Suppose that
$\theta:A\ra \BH$ is a \rep. The idempotents $q_i=\theta(p_i)$ form a
commuting disjoint family of projections, with sum $\sum
q_i=\theta(1)$. We add to this family the idempotent
$q_0=1-\theta(1)$. Then if $\V\subseteq\H$ is a submodule, we have
$q_i\V\subseteq\V$, and so $\V$ may be written as $\V=\sum^\oplus
\V_i$, where $\V_i\subseteq q_i\H$ are all submodules. Conversely, if
we have a family $\{\V_i\}$ of submodules with $\V_i\subseteq q_i\H$
for all $i$, then $\sum \V_i$ is a submodule. This observation implies
that $\H$ has the \rp\ \iff\ each of $q_i\H$ has the \rp.

Since $\theta(A)q_0\H=\{0\}$ 
the submodule $q_0\H$ trivially has the \rp. In all
other cases the \rep\ of $A$ induced on the submodule $q_i\H$ has
kernel $\sum_{j\neq i}^\oplus M_{k_j}$, and so is effectively a \rep\ of
$M_{k_j}$. Thus, we need only understand the \reps\ of $M_{k_j}$ in
order to understand \reps\ of $A$. Moreover, we may assume $\H=q_j\H$;
that is, we need only consider unital \reps\ of $M_{k_j}$. 

Consider then a fixed $M_k$ and unital \rep\ $\theta:M_k\ra\BH$.
We write $e_{ij}$ for the usual matrix
units of $M_k$. The operators $\theta(e_{ii})$ are a commuting
disjoint family of
idempotents, and the operators $\theta(e_{ij})$ are partial
isomorphisms, mapping $\theta(e_{jj})\H$ onto $\theta(e_{ii})\H$. We
may define a new inner product on $\H$ by 
$$
(\xi|\zeta)_{\rm new}=\sum_i(e_{1i}\xi|e_{1i}\zeta).
$$
This is readily checked to define a positive definite inner product on
$\H$, producing a norm equivalent to the original norm. It has
the effect of making all the $\theta(e_{ii})$ become \sa\ projections,
and all the $\theta(e_{ij})$ partial isometries. 
The usual polar decomposition argument shows that there
exists a similarity $S$ on $\H$ with
$(\xi|\zeta)_{\rm new}=(S\xi|S\zeta)$.
By conjugating $\theta$ by $S$ we may assume that $\theta(e_{ii})$ is
\sa\ for all $i$ and $\theta(e_{ij})$ is a partial isometry for all
$i\neq j$. 

If we now denote by $\H'$ the subspace $\theta(e_{11})\H$ and by
$1_{\H'}$ the identity operator on $\H'$, then $\H$ is isometrically
identified with $\C^k\ot_2\H'$, and the \rep\ $\theta$ becomes
$\id_{M_k}\ot 1_{\H'}$ under this identification. Moreover, the
submodules of $\C^k\ot_2 \H'$ are exactly those of the form $\C^k\ot_2 \V$,
where $\V\subseteq\H'$ is an arbitrary closed subspace. This implies
that $\H$ has the \rp, since $\C^k\ot_2 \V$ is complemented by many
submodules, but in particular it is complemented by $\C^k\ot_2
\V^\perp$.

From this discussion we get the following result.
\begin{proposition}\label{finite dimensional result}
  Let $A\subseteq\BH$ be a \fd\ \oa. Then the following are equivalent:
  \begin{enumerate}
  \item $A$ is a \tra.
  \item $A$ is a \RA.
  \item $A$ is isomorphic to a finite direct sum of full matrix
    algebras.
  \item $A$ is isomorphic to a \csalg.
  \end{enumerate}
\end{proposition}
\proof The only implication remaining to be proved is that if
$A\subseteq\BH$ is a finite-dimensional \RA\ then $A$ is isomorphic to
a finite direct sum of full matrix algebras. If there is a
finite-dimensional submodule $\V\in\lat A$ such that the restriction
of $A$ to $\V$ is faithful, then \nref{Finite-dimensional plus RA
  implies ss} gives us the desired result.  Suppose then that there is
no faithful finite-dimensional submodule in $\lat A$. Let $\V$ be a
finite-dimensional submodule such that the kernel of the
sub\rep\ of $A$ on $\V$ has minimal dimension. Choose $a\in A$ such that
$a|_{\V}=0$ but $a\neq 0$, and find a vector $\xi\in\H$ with $a\xi\neq
0$. Then the submodule generated by $\V$ and $\xi$ is
finite-dimensional and induces a smaller kernel. This contradiction
implies that there are faithful, finite-dimensional submodules, and
the result is proved.  \qed

\subsection{CSL algebras}
\index{CSL algebras}

A lattice $L$ of subspaces of a \hs\ $\H$ is said to
be a commutative subspace lattice (a CSL for short) if the orthogonal
projections onto elements of $L$ form a commutative family of projections.
An \oa\ $A$ is said to be a CSL algebra if it is of the form $\alg L$
for some CSL $L$ \cite{Davidson}. 

A well-known result in the field of CSL algebras states that
commutative subspace lattices are reflexive in the sense defined on
page~\pageref{Page with reflexive} \cite[Corollary 22.11]{Davidson}.
The next result shows that the class of CSL algebras introduces
no new \RA{}s.
\begin{proposition}\label{CSL RAs are sa}
  Let $A$ be a CSL algebra. If $A$ is a \RA\ then $A$ is \sa.
\end{proposition}
\proof
Since commutative subspace lattices are reflexive, $\lat A$ is a
CSL. Suppose that $A$ is a \RA, and take $\V\in\lat A$. Then there
is a complementing submodule $\W\in\lat A$, and if $p, q$ denote
the orthogonal projections onto $\V$ and $\W$ respectively we
have that $p$ commutes with $q$. The spectral theorem then implies
that $p=1-q$, and so $\W=\V^\perp$. This means that $L$ is
orthogonally complemented, and since $A=\alg L$, it follows that $A$ is
\sa.
\qed

Finite-dimensional CSL algebras are particularly pleasant algebras,
which do not, however, generally have the \trp.  These algebras are
discussed in \cite{Davidson}. We use them here to distinguish between
the \trp\ and the \cbrep\ property discussed in section~\ref{Section
  with csalgs}. The following simple calculation demonstrates that any
finite-dimensional CSL algebra has the \cbrep\ property.

It is shown in \cite{Davidson} that if $A$ is a finite-dimensional CSL
algebra then there is an integer $n$ and a finite set $E$ of matrix
units of $M_n$ such that $A$ is \isoc\ to the subalgebra of $M_n$
generated by $\{e_{ij}: e_{ij}\in E\}$.  The set $E$ can be assumed to
contain $e_{ii}$ for all $i$ and to be closed under the natural groupoid
operation; that is, $e_{ij}, e_{jk}\in E$ implies $e_{ik}\in
E$. The set $E$ can be conveniently thought of as a transitive
directed graph on $n$ nodes, where an edge between nodes $i$ and $j$
indicates $e_{ij}\in E$. This directed graph is then a complete invariant
for the algebra $A$ amongst the class of finite-dimensional CSL
algebras. For this reason such algebras 
are often called digraph algebras.

\comment{
Suppose that
$A\subseteq M_n$ is a finite-dimensional CSL algebra, 
and let $\{e_{ij}\}$ be the usual set of matrix units for $M_n$. 
Without loss of generality
we may assume that the orthogonal projections onto the
invariant subspaces for $A$ are diagonal. Then the
projections $e_{ii}$ lie in $\alg \lat A=A$, for each $i\leq
n$. Thus every element $a\in A$ can be decomposed as $a=\sum
e_{ii}ae_{jj}=\sum a_{ij}e_{ij}$, where $a_{ij}\in\C$, and the structure of 
$A$ is given by the collection of matrix units
$e_{ij}$ which lie inside $A$. This collection is frequently
represented by a directed graph on $n$ nodes, with an edge from $j$ to
$i$ representing the presence of $e_{ij}$ in $A$. The only constraint
that a directed graph need satisfy for it
to arise in this way is that it should be
reflexive and transitive. It is proved in \cite{Davidson} that this graph is a
complete invariant for $A$ amongst the \fd\ CSL algebras. For this
reason, such algebras are often called digraph algebras.
}

It follows from 
\nref{Finite-dimensional plus RA implies ss} that if $A$ is a
digraph algebra with the \rp, the corresponding digraph must be
symmetric. However, all digraph algebras have the \cbrep\ property.
\begin{lemma}
  Let $A$ be a digraph algebra, and $\theta:A\ra \BH$ a \rep\ of
  $A$. Then $\theta$ is \cb.
\end{lemma}
\proof We assume that $A$ is a subalgebra of $M_n$ of the form
discussed above, with $E$ denoting the set $\{e_{ij}: e_{ij}\in A\}$.
Since $\{e_{ii}\}\subseteq A$ is a commuting family of idempotents we
may use \nref{Group result} to apply a similarity on $\H$ to make each
$\theta(e_{ii})$ \sa.  Since $A$ is unital, we may without loss of
generality reduce to the case where $\sum\theta(e_{ii})=\H$.  Having
done this, $\H$ decomposes into an orthogonal direct sum
$\H=\sum^\oplus \H_i$, where $\H_i=\theta(e_{ii})\H$. We write
$T_{ij}$ for the map $\theta(e_{ij}):\H_j\ra \H_i$.

Every
element of $a\in A$ can be written as $a=\sum_{e_{ij}\in E}
a_{ij}e_{ij}$ where $a_{ij}\in \C$, 
whereupon the image of $a$ under $\theta$ is
$\theta(a)=\sum_{e_{ij}\in E} a_{ij}T_{ij}$. Moreover, every element
$b\in A\ot M_m$ can be written as $b=\sum_{e_{ij}\in E} e_{ij}\ot b_{ij}$,
where $b_{ij}\in M_m$. Consequently, the image of $b$ under $\amp
\theta m$ can be written as
$$
\amp \theta m(b)=\sum_{e_{ij}\in E} T_{ij}\ot b_{ij}.
$$
A crude estimate will suffice to show complete boundedness. Observe
that $\norm{b_{ij}}=\norm{e_{ij}\ot b_{ij}}=\norm{e_{ii}\ot 1_m\cdot
  b\cdot e_{jj}\ot 1_m}\leq\norm{b}$. 
Furthermore, $\norm{\amp \theta m(e_{ij}\ot b_{ij})}
=\norm{T_{ij}\ot b_{ij}}=\norm{T_{ij}}\norm{b_{ij}}$. If we put
$M=\max \{\norm{T_{ij}}\}$, then 
$$
\norm{\amp \theta m(b)} \leq \sum M\norm{b_{ij}}\leq Mn^2\norm{b}.
$$
Thus $\norm{\theta}_{cb}\leq Mn^2$, and $\theta$ is \cb.
\qed
Since there are certainly digraph algebras which are not \tra{}s, the
cb-\rep\ property is distinct from the \trp.
For a specific example, the directed graph $\bullet \ra \bullet$ gives the
algebra from section~\ref{Section with triangular algebra} of matrices of the form
$$
\TxT *, *, 0, *;;
$$
this has the \cbrep\ property but not the \trp.

\subsection{Some commutative operator algebras}

The algebra $A_\lambda$ of matrices of the form
$$
\TxT a, \lambda(b-a), 0, b;
$$
is \isoc\ to $\C^2$, and so has the \rp. The submodule $\C e_1\subseteq\C^2$ is
uniquely complemented by $\{(\lambda \eta, \eta) : \eta\in \C\}$, and
so the projection constant for $A_\lambda$ is no less than $\lambda$. 
Exploiting this
simple observation we 
can demonstrate that several other
natural \oas\ do not have the \trp.

If $S$ is a metric space, the \Index{Lipschitz algebra} $\Lip_\alpha(S)$ is the algebra
of bounded \cts\ functions $f:S\ra \C$ such that 
$$
p_\alpha(f)=\supof{{|f(x)-f(y)|\over d(x,y)^\alpha}: x, y\in S}< \infty,
$$
with pointwise algebra operations and the norm
$$
\norm{f}=\supof{|f(x)|: x\in S}+p_\alpha(f).
$$ 

$\Lip_\alpha(S)$ is in fact an \oa, as can be seen by considering the
\reps\ $\theta_{x,y}:\Lip_\alpha(S)\ra M_2$ given by
$$
\theta_{x,y}(f)=\TxT f(x), {d(x,y)}^{-\alpha}(f(y)-f(x)), 0, f(y);.
$$
These are all contractive \reps\ of $\Lip_\alpha(S)$, and if we
define $\theta=\sum^\oplus \theta_{x,y}$ then $\theta$ is a \cts\
faithful \rep. On the other hand, $\theta$ is bounded below and so
the \oa\ $\theta(\Lip_\alpha(S))$ is isomorphic to $\Lip_\alpha(S)$.

The \reps\ $\theta_{x,y}$ also show that $\Lip_\alpha(S)$ is not a
\tra\ when $S$ is not uniformly discrete. For, suppose that
$\Lip_\alpha(S)$ is a \tra; by \nref{Bound on \projconst{}} there is
$K>0$ so that if $\phi:\Lip_\alpha(S)\ra \BH$ is a contractive \rep\
then the \projconst{} for $\H$ is less than $K$. On the other hand,
the only proper invariant subspaces for $\theta_{x,y}$ are the
mutually complementing subspaces
$$
\left\{\Tvec \xi, 0; \right\} \hbox{\quad and \quad} \left\{\Tvec
    {d(x,y)}^{-\alpha}\eta, \eta;\right\}.
$$
If $S$ is not uniformly discrete, then there will be $x, y\in S$ with
$d(x,y)$ small enough 
that the norm of the projection induced by these two complements is
larger than $K$. This contradiction implies that $\Lip_\alpha(S)$ is
not a \tra.
When $S$ is uniformly discrete $\Lip_\alpha(S)$ is \isoc\ to
$\ell^\infty(S)$, which is a \tra\ by \nref{Amenable implies trp}. 

If $A$ is an \oa\ similar arguments can be applied to the algebra
$\Lip_\alpha(S; A)$ of $A$-valued Lipschitz functions on $S$, to show
that $A$-valued Lipschitz algebras never have the \trp\ if $S$
is not uniformly discrete.  These results continue results of
\cite{Gourdeau}, where it is shown that $\Lip_\alpha(S)$ is amenable
\iff\ $S$ is uniformly discrete, and that if $A$ is a \ba\ and
$\Lip_\alpha(S; A)$ is amenable then $S$ is uniformly discrete.

The algebra \Index{$C^1(\R)$} of continuously differentiable bounded functions
with bounded derivative is a \ba\ under pointwise algebra operations and norm
$$
\norm{f}=\supof{|f|}+\supof{|f'|}.
$$
$C^1(\R)$ is an \oa, since it is a closed subalgebra of
$\Lip_1(\R)$. 
Denoting by $M_f$ the multiplication operator on $L^2(\R)$ obtained
from $f\in C^1(\R)$, 
there is a contractive \rep\ $\theta:C^1(\R)\ra
\B(L^2(\R))\ot M_2$ given by 
$$
\theta(f)=\TxT M_f, M_{f'}, 0, M_f;.
$$
This is a contractive \rep\ which does not yield the \rp, since the subspace 
$$
\left\{\Tvec \xi, 0; : \xi\in L^2(\R)\right\} 
$$
is a submodule which is not complemented by any submodule. From the
point of view of the cohomological description of the \trp, this
reflects the fact that the map $C^1(\R)\ra L^2(\R)$ given by $f\mapsto
f'$ is a non-inner bounded derivation of $C^1(\R)$.
Thus $C^1(\R)$ is an \oa\ without the \trp.

A similar failure of the \trp\ is found with the \Index{disc algebra}
\Index{$A(D)$} of bounded \cts\ functions on the unit disc $D$ which
are analytic on the interior of $D\subseteq\C$. 
For $x$ in the interior of $D$
the point derivation
$f\mapsto f'(x)$ is a  non-inner derivation into $\C$ with respect to
the module action $f\cdot \xi=f(x) \xi$. Alternatively, the natural
\rep\ of $A(D)$ on $L^2(S^1)$ has many non-complemented invariant
subspaces.

\makeatletter
\let\oldsect\@sect
\def\@sect#1#2#3#4#5#6[#7]#8{%
  \ifnum #2>\c@secnumdepth
    \let\@svsec\@empty
  \else
    \refstepcounter{#1}%
    \protected@edef\@svsec{\@seccntformat{#1}\relax}%
  \fi
  \@tempskipa #5\relax
  \ifdim \@tempskipa>\z@
    \begingroup
      #6{%
        \@hangfrom{\hskip #3\relax\@svsec}%
          \interlinepenalty \@M #8\@@par}%
    \endgroup
    \csname #1mark\endcsname{The \rp\ does not imply the complete r.p.}%
    \addcontentsline{toc}{#1}{%
      \ifnum #2>\c@secnumdepth \else
        \protect\numberline{\csname the#1\endcsname}%
      \fi
      #7}%
  \else
    \def\@svsechd{%
      #6{\hskip #3\relax
      \@svsec #8}%
      \csname #1mark\endcsname{The \rp\ does not imply the complete r.p.}%
      \addcontentsline{toc}{#1}{%
        \ifnum #2>\c@secnumdepth \else
          \protect\numberline{\csname the#1\endcsname}%
        \fi
        #7}}%
  \fi
  \@xsect{#5}}

\dealwithsectionbreaks
\section{The \rp\ does not imply the \crp}
\label{rp does not imply crp}
\def\M{{\cal M}}
\let\@sect\oldsect
\makeatother

In this section\galley{Careful} we construct an \oa\ $A\subseteq\BH$
such that $\pma\H n$ has the \rp\ for all finite $n$, but $\H$ does
not have the \crp.  The demonstration that $\H$ has the \rp\ but not
the \crp\ is straightforward. To show that $\pma\H n$ has the \rp\ for
all finite $n$ involves some slightly fiddly but essentially
straightforward calculations with unbounded operators.

Let $\H_1$ be an infinite-dimensional \hs, $\H_2=\H_1$ an identical
copy of $\H_1$, and put $\H=\H_1\oplus
\H_2$. We write $p_1$ for the coordinate projection $\xi_1\oplus\xi_2\mapsto
\xi_1$ and similarly write $p_2$ for the coordinate projection 
$\xi_1\oplus \xi_2\mapsto \xi_2$.
Let $T:\H_1\ra\H_2$ be any \bded, one-one, dense-ranged
operator which is not bounded below.  Denote by $\gr T$ the subspace
$\{\xi\oplus T\xi\}\subseteq \H$, and consider the algebra
$$
A=\alg\{\H_1\oplus 0, 0\oplus\H_2, \gr T\}.
$$
This algebra will
satisfy our requirements.  For the remainder of this section 
$\H_1$, $\H_2$, $\H$, $T$ and $A$ will keep these fixed meanings.

\begin{lemma}\label{Matrix form}
  Every element of $A$ is a matrix of the form
$$
\TxT a_{11}, 0, 0, a_{22};
$$
where $T a_{11}  =  a_{22} T$.
  For all $\xi, \zeta\in \H_1$  the operator
$$
a_{\xi,\zeta}=\TxT \xi\ot T^*\zeta, 0, 0, T\xi\ot \zeta;
$$
is an element of $A$, where $\alpha\ot \beta$ denotes the rank one
operator $\eta\mapsto (\eta|\beta)\alpha$.
\end{lemma}
\proof
Since $\H_1\oplus 0$ and $0\oplus \H_2$ are submodules, elements
of $A$ are of the form $\TxT a_{11}, 0, 0, a_{22};$. The remaining
assertions are verified by the equalities 
$$
\TxT a_{11}, 0, 0, a_{22}; \Tvec \eta, T \eta; = \Tvec a_{11}\eta,
a_{22} T \eta;\in \Gr  T
$$
and
$$
\TxT \xi\ot T^*\zeta, 0, 0, T\xi\ot \zeta; \Tvec \eta, T\eta; 
= \Tvec \xi (T\eta|\zeta), T\xi(T\eta|\zeta);
$$
\qed

\begin{lemma}\label{Lattice is what it is}
  $\lat A=\{0\oplus\H_2\} \cup \{ \gr \lambda T: \lambda\in \C\}$.  If
  $\M\subseteq\H$ is a non-closed $A$-invariant submanifold, then
  $p_1\M\subseteq\M$ and $p_2\M\subseteq\M$. Moreover, $p_1\M\in\{0,\H_1\}$
  and  $T\H_2\subseteq p_2\M\subseteq\H_2$.
\end{lemma}
\proof Firstly, observe that the matrix form given for elements of $A$
in the above \nameref{Matrix form} implies that $\gr \lambda T\in \lat
A$ for all $\lambda\in \C$.  The operators $a_{\xi,\zeta}$ show that
$A$ acts strictly transitively on $\gr\lambda T$ for all $\lambda$,
and transitively on $0\oplus\H_2$.  Therefore these submodules are
irreducible.

If $\V\in\lat A$ contains any two of the subspaces $\{\gr \lambda T,
0\oplus \H_2\}$, it also contains a vector from $\H_1\oplus0$.  Since
$A$ acts strictly transitively on $\H_1\oplus0$, $\V$ must then contain
$\H_1\oplus0$, and hence it contains a vector from $0\oplus\H_2$. As $A$
acts transitively on $0\oplus\H_2$, $V$ also contains $0\oplus\H_2$,
so $\V=\H$ in this case.

Now let $\V\in\lat A$ be an arbitrary proper submodule of $\H$, and let
$ \eta_1\oplus \eta_2$ denote a general element of $\V$. For every
$\xi, \zeta\in\H_1$ we have
$$ a_{\xi,\zeta}\Tvec \eta_1,\eta_2; = \TxT \xi\ot T^*\zeta, 0, 0,
T\xi\ot \zeta; \Tvec \eta_1, \eta_2; = \Tvec \xi(T\eta_1|\zeta),
T\xi(\eta_2|\zeta); \in \V.
$$
This vector is an element of either $0\oplus\H_2$ or $\gr \lambda T$
for some $\lambda\in \C$, depending on the 
point in $\C P_1$ with homogeneous 
coordinates
$((T\eta_1|\zeta),(\eta_2|\zeta))$. If there are two vectors in
$\V$ which yield different points in $\C P_1$ then by the comments
above $V=\H$. 
Since $\V$ is assumed proper,  all vectors in $\V$ yield the same point in
$\C P_1$; it is easy to see that this implies that either $\V=0\oplus\H_2$ or $\V=\gr \lambda
T$ for some $\lambda\in \C$. This establishes the claim about $\lat A$.

Let $\M\subseteq\H$ be a non-closed invariant submanifold. The closure
$\clos\M$ is a submodule of $\H$.
The strict transitivity of $A$ on $\gr\lambda T$ and the assumption
that $\M$ is not closed implies 
$\M\not\subseteq\gr \lambda T$ for any
$\lambda$. 
In particular, $\M\not\subseteq\H_1\oplus0$. On the other hand, if
$\M\subseteq 0\oplus\H_2$, then multiplication by a suitable
$a_{\xi,\zeta}$ shows that $\M\cap (0\oplus T\H_1)\neq\{0\}$, and the
strict transitivity of $A$ on $0\oplus T\H_1$ shows that $0\oplus
T\H_1\subseteq \M \subseteq 0\oplus\H_2$. 

Suppose that $\M\not\subseteq 0\oplus \H_2$.  Let
$\xi=\xi_1\oplus\xi_2$ be an arbitrary element of $\M$ with $\xi_1\neq
0$ and $\xi_2\neq 0$. For suitable
$\gamma$ we have $a_{\xi_1,\gamma}\xi=\xi_1\oplus\lambda T\xi_1\in\M$
for some $\lambda\in\C$. 
If $\xi_2\not\in\C T\xi_1$, then subtracting gives
$0\neq\xi-a_{\xi,\gamma}\xi\in 0\oplus \H_2$, and so $0\oplus
T\H_1\subseteq\M$. But then $\xi_1\oplus 0\in\M$ and so $\H_1\oplus
0\subseteq \M$. In this circumstance $\M$ satisfies the statement in
the \nameref{Lattice is what it is}.

Finally, suppose $\xi=\xi_1\oplus\xi_2\in\M$ implies $\xi_2\in\C
T\xi_1$ for all $\xi\in\M$. Since we know $\M\not\subseteq\gr \lambda
T$ for all $\lambda\in\C$, there are non-zero vectors
$\eta=\eta_1\oplus \lambda_1 T\eta_1\in \M$ and $\rho=\rho_1\oplus
\lambda_2 T\rho_1\in\M$ with $\lambda_1\neq \lambda_2$. Choosing a
suitable $\gamma$ from $\H_1$, we obtain
$\rho'=a_{\eta_1,\gamma}\rho\in\gr \lambda_2\cap \M$ and $0\neq
\eta-\rho'\in 0\oplus T\H_1$. Then as before $0\oplus
T\H_1\subseteq\M$ and $\H_1\oplus0\subseteq\M$.  Again, $\M$ splits
into $\M=p_1\M + p_2\M$, and the components $p_1\M$ and $p_2\M$ are as
stated in the \nameref{Lattice is what it is}.  \qed

\begin{proposition}\label{RP but not CRP}
  The algebra $A$ has the \rp, but not the \crp.
\end{proposition}
\proof
The submodule $0\oplus\H_2$ is complemented by $\H_1\oplus 0$ in $\lat
A$. The submodules $\gr \lambda T$ are complemented (uniquely) by $0\oplus
\H_2$ in $\lat A$. Thus $A$ has the \rp.

\Nref{Existence of \projconst{}s for crp} shows that if $A$ has the
\crp\ then there is $K\geq 1$ so that every submodule in $\lat A$ is the range of a
module projection $p$ with $\norm{p}\leq K$. However, as
$\lambda\ra\infty$ the norm of the projection onto $\gr \lambda T$
along the unique complement $0\oplus\H_2$ increases without
bound. Hence $A$ does not have the \crp.
\qed


This example shows that the \rp\ does not imply the \crp. In fact, we
can show that $\pma\H n$ has the \rp\ for all finite $n$. To motivate
the following calculation, consider the related \oa\ $B=\alg\{\H_1\oplus0,
0\oplus\H_2, \gr Q\}$ where $Q\in\B(\H_1,\H_2)$ is invertible. Easy
calculations akin to those for $A$ 
show that the operators in $B$ have the matrix form 
$$
\TxT b, 0, 0, Qb\inv Q;
$$
where $b\in\B(\H_1)$ is arbitrary. In fact, $B$ is similar to the
algebra $D=\discretionary{\hbox{\vadjust{\nobreak}}}{}{}\pma{\B(\H_1)} 2 = \B(\H_1)\ot 1_2$ via the similarity 
$$
\TxT 1, 0, 0, Q;,
$$
and so the invariant subspace lattice for $B$ is
isomorphic to the invariant subspace lattice of $D$. 

It is easily shown that the submodules of $D$ are the subspaces of the
form $\H_1\ot_2 \W\subseteq\H_1\ot_2\C^2$, where $\W\subseteq\C^2$ is
any subspace. Similarly, the invariant subspaces of $\pma Dn$ are the
subspaces of the form $\H_1\ot_2 \W'$ where $\W'\subseteq \C^{2n}$ is
again an arbitrary subspace. Accordingly, the invariant subspaces of
$\pma B n$ are the subspaces of the form
$$
\left[\begin{array}{ccccc}
            1&&&&\\
            &Q&&&\\
            &&1&&\\
            &&&\ddots&\\
            &&&&Q
\end{array}\right] (\H_1\ot_2 \W')=\TxT 1, ,, Q;\ot 1_n \;(H_1\ot_2\W').
$$
Returning to $A$, similar manipulations show that the subspaces
$$
\TxT 1, 0, 0, T; \ot 1_n \;(\H_1\ot_2 \W')
$$
(or their closures when they are not closed) 
are submodules for $\pma A n$. In fact, the following
calculations will show that these account for all submodules of $\pma A n$.
For brevity, let us write $T_n$ for the operator
$$
\TxT 1, 0, 0, T; \ot 1_n:\pma \H n\ra \pma \H n.
$$

\index{graph transformation}
The first step is to determine the graph transformations for
$A$. Recall that a linear operator $S:X\ra \H$ defined on some linear
submanifold $X\subseteq \H$ is a graph transformation for $A$
if $AX\subseteq X$ and $aS=Sa$ for all $a\in A$. 
\begin{lemma}\label{Decompose graph transforms}
  Let $S$ be a graph transformation for $A$. Then there are complex numbers
  $\lambda_{ij}$ such that $S$ is given by the restriction of the
  operator 
$$
\TxT \lambda_{11}, \lambda_{12}\inv T, \lambda_{21}T, \lambda_{22};
$$ to $\dom S$.  
\end{lemma}
\proof
Note that the domain of $S$ is an $A$-invariant manifold, and so
it is of the form described in \nref{Lattice is what it is}.

To begin we show that there are always (possibly unbounded) operators
$S_{ij}$ such that $S$ is the restriction of the operator
$$
\TxT S_{11}, S_{12}, S_{21}, S_{22}; 
$$
to $\dom S$. In the cases where $\xi_1\oplus \xi_2\in \dom S$ implies that
$\xi_1\oplus 0, 0\oplus \xi_2 \in \dom S$, then this is
elementary. 
By \nref{Lattice is what it is}
the only situation where this does not happen is when $\dom S=\gr
\lambda T$ for some $\lambda\in \C$. 
In this case there are linear operators $S_{11}$
and $S_{21}$ such that $S(\xi\oplus\lambda T\xi)=S_{11}\xi\oplus
S_{21}\xi$. Then $S$ equals the restriction of the matrix operator
$$
\TxT S_{11}, 0, S_{21}, 0;
$$
to $\gr \lambda T$.

This means that in all cases $S$ may be written in terms of matrix
elements.
Moreover, if $\dom S=\gr\lambda T$ for some $\lambda\in\C$ we
can choose $S_{12}=S_{22}=0$, and if $\dom
S\subseteq 0\oplus\H_2$ we can choose $S_{11}=S_{21}=0$. 
To show that the matrix elements are of the required form, we 
calculate the products $a_{\xi,\zeta}S$ and $Sa_{\xi,\zeta}$. 
For convenience of
notation, we refer to the linear map $\eta\mapsto (S_{ij}\eta|\zeta)$ as
$S'_{ij}\zeta$. We have 
\begin{eqnarray*}
a_{\xi,\zeta}S&=&\TxT \xi\ot S'_{11}T^*\zeta, \xi\ot S_{12}'T^*\zeta, 
T\xi\ot S'_{21}\zeta, T\xi\ot S'_{22}\zeta;,\\
Sa_{\xi,\zeta}&=&\TxT S_{11}\xi\ot T^*\zeta, S_{12}T\xi\ot \zeta,
S_{21}\xi\ot T^*\zeta, S_{22}T\xi\ot \zeta;.
\end{eqnarray*}
Since the operators $a_{\xi,\zeta}S$ and $Sa_{\zeta,\xi}$ agree on
$\dom S$ for any $\xi, \zeta\in \H_1$, we may read off several 
equations involving the operators $S_{ij}$. 

Consider firstly the case
where $p_1(\dom S)\subseteq\dom S$ and $p_2(\dom S)\subseteq\dom S$.
\begin{enumerate}
\item Since $T\xi\ot S'_{22}\zeta=S_{22}T\xi\ot \zeta$
  on $p_2(\dom S)$ for all $\xi, \zeta\in \H_1$, we have $S_{22}\in \C$.
\item Similarly, the equation $\xi\ot S'_{11}T^*\zeta=S_{11}\xi\ot
  T^*\zeta$ implies that $S_{11}\in \C$.
\item Since $\xi\ot S'_{12}T\zeta=S_{12}T\xi\ot \zeta$, we have
  $S_{12}T=S'_{12}T^*\in \C$, so $S_{12}\in \C\inv T$. 
\item Finally, $T\xi\ot S_{21}'\zeta=S_{21}\xi\ot T^*\zeta$ implies
  that $S_{21}\in\C T$.
\end{enumerate}

In the case where $\dom S\subseteq 0\oplus\H_2$ the matrix form of $S$
is
$$
S=\TxT 0, S_{12}, 0, S_{22};
$$
and similar arguments give the same equations. Finally, in all other
cases we have written $S$ in the form
$$
S=\TxT S_{11}, 0, S_{21}, 0;$$
and again similar arguments give the same equations for the $S_{ij}$.
This establishes the matrix form
$$
S=\TxT \lambda_{11}, \lambda_{12}\inv T, \lambda_{21} T, \lambda_{22};
$$
for suitable $\lambda_{ij}$.
\qed

If $\H\supseteq\V_1\supseteq\V_2$ is a chain of \hs{}s, we write
$\V_1\ominus\V_2$ for $\V_2^\perp\cap\V_1$.  We will need the
following elementary linear algebraic fact.
\begin{lemma}
\label{Linear algebraic fact}
  Suppose that $\Y$ is a finite dimensional \hs\ with
  complementary subspaces $\Y_1$ and $\Y_2$, and let $p$ denote
  the projection onto $\Y_1$ along $\Y_2$. 
  If $\V\subseteq\Y$ is any subspace, there
  is a projection $q:\Y\ra \V$ with $(1-p)qp=0$. That is,
  $qY_1\subseteq Y_1$.
\end{lemma}
\proof
Let $\V_1=\V\cap \Y_1$, and $\V_2=\V\cap
\big((\Y_1\ominus\V_1)\oplus\Y_2\big)$. Then $\V=\V_1\oplus\V_2$, and $\V_2$ is
a graph subspace over $(1-p)\V\subseteq\Y_2$. For brevity we write
$\dom \V_2$ for $(1-p)\V$.
The subspace $\W=(\Y_1\ominus\V_1)\oplus(\Y_2\ominus \dom\V_2)$ is a
complement to $\V$, and if we decompose $\Y_1$ as $\V_1 \oplus
(\Y_1\ominus \V_1)$ and $\Y_2$ as $\dom\V_2 \oplus (\Y_2\ominus \dom\V_2)$ then the
projection $q$ onto $\V$ along $\W$ has the matrix form
$$
q= \TxT {\TxT 1, 0, 0, 0; }, {\TxT 0, 0, Q_2, 0; }, 
     {\TxT 0, 0, 0, 0; }, {\TxT 1, 0, 0, 0; };
\begin{array}{l}
  \left.\mathstrut\right\}\V_1\\
  \left.\mathstrut\right\}\Y_1\ominus\V_1\\
  \left.\mathstrut\right\}\dom \V_2\\
  \left.\mathstrut\right\}\Y_2\ominus\dom\V_2
\end{array}
$$
where $ Q_2 :\dom\V_2\ra \Y_1\ominus\V_1$ is the graph transformation
associated with $\V_2$.  The operator $(1-p)qp$ corresponds to the
bottom left hand corner of the matrix, which is zero.  \qed

\begin{lemma}
  \label{silly subspace is closed}
  Let $\M\subseteq\H$ be a not-necessarily-closed invariant
  submanifold of $\H$. The set $\W=\{\xi\in \H: T_1\xi\in\M\}$ is
  a closed subspace of $\H$.
\end{lemma}
\proof If $\M$ is closed this is clear. Otherwise, the conclusion
follows from the form for $\M$ given by \nref{Lattice is what it is}.
\qed

\begin{proposition}
  The $A$-module $\pma\H n$ has the \rp\ for all $n\in \N$. 
\end{proposition}
\proof We proceed by induction on $n$.  The inductive hypothesis will
be: 
\begin{enumerate}
\item $\pma\H n$ has the \rp.
\item if $\V\subseteq\pma \H n$
is $A$-invariant then there is a dense submanifold $\U\subseteq \V$ with
$\U\subseteq\dom \inv {T_n}$, such that the subspace
$
\W=\inv {T_n} (\U)
$ is a closed invariant subspace for $\pma D n$.
\end{enumerate}

\Nref{RP but not CRP} shows that $\pma \H 1=\H$ has the \rp, and 
\nref{Lattice is what it is} shows that the statement about dense
submanifolds is true for $n=1$. We consider $\pma\H{n-1}$ as a
submodule of $\pma\H n$ via the embedding
$j:(\xi_1,\ldots,\xi_n)\mapsto(0,\xi_1,\ldots,\xi_n)$. 

Let $\V\subseteq\pma\H n$ be $A$-invariant. Define 
$\V_0=\V \cap j(\pma \H {n-1})$. Then $\V_0$ can be considered
as a submodule of $\pma \H {n-1}$, and by the inductive hypothesis
$\V_0$ is complemented in $\pma\H {n-1}$ by a submodule $\V_0'$. 
We write $\H\oplus\V_0'$ for the submodule
$(\H\oplus0\oplus\cdots\oplus0)+j\V_0'$, and 
define $\V_1=\V \cap (\H\oplus \V_0')$. Then $\V$ is the internal direct sum of
$\V_0$ and $\V_1$.

By construction, if $\xi=(0,\xi_2,\ldots,\xi_n)\in\V_1$ then $\xi=0$.
Thus there are linear transformations $S_2,\ldots,S_n$ defined on some
common domain in $\H$ with
$\V_1=\{(\xi_1,S_2(\xi_1),\ldots,S_n(\xi_1)\}$. Let us denote this
common domain by $\dom \V_1$. The invariance of $\V_1$ implies that
the maps $S_i$ are graph transformations over $\dom \V_1$.

Using 
\nref{Decompose graph transforms} we may find numbers
$\lambda^i_{kl}\in \C$ so that 
$$
S_i=\TxT \lambda^i_{11}, \lambda^i_{12} \inv T, \lambda^i_{21} T, \lambda^i_{22};.
$$
We define a $D$-invariant subspace $\W_1\subseteq \H_1\ot_2\C^{2n}$ by
$$
\W_1=\left\{\left(
\xi, \TxT \lambda^2_{11},\lambda^2_{12},\lambda^2_{21},\lambda^2_{22};\xi, 
      \cdots, \TxT
      \lambda^n_{11},\lambda^n_{12},\lambda^n_{21},\lambda^n_{22};\xi
\right): T_1\xi\in\dom\V_1\right\}
$$ 
By \nref{silly subspace is closed}, $\W_1$ is a closed graph subspace
over $\dom\W_1=\{\xi:T_1\xi\in\dom\V_1\}$.  Furthermore, $\W_1$ is a
$D$-invariant submodule of $\pma \H n$.

If none of the graph transforms $S_i$ are unbounded then 
$\dom \V_1$ is closed and 
$T_n\W_1$ is
a dense submanifold of $\V_1$. On the other hand, if any of the $S_i$
are unbounded, this implies that one of the $\lambda_{12}^i\neq 0$, in
which case $T_n\W_1=\V_1$. In any case $T_n\W_1$ is a dense
$A$-invariant submanifold of $\V_1$.

By induction there is a dense submanifold
$\U_{0}\subseteq\V_0$ and a $D$-invariant subspace
$\W_{0}\subseteq\pma \H {n-1}$ 
such that $\U_{0}=T_n \W_{0}$. 
We put $\W=\W_0 + \W_1$, making
$T_n\W$ a dense submanifold of $\V$.

To complete the proof, we need only find a 
 $D$-module projection  $q:\pma\H n\ra \W$
such that $p=T_n q \inv {T_n}$  is bounded. 
If this can be done, 
then $p$ fixes $T_n\W$ and is thus a projection
onto a subspace containing $\V$. Also, on the dense submanifold
$T_n\pma\H n\subseteq\pma\H n$, the projection $p$ maps into $T_n\W$,
and so $p$ is a projection onto $\V$. Finally,  an
easy calculation shows that $p$ commutes with $A$ on the dense
submanifold $T_n\pma\H n$, and so $p\in A'$.
This will show  that $\V$ is complemented by an $A$-submodule, whence 
$\pma \H n$ has the \rp\ for all $n\in\N$.

In order to obtain the desired $D$-module projection $q$, observe that
$\W=\H_1\ot_2 \W'$ for some $\W'\subseteq \C^{2n}$, and a $D$-module
projection onto $\W$ is any operator of the form $q=1_{\H_1}\ot q'$
where $q':\C^{2n}\ra \W'$ is an arbitrary linear projection onto
$\W'$.  If we decompose the space $\C^{2n}$ as the direct sum of
$\Y_{odd}=\C\oplus 0 \oplus \C \cdots 0$ and
$\Y_{even}=0\oplus\C\oplus0\oplus\cdots\oplus\C$, then matrix
multiplication shows that the condition on $q'$ to ensure
the boundedness of $T_n q \inv {T_n}$ is precisely that
$q'\Y_{even}\subseteq\Y_{even}$.  \Nref{Linear algebraic fact} ensures
that we can find such a $q'$.  \qed

The existence of \RA{}s which are not \cra{}s is clearly related to
the failure of the \rp\ to imply the existence of a finite projection
constant. This contrasts with a corollary of the complemented
subspaces theorem: 
if $X$ is a \bs\ 
every subspace of which is complemented, then there is $\lambda\in\C$
such that every subspace of $X$ is the range of a projection with norm
$\leq \lambda$ \cite{DavisDeanSinger},
\cite{LindenstraussTzafriri}. The proof of the complemented subspaces
theorem uses the abundance of
finite-dimensional subspaces of a \bs. 
It is the lack of an equivalently rich invariant subspace lattice
which thwarts any possibility of a similar argument working here.

\comment{ Theorem is busted for the time being
The next result clarifies the connection between the \crp\ and the
boundedness of the \projconst{}s for $\pma\H n$.
\begin{lemma}
  Let $A$ be an \oa\ with \hm\ $\H$. If there exists a $K>0$ such that
  $\pma\H n$ has the \rp\ and a finite \projconst{} $K(\pma\H
  n)\leq K$ for each $n\in \N$, then $\H$ has the \crp.
\end{lemma}
\proof
Let $\V\subseteq\pma\H \infty$, and define $\V_n=\V\cap\pma\H n$. Then
$\clos{\bigcup \V_n}=\V$. Since each $\V_n$ is a submodule of $\pma\H n$,
there is a module projection $p_n:\pma\H n\ra \V_n$ with
$\norm{p_n}\leq K$. If $q_n:\pma\H \infty\ra \pma\H n$ denotes the
orthogonal projection onto $\pma\H n$ then $\{q_np_n\}$ is a bounded
sequence of projections. Since the closed ball of $\B(\pma\H\infty)$
is \ws-closed, this sequence has a cluster point, $p$ say. Then $p$
is a module projection onto $\V$.
\qed
}

\vfill
\break

\Nref{RP but not CRP} and the example on page~\pageref{Page with RP
  example} show that each of the implications
$$
\begin{array}{c}\hbox{total reduction}\\\hbox{property}\end{array}\implies
\begin{array}{c}\hbox{complete reduction}\\\hbox{property}\end{array}\implies
\begin{array}{c}\hbox{reduction}\\\hbox{property}\end{array}
$$
is strict.

\chapter{Properties of \RAs}
\label{Chapter with properties}

In this chapter we study properties of algebras with one of the three
reduction properties. Since the \rp\ is defined in terms of
modules, it is not surprising that the first results are spatial in
character. Using the GNS construction, we are able to apply these
spatial results to obtain algebraic results, chiefly concerning the
ideal structure of \tras. The stability of the \trp\ is investigated
in section~\ref{Section with stability}. In section~\ref{Section with
  twine} some more spatial results are developed, which concentrate
on the relationships between modules of a \cra.

\dealwithsectionbreaks
\section{Spatial properties of \RAs}
\label{Section with properties}

It is the central theme of this work that \TRAs\ behave in many ways
like \csalgs. In this section, we develop some 
\csalg-like properties of \TRAs\ and, more generally, of \cra{}s. We
start with some technical lemmas which lead up to a proof of the
`double commutant theorem' for \oas\ with the \crp\ (\nref{Double
  commutant theorem}). 

\csalgs\ and their closed ideals always contain \bais\ \cite{DixmierC}. In
fact this is also true for \tra{}s (\nref{Ideals have bai}), but until
this is established the following \nameref{Nondegeneracy lemma} will
suffice as a replacement.

It is worth noting that while  we are mainly interested in closed
ideals, many of the spatial results of this section apply to non-closed ideals as
well. For this reason any closure requirements on ideals will be
declared explicitly in this and the subsequent section.

\begin{lemma} \label{Nondegeneracy lemma}
  Let $A$ be an \oa\ and suppose $J\ideal A$ is a two-sided ideal.
  Suppose further that $\H$ is a \hAm A with the \rp, and
  $\clos{J\H}=\H$.  Then $\xi\in\closure{J\xi}$ for all $\xi\in\H$.
\end{lemma}
\proof For every $\xi\in\H$, the space $\closure{J\xi}$ is
$A$-invariant.  The \rp\ yields $\V\in\lat A$ with
$\H=\closure{J\xi}\oplus\V$, giving a decomposition
$\xi=\xi_j\oplus \xi_v$ into the respective components.  We have
$J\xi=J(\xi_j\oplus\xi_v)\subseteq \clos{J\xi}$, so $J\xi_v=0$.  We
wish to show that $\xi_v=0$. Consider the subspace $\W=\{\zeta:
J\zeta=0\}$.  Since $J$ is a right ideal, this is a closed $A$-module
containing $\xi_v$, and using the \rp\ there is $\U\in\lat A$ with
$\H=\W\oplus \U$.  Then by the nondegeneracy assumption
$\H=\clos{J\H}=\pclos{J\W\oplus J\U}\subseteq \clos{J\U}$, implying
that $\U=\H$. Thus $\W=\{0\}$ and so $\xi_v=0$.  \qed
\comment{ Deleted because never used.
\begin{corollary}
  Let $A$ be an \oa\, with a \hm\ $\H$, and a two-sided ideal
  $J\ideal A$. Suppose that $\clos{J\H}=\H$ and that $\H$ has the
  \rp\ as a $J$-module. Then $\xi\in \clos{J\xi}$ for all $\xi\in \H$.
\end{corollary}
\proof
Apply the above \nameref{Nondegeneracy lemma} with $A=J$.
\qed
}
\begin{corollary}
\label{First nondegeneracy}
Let $A\subseteq\BH$ be an \oa\ 
such that $\H$ has the \crp.
Suppose that 
$J\ideal A$ is a two-sided ideal of $A$ with $\clos{J\H}=\H$.
For any $1\leq n\leq\infty$ and any $\xi\in\H^{(n)}$ we have 
$\xi\in\closure{\pma Jn\xi}$.
\end{corollary}
\proof
It follows from $\clos{J\H}=\H$ that $\clos{\pma
  J n\pma\H n}=\pma\H n$. Then applying
\nref{Nondegeneracy lemma} to the amplified \rep\ gives the result.
\comment{
The subspace $\clos{J^{(n)}\xi}$ is an
$A$-module.
Using
 the \rp\ for $\H^{(n)}$, we have
$\H^{(n)}=\closure{J^{(n)}\xi}\oplus \V$, and again we get a decomposition
$\xi=\xi_j\oplus\xi_v$, where $J^{(n)}\xi_v=0$. If
$\xi_v=(\xi_v^i)_{i=1}^n$ is the decomposition of $\xi_v$ into
coordinate components, this says that $J\xi_v^i=0$ for all $i$. But
from the proof of \nref{Nondegeneracy lemma} this implies that
$\xi_v^i=0$ for all $i$, and so $\xi_v=0$. Thus
$\xi=\xi_j\in\clos{J^{(n)}\xi}.$
}
\qed
\begin{corollary}
\label{Second nondegeneracy}
  Let $A$ be an \oa\ with \hm\ $\H$, and let 
  $J\ideal A$ be a  two-sided ideal. If
  $\clos{J\H}=\H$ and $\H$ has the \crp\ as a $J$-module, 
  then for any $1\leq n\leq\infty$ and  $\xi\in\H^{(n)}$ we have
  $\xi\in\clos{J\xi}$. 
\end{corollary}
\proof
Apply the above \nameref{First nondegeneracy} with $J$ in place of
$A$.
\qed
\begin{theorem} \label{Bare Double Commutant}
\label{Bare double commutant}
\label{Double commutant}
\label{Double commutant theorem}
Let $A$ be an \oa\ 
and let $J\ideal A$ be a two-sided ideal.
Suppose that $\theta: A\ra \BH$ is a \rep\ of $A$ such that
$\H$ has the \crp\ as either an $A$-module or a $J$-module (not
necessarily both).
If $\clos{J\H}=\H$, then
$\swclos{\theta(J)}=J''$. 
%
%
\end{theorem}
\proof It is immediate from the separate \sw\ continuity of
multiplication in $\BH$ that $\swclos{\theta(J)}\subseteq J''$. For
the converse, take $T\in J''$ and a sequence
$\xi=(\xi_i)\in\pma\H\infty$. We will show that for any $\eps>0$ there
is $a\in J$ with $\norm{T^{(\infty)}\xi-a^{(\infty)}\xi}<\eps.$ This
will demonstrate that the \sstrong\ closure of $\theta(J)$ contains
$J''$; since the \sweak\ topology is weaker than the \sstrong\ 
topology, the result will follow.

Using \nref{First nondegeneracy} when $\pma\H\infty$ has the \rp\ as an
$A$-module and \nref{Second nondegeneracy} when $\pma\H\infty$ has the
\rp\ as a $J$-module,
the vector  $\xi$ lies in $\closure{J^{(\infty)}\xi}$.
The \rp\ gives a 
projection in $(J^{(\infty)})'$ onto
$\closure{J^{(\infty)}\xi}$. 
Since
$T\in J''$, it follows that 
$T^{(\infty)}\in (J^{(\infty)})''$, and so
$T^{(\infty)}\xi\in \pclos{J^{(\infty)}\xi}$. 
  Consequently, there is $a\in J$ with 
$\norm{T^{(\infty)}\xi-a^{(\infty)}\xi}<\eps$.
\qed

\begin{corollary}\label{Double commutant is closure}
  Let $A\subseteq\BH$ be a nondegenerate \oa\ with the \crp.
  Then $A''=\swclos{A}$.
\qed
\end{corollary}

For \csalgs\ this is the famous von Neumann \Index{double commutant
theorem} \cite{ArvesonBook}. 
The proof given here is a modification of the usual \csalg\ proof.

We finish this section with a series of technical \nameref{Restriction
  is closed}s which  allow us to deal with the possibility of
nondegenerate \reps.

\begin{lemma}\label{nondegenerate on AH}
  Let $A\subseteq\BH$ be a \RA, and let $J\ideal A$ be a left ideal.
  Then $J$ acts nondegenerately on $\clos{J\H}$.
\end{lemma}
\proof
The subspace $\clos{J\H}$ is an $A$-module so
there is $\W\in\lat A$ with $\clos{J\H}\oplus\W=\H$. 
We have $J\W\subseteq A\H\cap \W=\{0\}$, so $J\W=\{0\}$. 
Then
$\clos{J\H}=\clos{J(\clos{J\H})\oplus
  J\W}=\clos{J(\clos{J\H})}$. 
Thus $J$ acts nondegenerately
on $\clos{J\H}$. 
\qed

As a matter of notation, if $B\subseteq\BH$ is a set of operators with
an invariant subspace $\V\subseteq\H$, we write $B|_{\V}$ for the set
$\{b|_\V:b\in B\}\subseteq\B(\V)$.

\begin{lemma}\label{Restriction is closed}
  Let $A\subseteq\BH$ be a \RA, and let $J\ideal A$ be a left ideal.
  Then $\swclos{{J}}|_{\ssclos{J\H}}=
  \swclos{{J}|_{\ssclos{J\H}}}$.  Consequently,
  $\swclos{{J}|_{\ssclos{J\H}}}$ is \sw-homeomorphically
  isomorphic to $\swclos{{J}}$.
\end{lemma}
\proof The submodule $\clos{J\H}$ is complemented by some $\W\in\lat A$
with $J\W\subseteq\clos{J\H}\cap\W=\{0\}$.  We may apply a
similarity to $\H$ to orthogonalise $\W$ and $\clos{J\H}$, since this
induces a \sw\ homeomorphism on $\BH$ and will not affect the
conclusion.

The restriction map 
$\swclos{{J}}\ra \swclos{{J}|_{\clos{J\H}}}$
is a \sw\ \cts\ contractive algebra
isomorphism. We need only to show that it is open.
Suppose that $a_\nu$ is a net in $J$ and $a_\nu|_{\clos{J\H}}\ra b\in
\B(\clos{J\H})$ in the \sw\ topology on $\B(\clos{J\H})$. Define an
extension of $b$ to an element of $\B(\H)$ by $b(\W)=0$. We check that
$a_\nu\ra b$ \sw{}ly in $\BH$.  To do this, let $\{\xi_i\}$ and $\{\zeta_i\}$
be square-summable sequences in $\H$. Let $\xi_i=\xi_i^0 \oplus
\xi_i^1$ and $\zeta_i=\zeta_i^0 \oplus \zeta_i^1$ be the
decompositions obtained from $\H=\clos{J\H}\oplus \W$.  Then $\sum_i
(a_\nu \xi_i | \zeta_i) = (a_\nu \xi_i^0 | \zeta_i ^ 0) \ra \sum_i (b
\xi_i^0| \zeta_i^0) = \sum_i (b \xi_i| \zeta_i)$.
\qed

\comment{
The following simple 
\nameref{Degenerate case} extends the double commutant theorem to the
degenerate case.
}
\begin{proposition}\label{Degenerate case}
  Let $A\subseteq\BH$ be a \cra. Then $\clos{A\H}$ is 
  the range of a unique central projection $p\in A'\cap A''$ and 
  $A''=\swclos{A}\oplus \C(1-p)$. 
\end{proposition}
\proof The submodule $\clos{A\H}$ is complemented by some submodule
$\W\in\lat A$ with $AW=\{0\}$. By \nref{Restriction is closed} we know
that $\swclos{A}|_{\ssclos{A\H}}$ is the \sw\ closure of
$A|_{\ssclos{A\H}}$. By \nref{Double commutant is closure} we have
$1_{\ssclos{A\H}}\in\swclos{A|_{\ssclos{A\H}}}$, so there is a unique
central element $p\in\swclos{A}$ with
$p|_{\ssclos{A\H}}=1_{\ssclos{A\H}}$.  Since $A\W=\{0\}$, it follows
that $p$ is the projection of $\H$ onto $\clos{A\H}$ along $\W$.
Since $A\subseteq A|_{\ssclos{A\H}}\oplus 0$ we have $A'=
(A|_{\ssclos{A\H}})'\oplus\B(\W)$ and $A''=
\swclos{A|_{\ssclos{A\H}}}\oplus \C(1-p)$.  \qed
\begin{corollary}
  \label{weak and sw are same}
  Let $A$ be a \sw{}ly closed \cra. Then $A$ is a weakly closed unital
  ideal in $A''$. \qed
\end{corollary}
 The fact that the \sw{}ly closed \cra{}s are also weakly closed
remedies the apparent discrepancy in the statements of the reductive
algebra problem (\nref{reductive algebra problem}) and its non-\sa\ 
analogue \nref{cra conjecture}.

\comment{
Observe that 
if $J$ possesses a \bai, then the same conclusion holds if
$J$ is only a left ideal instead of a two-sided ideal.
}

\dealwithsectionbreaks
\section{Algebraic properties of \RAs}
\label{Section with algebraic properties}

The spatial tools built up in the previous section, particularly the
double commutant theorem (\nref{Double commutant theorem}), allow us to
draw conclusions about the ideal structure of a \cra\ or \tra. We
start by considering the ideal structure of $A''$; the GNS
construction will provide the links between $A''$, $\bd A$ and hence $A$.

\begin{lemma}\label{Bicommutant has rp}
  If $A\subseteq\BH$ has the \rp\ then $A''$ has the \rp. If $A$ has
  the \crp\ then $A''$ has the \crp.
\end{lemma}
\proof Suppose that $A$ has the \rp\ and $\V\in\lat A''$. Since
$A\subseteq A''$ we have $\V\in\lat A$, and so there is a
projection $p\in A'$ onto $\V$.  Since $A'=(A'')'$, $p$ is also an
$A''$-module projection.  Hence $A''$ has the \rp.

If $A$ has the \crp\ then $\pma A\infty$ has the \rp, hence $\pma
{(A'')}\infty=(\pma A\infty)''$ has the \rp, hence $A''$ has the \crp.
\qed

\begin{theorem} \label{Ideals split} \label{Ideals are complemented}
  Let $A\subseteq\BH$ be a \cra.
  Suppose that $J\ideal A''$ is a \sweak{}ly closed two-sided ideal of
  $A''$. Then there is a unique central idempotent $j\in J\cap A'$
  with the property that $jk=k$ for all $k\in J$ (that is, $j$ is an
  identity for $J$).  \comment{ Moreover, the same conclusion holds
    for \sw-closed left ideals which possess \bais.  }
\end{theorem}
\proof The space $\clos{J\H}$ is an $A''$-module since
it is the closure of the $A''$-invariant manifold $J\H$.  Using
\nref{Bicommutant has rp}, the \rp\ yields an $A''$-module
decomposition $\H=\clos{J\H}\oplus \V$ with $J\V=\{0\}$.  By
\nref{nondegenerate on AH} $J$ acts nondegenerately on $\clos{J\H}$.
Since $\clos{J\H}$ is an $A''$-submodule of $\H$ it has the \rp\ as
an $A''$-module. Similarly, the $\infty$-fold amplification
$\pma{\clos{J\H}}\infty$ is a closed submodule of
$\pma\H\infty$, and so it also has the \rp\ as an $A''$-module.  Now
we may apply \nref{Double commutant} to conclude that
$(J|_{\ssclos{J\H}})''=\swclos{J|_{\ssclos{J\H}}}$.  By
\nref{Restriction is closed} we know that $J|_{\ssclos{J\H}}$ is \sw\ 
closed in $\B(\clos{J\H})$. Since the identity operator
$1_{\ssclos{J\H}}$ belongs to the bicommutant, this implies that
there is an element $j\in J$ such that
$j|_{\ssclos{J\H}}=1_{\ssclos{J\H}}$. Since $j\H\subseteq\clos{J\H}$,
we see that $j$ must be a projection of $\H$ onto $\clos{J\H}$.  
For any $k\in J$, we have $jk=k$. Further, since $J\V=\{0\}$, we have
$j\V=\{0\}$. Because $\clos{J\H}\oplus\V=\H$, $j$ must be the
projection onto $\clos{J\H}$ along $\V$. Being a projection onto an
$A''$-invariant space with an $A''$-invariant kernel means that $j\in
(A'')'=A'$, and so $j$ is a central projection. Consequently $j$ is an identity
for $J$, and it is the unique element of $A''$ with this property.
\qed
This \nameref{Ideals split} shows that if $J\ideal A''$ is a
\sw{}ly closed ideal then $A''$ splits into two \sw{}ly closed ideals,
$A''=J\oplus (1-j)A''$.  Conversely, if we have a central projection
$j\in A''$, then we obtain an algebra direct sum $A''=jA''\oplus
(1-j)A''$ by the \sw{}ly closed ideals $jA''$ and $(1-j)A''$. Thus
there is a one-to-one correspondence between the central projections
of $A''$ and its \sw{}ly closed ideals. This is exactly the same ideal
structure possessed by \vnas.

The central projections of $A''$ induce a similarly nice decomposition
of $\lat A$. If $p\in A''$ is a central projection and $\V\in\lat A$,
then $p\V\subseteq\V$ and $(1-p)\V\subseteq\V$. This means that
$\V=p\V\oplus(1-p)\V$ is a module decomposition of $\V$, and $\lats\H
A=\lats {p\H}A \oplus \lats {(1-p)\H} A$ is a lattice direct sum
decomposition of $\lats\H A$. Conversely, we will see in \nref{split
  lattice} that if $A$ is a \cra\ and $\lat A$ splits, as a lattice,
into the direct sum of two sublattices, then this split is induced by
a central projection in $A''$.

From spectral theory, the central projections of a \vna\ are always
\sa\ and hence have norm $\leq 1$. Thus, if $A$ is similar to a
\csalg\ the central projections of $A''$ will be uniformly \bded.  On
the other hand, if $A$ is an \oa\ and the central idempotents of $A''$
are bounded, then \nref{Orthogonalise idempotents} implies that there
is a similarity $S$ which makes these idempotents \sa.
It is a happy fact that  \cra{}s\ always have uniformly bounded central
idempotents:
\begin{lemma} \label{Central idempotents are bounded}
\label{Central projections are sa}
Let $A\subseteq\BH$ be a \cra\ with \projconst\ $K$, 
and let $P$ be the set of central
projections of $A''$. Then $P$ is  bounded by $K$. Moreover, there
exists a similarity $S$ of $\H$ which makes all the  central
projections \sa.
\end{lemma}
\proof
Suppose $p\in P$. Then $p\H$ is a submodule of $\H$, so by
\nref{Existence of \projconst{}s for crp} 
there is  a projection $q\in A'$ with
$\norm{q}\le K$ and $q\H=p\H$. Now $q\in A'$, $p\in A''$ and they
share the same range space, so $q=pq=qp=p$, giving 
$\norm{p}\le K$.
The existence of the orthogonalising \sy\ follows from 
\nref{Orthogonalise idempotents}.
\qed

\begin{proposition}\label{DCT on GNS}
  Let $A\subseteq\B(\H')$ be a \tra, and let $\pi:\B(\H')\ra \BH$ be
  the GNS \rep\ for the \csalg\ $\B(\H')$.  Let $\widehat\pi$ denote
  the nondegenerate sub\rep\ $\widehat\pi:A\ra \B(\clos{\pi(A)\H)}$.
  Then $\widehat\pi(A)''=\swclos{\widehat\pi(A)}$ and this algebra is
  \sw--\ws\ homeomorphically isomorphic to $\bd A$.
\end{proposition}
\proof
The point here is that $\pi$ might be a degenerate \rep\ of $A$. 
The properties of the GNS \rep\ mean that $\bd A$ can be 
identified with $\swclos{\pi(A)}$. 
\Nref{Restriction is closed} shows that $\swclos{\pi(A)}$
is \sw\ homeomorphically isomorphic with $\swclos{\widehat\pi(A)}$.
However, $\widehat\pi$ is nondegenerate and
\nref{Bare double commutant} shows that
$\widehat\pi(A)''=\swclos{\widehat\pi(A)}$. 
\qed

\begin{corollary}
  \label{Identity of A**}
  Let $A$ be a \tra, and let $K$ be the \projconst\ function of
  \nref{Bound on \projconst{}}. Then $\bd A$ has identity $e$ with
  $\norm{e}\leq K(1)$.
\end{corollary}
\proof
This follows from \nref{Double commutant} and \nref{DCT on GNS}.
\qed

We may use \nref{DCT on GNS} to apply \nref{Central idempotents are
  bounded} to $\bd A$, concluding that if $A$ is a \tra\ then the set
of central idempotents of $\bd A$ is uniformly bounded.
\begin{defn}
  Let $A$ be a \tra. The supremum of the norms of the central
  idempotents of $\bd A$ is called the \Index{ideal constant} of $A$.
\end{defn}

As another application of these results we can show that closed ideals
of \tras\ have \bais\ (as is the case for \csalgs).

\begin{proposition} \label{Ideals are amenable}\label{Ideals have bai}
Let $A$ be a \tra, and let $J\ideal A$ be a closed two-sided
ideal. 
Then $J$ possesses a \Index{\bai}.
\end{proposition}
\proof
Using \nref{DCT on GNS} to obtain a \rep\ of $A$ with $A''\cong\bd A$, 
\nref{Ideals split} says that the bidual ideal $J^{**}\ideal
A^{**}$ has an identity $j\in J^{**}$. Letting $j_\nu\ra j$ be a
bounded net of elements of $J$ \ws\ convergent to $j$, 
the Mazur theorem 
shows that suitable convex combinations of the
$j_\nu$ will form a \lbai\ for $J$. Similarly we may construct a
\rbai\ for $J$, and hence a \bai, by the observation of Dixon
\cite{DixonTrick}.
\qed
\begin{corollary}
  Let $A$ be an \aoa\ and suppose that $J\ideal A$ is a closed
  two-sided ideal. Then $J$ is amenable.\index{amenable algebra}
\end{corollary}
\proof
It is known that an ideal of an amenable \ba\ is itself amenable \iff\
it possesses a \bai\ \cite{CurtisLoy}.
\qed

Suppose that we have an \oa\ $A$ and a finite set
of closed two-sided ideals $\{J_i\}$ with the property that $\sum
J_i=A$, meaning that for any $a\in A$ there are elements $a_i\in J_i$
with $\sum a_i=a$. 
In the case where $A$ is a \csalg, it is known that the $a_i$ can
be chosen with $\norm{a_i}\leq \norm{a}$ \cite{StormerIdeals}, 
\cite{ElliotOleson}. We show
that a similar fact is true for \tras.
\begin{proposition}\label{Bounded decomposition}
  Let $A$ be a \tra, and $\{J_i\}$ a finite set of closed two-sided
  ideals with $\sum J_i=A$. Let $M$ be the ideal constant
  of $A$. Then for any $a\in A$ and $\epsilon>0$ 
  there are $a_i\in J_i$ with $\sum
  a_i=a$ and $\norm{a_i}\leq (M+\epsilon)\norm{a}$.
\end{proposition}
\proof\def\ov{\expandafter\widehat}%
The $\ell^\infty$~direct sum of the algebras $J_i$ forms a \bs\ which
we denote $\sum^{\ell^\infty} J_i$.  Let $\varphi: \sum^{\ell^\infty} J_i \ra A$ be
given by $\varphi(a_i)=\sum a_i$. By assumption, $\varphi$ maps onto $A$
and hence induces a \bs\ isomorphism $\ov\varphi$ of $\sum^{\ell^\infty} J_i /
\ker \varphi$ onto $A$. It is readily verified that the result to be
proved holds \iff\ the inverse isomorphism ${\ov\varphi}^{-1}$ enjoys
the bound $\norm{\ov\varphi^{-1}}\leq M$. 

We make the following natural isometric identifications:
\begin{eqnarray*}
\bd{\left({\sum}^{\ell^\infty} J_i\right)} & = & {\sum}^{\ell^\infty} \bd {J_i}, \\
\bd{\left({\sum^{\ell^\infty} J_i \over \ker \varphi}\right)}
& = & {\sum^{\ell^\infty} \bd {J_i} \over \bd{(\ker \varphi)}}, \mbox{\rm\quad and} \\
{\sum^{\ell^\infty} \bd{J_i} \over \bd{(\ker \varphi)}} & = & {\sum^{\ell^\infty}
    \bd{J_i}\over \ker(\bd \varphi)}.
\end{eqnarray*}
The last of these identifications is possible because $\varphi$ has
closed range (see the discussion on page~\pageref{Page with closed range}).
These identifications take the bidual isomorphism
$$\bd{\ov\varphi}:\bd{\left({\sum^{\ell^\infty} J_i\over \ker
      \varphi}\right)}\ra \bd A
\hbox{\quad to \quad}
\ov{\bd \varphi}: {\sum^{\ell^\infty} \bd {J_i} \over \ker (\bd\varphi)} \ra
\bd A.
$$
Since $
{\bd{\ov\varphi}}{}^{-1}=
{{{\ov\varphi}^{-1}}}\bd{}
$ and
${{\ov\varphi}^{-1}}\bd{}$ extends ${\ov\varphi}^{-1}$, it
suffices to show that $\norm{{\ov{\bd\varphi}}{}^{-1}}\leq M$.

The ideals $\bd {J_i}$ are \ws-closed two-sided ideals in $\bd A$,
and so by \nref{Ideals split} and \nref{Central idempotents are
  bounded}  they correspond to central idempotents
$j_i\in \bd A$ with $\norm{j_i}\leq M$.
We define another set of central idempotents inductively as follows:
let $k_1=j_1$ and $k_{n+1}=j_{n+1}(1-\sum_{i\leq n} k_i)$. Then
$k_i\in\bd{J_i}$, $k_i k_j = 0$ for $i\neq j$ and $\sum
k_i=\bigvee j_i$.  We put $k=\sum k_i$.  
Since $j_i\leq k$ and $k\in\sum \bd{J_i}$ it follows that 
$\sum \bd{J_i}=k\bd A$.  This is a \ws-closed ideal. If
$\sum\bd{J_i}\neq \bd A$, then by the Hahn-Banach theorem there is
$f\in A^*$ with $f\neq 0$ and $\<f, \bd J_i>=0$ for all $i$, which
implies that $\<f, J_i>=0$ for all $i$ and thus $f=0$. This
contradiction implies that $\sum \bd{J_i}=\bd A$, and that $k=1$.

The set $\{k_i\}$ allows us to decompose any
$a\in \bd A$ into $a=\sum ak_i$, where $ak_i\in \bd J_i$. 
Note that since $k_i$ is a central
idempotent, we have $\norm{k_i}\leq M$, and so this decomposition
demonstrates that ${\ov{\bd\varphi}}{}^{-1}$ is bounded by $M$.
\qed

Similar techniques can be used to establish further properties of the
ideal structure of $A$ which are familiar from \csalg\ theory.
\begin{proposition}\label{4th isomorphism theorem}
  Let $A$ be a \tra, and let $J_1, J_2\ideal A$ be closed two-sided
  ideals. Then the natural homomorphism 
$$
\varphi: {J_1\over J_1\cap
  J_2} \ra {\overline{J_1+J_2}\over J_2}
$$ 
is an onto isomorphism. If the
 ideal constant of $A$ is $M$, then $\norm{\varphi^{-1}}\leq M$.
\end{proposition}
\proof The map $\varphi$ is one-to-one, contractive and dense-ranged.
Consider the bidual map
$$
\bd \varphi: \bd{\left({J_1\over J_1\cap J_2}\right)}\ra 
  \bd{\left({\overline{J_1+J_2}\over J_2}\right)}.$$
Again, we may isometrically identify 
$$\bd{\left({J_1\over J_1\cap
      J_2}\right)}={\bd {J_1} \over \bd{(J_1 \cap J_2)}}$$
 and 
$$ \bd{\left( {\overline{J_1 + J_2} \over J_2}\right)} =
{\bd{(\clos{J_1+J_2})}\over \bd{J_2}}.$$

These identifications allow $\bd \varphi$ to be written as
$$
\bd\varphi: {\bd J_1\over \bd{(J_1\cap J_2)}}\ra
{\bd{\clos{J_1+J_2}}\over \bd J_2}.
$$

The ideals $\bd {J_1}$, $\bd{J_2}$, $\bd {\overline{J_1+J_2}}$ and
$\bd{(J_1\cap J_2)}$ are \ws-closed ideals of $\bd A$, corresponding
to central idempotents. We write $j_i$ for the central idempotents
corresponding to $\bd J_i$, so that $\bd {J_i}=j_i\bd A$.

Observe that $(j_1+j_2-j_1j_2)\bd A \subseteq \bd J_1+\bd J_2
\subseteq \bd{\clos{J_1+J_2}}$. Also, since $j_i \leq j_1+j_2-j_1j_2$
for $i=1,2$ we have $J_1+J_2\subseteq
\bd{J_1}+\bd{J_2}\subseteq(j_1+j_2-j_1j_2)\bd A$, and consequently
$\bd{\clos{J_1+J_2}}=\bd J_1+\bd J_2 = (j_1+j_2-j_1j_2)\bd A$.

On the other hand, we have $\bd{(J_1\cap J_2)}\subseteq \bd J_1\cap
\bd J_2$, and so if $k$ is the central projection corresponding to
$\bd{(J_1\cap J_2)}$ then $k\leq j_1j_2$. 

The splittings of $\bd A$ induced by these central projections allow
us to rewrite the action of $\bd \varphi$ in terms of ideals rather
than quotients. The map $\overline\varphi$ appearing in the bottom
line of the commutative diagram
$$
\begin{CD}
{\displaystyle {\bd {J_1} \over \bd{(J_1\cap J_2)}}} @>\bd\varphi>> {\displaystyle{ \bd{\clos{J_1+J_2}}
  \over \bd{J_2}}} \\
@A\cong AA                                                @V\cong VV \\
j_1(1-k)\bd A @>\overline\varphi>> (j_1+j_2-j_1j_2)\bd A
\end{CD}
$$
is given by 
$$
\overline\varphi:j_1(1-k)a\mapsto(j_1+j_2-j_1j_2)(1-j_2)a.$$
Since $(j_1+j_2-j_1j_2)(1-j_2)=j_1(1-j_1j_2)\leq j_1(1-k)$ it follows
that the range of $\overline\varphi$  is all of
$(j_1+j_2-j_1j_2)(1-j_2)\bd A$ and so is \ws\ closed. This means 
that the range of $\bd \varphi$ is \ws\ closed, and then
a standard theorem of functional analysis \cite[corollary
8.7.4]{Edwards} implies that the range of $\varphi$ is norm closed in
$\clos{J_1+J_2}/J_2$. Thus $\varphi$ is an isomorphism onto
$\clos{J_1+J_2}/J_2$. 

Since $\varphi$ is an isomorphism, $\bd \varphi$ and hence
$\overline\varphi$ are also isomorphisms. Consequently $k=j_1j_2$
and $\overline\varphi$ is the identity map. The norm of $\inv{(\bd
  \varphi)}$ can be calculated from the norms of the splitting maps.
The map from 
$\bd{\clos{{J_1}+{J_2}}}/\bd{J_2}\ra(j_1+j_2-j_1j_2)\bd A$
has norm at most $\norm{j_1+j_2-j_1j_2}\leq M$. The map from
$j_1(1-k)\bd A$ to 
${\bd{J_1}}/{\bd{(J_1\cap J_2)}}$ is contractive. Hence
$\norm{\inv\varphi}\leq\norm{\inv{\bd\varphi}}\leq M$.
\qed

The following result is elementary,
but a proof is included here for completeness.
\begin{lemma}
\label{Elementary normed space result}
  Let $0\ra X\stackrel{f}\ra Y\stackrel{g}\ra Y/X\ra 0$ 
  be an exact sequence of normed spaces, with $f$ an isometry and $g$
  the natural quotient map. If $X$ and $Y/X$ are complete, then so is $Y$.
\end{lemma}
\proof
Let $\{y_i\}$ be a sequence in $Y$ with $\sum\norm{y_i}< \infty.$ Then
the image sequence $\{y_i+X\}$ is summable in $X/Y$ with sum
$y+X$. Denote by $s_n$ the partial sum $\sum_{i\leq n}y_i$. Then
$\norm{s_n - y + X}\ra 0$. Thus, there is a subsequence $\net{y_{k_i}}$ and a
corresponding sequence of elements $x_i\in X$ with 
$\norm{s_{k_i}-y-x_i}< 2^{-i}$ for all $i$.
Since the sequence $\{s_{k_i}\}$ is Cauchy, the sequence $\{x_i\}$ is
Cauchy and converges in $X$ to $x$, say. Then it is readily verified
that $s_i-y-x\ra 0$, and so $Y$ is complete.
\qed

\begin{proposition}\label{Sum is closed}
  Let $A$ be a \tra, and let $J_1, J_2\ideal A$ be closed two-sided
  ideals. Then $J_1+J_2$ is closed.
\end{proposition}
\proof
As associative algebras, we have the natural algebraic isomorphism
$$
{J_1\over {J_1\cap J_2}} \cong {J_1+J_2\over J_2},
$$
where $J_1+J_2$ denotes the algebraic sum of $J_1$ and $J_2$ in $A$ (\ie, no
closure is taken).
\Nref{4th isomorphism theorem} then implies that $(J_1+J_2)/J_2$ is
a closed subspace of $A/J_2$. 
Consider the natural short exact sequence of normed spaces:
$$
0\ra J_2\ra J_1+J_2\ra {J_1+J_2\over J_2}\ra 0.
$$
Since the first and last terms are complete normed spaces, 
\nref{Elementary normed space result} shows that $J_1+J_2$ is complete
as a normed space. That is, $J_1+J_2$ is a closed subspace of $A$.
\qed

\begin{lemma}\label{Intersection is product}
  Let $A$ be a \tra, and let$J_1, J_2\ideal
  A$ be closed two-sided ideals. Then $J_1\cap J_2=J_1 J_2$.
\end{lemma}
\proof
Since the $J_i$ are two-sided ideals, we have $J_1 J_2 \subseteq
J_1\cap J_2$. For the converse, recall that \nref{Ideals have bai}
shows that $J_1$ contains a \bai, $\{e_\nu\}$ say. By Cohen's
factorisation theorem we know that $J_1J_2=\clos{J_1J_2}$. Suppose
that $j\in J_1\cap J_2$. Then $e_\nu j\in J_1J_2$ and $e_\nu j\ra j$,
so $j\in J_1J_2$.
\qed
\begin{lemma}\label{pairwise disjoint implies disjoint}
  Suppose that $A$ is a \tra\ and let $K\ideal A$ be a closed
  two-sided ideal.
  Let $\{J_i\}$ be a finite collection of  two-sided ideals
  with
  $J_i K=\{0\}$ for $i\neq j$. Then $\left(\sum J_i\right) \cap K=\{0\}$. 
\end{lemma}
\proof
Since the sum $\sum J_i$ is the algebraic sum, we have 
$\left(\sum J_i\right) \cap K = \left(\sum J_i\right) \cdot K 
= \sum (J_i\cdot K) = \{0\}$.
\qed
\begin{proposition}\label{l-infty norms}
  Let $A$ be a \tra\ with ideal constant $M$
  and suppose that $\{J_i\}$ is a finite collection
  of closed two-sided ideals with $J_i\cap J_j=\{0\}$ for $i\neq j$. 
  Then if $a_i\in J_i$ for all $i$, we have $\norm{\sum a_i}\leq
  (1+2M)\sup_i\{\norm{a_i}\}$.
\end{proposition}
\proof Using \nref{Sum is closed} and \nref{pairwise disjoint implies
  disjoint} an inductive argument shows that $\sum
J_i\cong\sum^{\ell^\infty} J_i$.  Applying the bidual functor yields
$\bd{(\sum J_i)}=\sum \bd{J_i}\cong \sum^{\ell^\infty}\bd {J_i}$,
where $\sum \bd{J_i}$ denotes the sum inside $\bd A$.  Using
\nref{Ideals split} we obtain a collection of central idempotents
$\{j_i\}$ with $j_i\bd A=\bd{J_i}$. The isomorphism
$\sum^{\ell^\infty}\bd {J_i} \cong \sum \bd {J_i}$ implies that
$j_ij_k = 0$ for $i\neq k$. Using the GNS construction to represent
$\bd A$ isometrically as a concrete \oa, \nref{Orthogonalise
  idempotents} implies these idempotents can be conjugated to
self-adjoint projections by a similarity $S$ with $\norm{S}\norm{\inv
  S}\leq 1+2M$.  Once the idempotents are self-adjoint the addition
map $\sum^{\ell^\infty} \bd{J_i}\ra \bd A$, which extends the map
$\sum^{\ell^\infty} J_i\ra A$, is contractive; the similarity $S$
degrades this norm to $1+2M$.  \qed

\dealwithsectionbreaks
\section{Stability of the \trp}
\label{Section with stability}

It is interesting to consider the stability of the \rp\
and its relations. 
The
results of this section show that 
the \trp\ is fairly stable, in the sense that
many of the standard operations for
creating new algebras from old preserve the \trp.
The main interest is in the stability of the \trp, because
it is the only one of the three \rps\ which is isomorphism invariant.

It is satisfying that the category of \oas\ enjoys a certain amount of
stability itself. In particular, the theory of \oas\ as developed by
Effros and Ruan \cite{EffrosRuan} shows that quotients of \oas\ are
again \oas.

We start with the simple observation that the \trp\ is stable under homomorphisms.

\begin{proposition}\label{Image has cp}
  Let $A$ be an \oa\ with the \trp, 
  and let $\theta:A\ra \BH$ be a \rep. Then the norm closure 
  $\clos{\theta(A)}$
  has the \trp.
\end{proposition}
\proof
Denote by $B$ the closure of $\theta(A)$, and let $\varphi:B\ra \B(\H)$ be
a \rep. Then $\varphi\theta$ is a \rep\ of $A$. Moreover,
since $\theta(A)$ is dense in $B$ it follows that $\lat \varphi\theta(A)$
coincides with $\lat \varphi(B)$. Since $A$ has the \trp\ $\lat
\varphi\theta(A)$ is topologically complemented, so $\lat \varphi(B)$
is topologically complemented, and
$B$ has the \trp.
\qed
\begin{corollary}\label{Quotient has cp}
  Let $A$ be an \oa\ with the \trp, and let $B$ be a quotient of $A$.
  Then $B$ has the \trp.  \qed
\end{corollary}

Recall from the example on page~\pageref{Page with RP example} that
neither the \rp\ nor the \crp\ are stable under
quotients. Furthermore, if
$A\subseteq\BH$ is a \cra\ and $\theta:A\ra\B(\H')$ is a \rep\ such that
$\H'$ does not have the \rp, then $\id_A\oplus\,\theta$ is a faithful
\rep\ and $\H\oplus\H'$ does not have the \rp. Consequently neither
the \rp\ nor the \crp\ are \oa\ isomorphism invariants.

\begin{proposition}\label{TRP is hereditary}\label{CP is hereditary}
  Let $A$ be an \oa\ with the \trp\ and let $J\ideal A$ be a closed
  two-sided ideal. Then $J$ has the \trp.
\end{proposition}
\proof Let $\theta: J\ra \BH$ be an arbitrary \rep\ of $J$. We need to
show that $\H$ has the \rp\ as a $J$-module.

\Nref{Ideals split} shows that there is a central projection $j\in \bd A$ such
that $\bd J=j\bd A$.
By the universal property of the bidual the homomorphism $\theta$ may
be extended to a \ws\ \cts\ \rep\ $\bar\theta:\bd J\ra \BH$.  We
define a representation $\widehat\theta: A\ra \BH$ by $a\mapsto
\bar\theta(ja)$. Since for $a\in J$ we have $ja=a$, $\widehat\theta$
is an extension of the \rep\ $\theta$ to $A$. Now suppose
$\V\subseteq\H$ is $\theta(J)$-invariant.  By the \ws--\sw\ continuity
of $\bar\theta$, $\V$ is $\bar\theta(\bd J)$-invariant, hence
$\widehat\theta(A)$-invariant. Since $A$ has the \trp, there is a
complementing $\widehat\theta(A)$-invariant subspace $\W$. As
$\widehat\theta$ extends $\theta$, it follows that $\W$ is
$\theta(J)$-invariant.  Hence $\H$ has the \rp\ as a $J$-module, and
$J$ has the \trp.  \qed

A corollary of this result is that the radical of a \tra\ is also a
\tra.  This contrasts with the case of amenable \bas, where it is
rather more difficult to find radical examples than it is to find
non-\ss\ examples. For instance, if $G$ is a locally compact abelian
group, and $J\ideal L^1(G)$ is an ideal with $J\varsubsetneq \ker\hul J$, then
$L^1(G)/J$ is amenable and non-semisimple. However, no amenable radical
abelian \ba\ is known \cite{GronbaekWillis}.

On the other hand, no non-semisimple \tras\ are known to the author.
In fact, since \csalgs\ are always semisimple, if \nref{tra
  conjecture} is correct then every \tra\ is semisimple.

\comment{
\begin{lemma}
  Suppose that the diagram 
  $$0\ra A\ra C\ra B\ra 0$$ 
  is a short exact sequence of \oas, 
  such that both $A$ and $B$ have the \trp. Then
  $C$ is an \oa\ and has the \trp.
\end{lemma}
\proof Consider the bidual diagram 
$$
0\ra \bd A\ra \bd C\ra \bd B\ra 0,
$$ 
where $\bd A$ is a \ws\ closed ideal in the \oa\ $\bd C$.  We
represent $C$ via the universal representation; that is, we realise
$C\subseteq\BH$ such that $\clos{C}=\bd C$. Consequently, by
\nref{Bare double commutant}, $A\subseteq\BH$ is \sw\ dense in $\bd
A$.  The subspace $\clos{\bd A\H}$ is $A$-invariant, and so there is a
projection $p\in (\bd A)'$ onto $\clos{\bd A\H}$. It follows that $\bd
A$ must annihilate the kernel of $p$, and that algebra $\bd
A|_{\ssclos{\bd A\H}}$ is \sw\ closed. Applying \nref{Bare Double
  Commutant} again, we find that there is a projection $e\in \bd A$
which is an identity for $\bd A$, with range $\clos{\bd A\H}$. Since
$\bd A$ is a two-sided ideal, $e$ is central in $\bd C$. This central
projection $e$ decomposes $\bd C$ into the direct sum of two ideals,
which are isomorphic to $\bd A$ and $\bd B$ respectively. We write
$\bd C=\bd A\oplus\bd B$.  The ideal $\bd A$ acts on the space $e\H$,
and $\bd B$ acts on $(1-e)\H$.

Now suppose that we have a \rep\ $\theta:C\ra \B(\L)$. This \rep\
extends to a \rep\ $\bar\theta:\bd C\ra \B(\L)$. Then the
projection $e\in \bd C$ maps under $\bar\theta$ to some
projection $\bar e$ 
on $\L$. If $\K\subseteq\L$ is a $C$-invariant subspace,
then it is invariant under $\bar e$ and $1-\bar e$. 
For any  element $\xi\in\K$, 
both $\bar e\xi$ and $(1-\bar e)\xi$ are in $\K$,  and so $\K$
splits as the direct sum of two subspaces, $\bar e\K$ and $(1-\bar
e)\K$.
The subspace $\bar e\K$ is $A$-invariant, and so is
complemented in $\bar e\H$ by an $A$-module, hence an $\bd
A$-module, hence a $\bd C$-module. Similarly, there is a $\bd
C$-module complement to $(1-\bar e)\K$ in $(1-\bar e)\H$. 
The direct
sum of these two modules is the desired complement to $\K$.
\qed
}
\begin{proposition}
  Suppose that the diagram $$0\ra A\ra C\ra B\ra 0$$ is a short exact
  sequence of \bas, such that both $A$ and $B$ are \oa{}s with the
  \trp. Then $C$ is an \oa\ and has the \trp.
\end{proposition}
\proof
Consider the bidual \bs\ diagram $$0\ra \bd A\ra \bd C\ra \bd B\ra 0.$$
We equip $\bd C$ with the left Arens product temporarily, until
we show that $C$ is Arens regular.

The left Arens product has the property that $x\mapsto xy$ and
$x\mapsto zx$ are \ws\ \cts\ maps on $\bd C$ for all $y\in \bd C$ and
$z\in C$. If $a\in \bd A$ and $c\in \bd C$, there are nets
$\net{a_\nu}\subseteq A$ and $\net{c_\mu}\subseteq C$ with $a_\nu\ra
a$ and $c_\mu\ra c$ in the respective \ws\ topologies. Then
$ac=\lim_\nu a_\nu c=\lim_\nu\lim_\mu a_\nu c_\mu\in \bd A$ and
$ca=\lim_\mu c_\mu a=\lim_\mu\lim_\nu c_\mu a_\nu\in \bd A$, showing
that $\bd A$ is a \ws\ closed ideal in $\bd C$. As $\bd A$ contains an
identity $e$ say, it follows that $e$ is central in $\bd C$.  Then
$\bd C=e\bd C \oplus (1-e)\bd C$ is an algebra decomposition of $\bd
C$. We have $e\bd C\cong\bd A$ and $(1-e)\bd C\cong \bd C/\bd
A\cong\bd B$, and so up to \ws\ homeomorphic isomorphism $\bd C\cong \bd
A\oplus \bd B$ as a \bs. Checking the definition of the Arens product
shows that this is also an algebraic direct sum. Hence $\bd C$ is an
\oa.  This implies that $C$ itself is an \oa\ and hence is Arens
regular.

Suppose now that we have a \rep\ $\theta:C\ra \B(\H)$. This 
extends to a \rep\ $\bar\theta:\bd C\ra \B(\H)$. Then the projection
$e\in \bd C$ maps under $\bar\theta$ to some projection $\bar e$ on
$\H$. If $\V\subseteq\H$ is a $C$-invariant subspace, then it is
invariant under $\bar e$ and $1-\bar e$.  Thus $\V$ splits
as the direct sum of the two subspaces $\bar e\V$ and $(1-\bar e)\V$.
The subspace $\bar e\V$ is $A$-invariant, and so is complemented in
$\bar e\H$ by an $A$-module, hence an $\bd A$-module, hence a $\bd
C$-module. Similarly, using the \trp\ for $B$ there is a $\bd
C$-module complement to $(1-\bar e)\V$ in $(1-\bar e)\H$.  The direct
sum of these two modules is the desired complement to $\V$.  \qed
So far we have shown that ideals, quotients and extensions preserve the \trp.
Another standard tool for constructing algebras is the tensor
product. 
The tensor product is generally tricky to define for operator
algebras---it is not too clear how to norm the algebraic tensor
product in a suitable way. 
However, when one of the \oas\ is
$M_n$ there is an elementary notion of the tensor product. 
If $A\subseteq\BH$ is an \oa, 
then the tensor
product $A\ot M_n$ can be considered as the algebra of $A$-valued
$n\times n$ matrices, given the operator norm derived from
the natural action on $\pma\H n$.
\begin{proposition}
  Let $A$ be an \oa\ with the \trp. Then the tensor product $A\ot M_n$
  has the \trp.
\end{proposition}
\proof
We consider first the case where $A$ is unital. This means that
we can consider both $A$ and $M_n$ to be embedded in $A\ot M_n$ via
$A=A\ot 1$ and $M_n=1 \ot M_n$. Suppose that
$\theta:A\ra \BH$ is a \rep.  Then restricting to the elements of the
form $1\ot T\in 1\ot M_n$ gives a \rep\ of $M_n$. The idempotents
$e_{ii}\in M_n$ map to a collection of commuting idempotents
$q_i=\theta(e_{ii})$. We put $q_0=1-\sum q_i$, and $\H_i=q_i\H$.

If $\V\in\lat A\ot M_n$, then $\V=\sum^\oplus q_i \V$. Note that the
subspaces $q_i \V$ are not $A\ot M_n$-modules. They are, however,
invariant under $A\ot 1$, and so are $A$-modules. Taking $q_1\V$ for
definiteness, there is an $A$-module complement $\W_1\subseteq\H_1$.
We set $\W_j=\theta(e_{j1})\W_1$. To deal with the degenerate part of
the \rep, we put $\W_0=\H_0\ominus q_0\V$. Then $\W=W_0\oplus
\sum^\oplus_{i\leq n} \W_i$ is an $A\ot M_n$-module complement
to $\V$.

When $A$ has no identity, we can lift the \rep\ to the bidual $\bd {(A\ot
  M_n)}=\bd A\ot M_n$. The invariant subspaces of $\H$ are not
changed by this transition. Since $\bd A$ is unital, the result
follows from the above argument and \nref{Bicommutant has rp}.
\qed

A simple but useful result is that the addition of a identity
does not affect the \trp. 
\begin{proposition}
  \label{Unitisation lemma}
Let $A\subseteq\BH$ be an \oa, and denote by $A_+$ the \oa\ generated
by $A$ and $1_\H$. Then $A_+$ has the \trp\ \iff\ $A$ does.
\end{proposition}
\proof Since $A$ is a two-sided ideal in $A_+$, \nref{TRP is
  hereditary} gives the `only if' direction. For the converse, suppose
that $A$ has the \trp\ and $\theta:A_+\ra\B(\H')$ is a \rep\ of $A_+$.
By \nref{Ideals have bai} $A$ has a \bai. Let $\bar e$ denote the
projection of \nref{bai lemma}.  The identity $1_\H\in A_+$ maps to
an idempotent $f\in\B(\H')$, with $\bar e\leq f$.

If $\V\in\lat A_+$, then $\H=\bar e \H \oplus f(1-\bar e)\H \oplus
(1-f)\V$ is a central decomposition of $\H$, which induces a central
decomposition of $\V$.
Since $A$ is a \tra, there is an
$A$-submodule $\W\subseteq\bar e\H$ with $\bar e\V \oplus \W=\bar
e\H$, which is $A_+$-invariant because
$\bar e\leq f$. Then the subspace 
$$
\W \oplus [f(1-\bar e)\H\ominus f(1-\bar e)\V] 
  \oplus [(1-f)\H \ominus (1-f)\V].
$$ is an $A_+$-submodule complement to $\V$.  \qed

The last\galley{careful} construction to investigate is the direct sum.  If
$\{A_\lambda\}$ is a family of \tras\ the direct sum
$\sum^{c_0} A_\lambda$ need not have the
\trp. For instance, let $A_n$ be the matrix algebra
$$
A_n=\left\{\TxT a, n(b-a), 0, b; :  a,b\in \C\right\}.
$$
Each $A_n$ is similar to $\C^2$, and so has the \trp. However, the
natural \rep\ of $A_n$ on $\C^2$ gives $\lat A_n=\{0, \C^2,
\C e_1, \C(ne_1-e_2)\}$, and so the only $A_n$-module projection onto $\C e_1$
is 
$$
p_n=\TxT 1, n, 0, 0;.
$$
As $n\ra \infty$ we have $\norm{p_n}\ra \infty$. Then the natural
contractive \rep\ of $\sum^{c_0}A_n$ on
$\sum^{\oplus}\C^2$ has the invariant subspace $\sum^{\oplus}\C e_1$, which
cannot be the range of a bounded module projection.

The problem is that the algebras $A_n$ do not have the \trp\
`uniformly'. We can formalise this notion by referring to the
\projconst{} function $K:\R^+\ra \R^+$ defined in \nref{Bound
  on \projconst{}}. 
\begin{proposition}
  \label{Direct sum is tra}
  Let $\{A_\lambda\}_{\lambda\in \Lambda}$ be a family of \tra{}s,
  with corresponding \projconst{} functions $K_\lambda$. Define
  $K(\alpha)=\sup_\lambda\{K_\lambda(\alpha)\}$.  Then
  $\sum^{c_0}A_\lambda$ has the \trp\ \iff\ $K(\alpha)<\infty$ for all
  $\alpha\in\R^+$.
\end{proposition}
\proof Suppose first that $K(\alpha)<\infty$ for all $\alpha\in\R^+$,
and define $A=\sum^{c_0}A_\lambda$. The bidual of $A$ can be
identified as $\bd A=\sum^{\ell^\infty}\bd {A_\lambda}$. Each $\bd A$
contains an identity $e_\lambda$ with norm bounded by the ideal
constant of $\bd {A_\lambda{}}$, which is in turn uniformly bounded
over $\lambda$ by assumption.  Suppose $\theta:A\ra\BH$ is a \rep. We
may extend this to a \rep\ of the bidual algebra $\bd
A=\sum^{\ell^\infty} \bd{A_\lambda}$.  The Boolean algebra of
idempotents in $\bd A$ generated by the $e_\lambda$ is uniformly
bounded, and so there is a similarity $S$ on $\H$ such that
$\theta^S(e_\lambda)$ is \sa\ for all $\lambda$. We may thus assume
without loss of generality that $\theta(e_\lambda)$ is \sa\ for all
$\lambda$.

If $a=(a_\lambda)\in A$ then $a=\sum
e_\lambda a$ and $\theta(a)=\sum \theta(e_\lambda)\theta(a)$, with
norm convergence of the sums.  Also, $\sum
\theta(e_\lambda)$ converges in the strong topology of $\BH$.  For
$\V\in \lat A$, define $\V_\lambda=\theta(e_\lambda)\V$. Put
$f_0=1-\sum \theta(e_\lambda)$, and $\V_0=f_0\V$. Then each
$\V_\lambda$ is an $A_\lambda$-module if we consider $A_\lambda$ as a
subalgebra of $A$ in the natural way, and so $\V_\lambda$ is
complemented in $\theta(e_\lambda)\H$ by an $A_\lambda$-module such
that the norm of corresponding projection $p_\lambda$ onto
$\V_\lambda$ is bounded by $K(\norm{\theta})$. We may choose the
orthogonal projection $p_0$ from $f_0\H$ to $\V_0$ since $A$ annihilates
$f_0\H$. Then the sum $p_0f_0+ \sum p_\lambda e_\lambda $ converges
in the strong topology to an $A$-module projection onto $\V$.

Conversely, 
if $K(\alpha)=\infty$ for some $\alpha\in\R^+$, then a simple direct
sum argument shows that $\sum^{c_0}A_\lambda$ does not have the \trp.
\qed

Is is an interesting question whether the \trp\ is stable under
$\ell^\infty$~direct sums. The difficulty is that there can be a
non-trivial module action `at infinity', and so the lattice structure
of $\sum^{\ell^\infty} A_\lambda$ cannot be decomposed into individual
$A_\lambda$-modules so easily. On the other hand, if \nref{tra
  conjecture} is correct, then \nref{Uniform bound on similarities}
implies that a result analogous to \nref{Direct sum is tra} is true
for $\ell^\infty$~direct sums.
\galley{comment on amenability results here?}

\dealwithsectionbreaks
\section[Calculations with {$2\times 2$ matrices}]{Calculations with
  {$\mathbf {2\times 2}$ matrices}}
\label{Section with twine}

Let $A\subseteq\BH$ be a \cra\ acting on \hs\ $\H$. If $\v$ and $\w$
are two submodules of $\H$, the set of module maps from $\v$ to $\w$
is a linear space. We define a relation $\twine$ on the set of
submodules of $\H$ as follows: $\v\twine \w$ \iff\ there is a non-zero
module map from $\v$ to $\w$. In accordance with standard usage we say
that such a map intertwines $\V$ and $\W$.

If $A$ is a \csalg\ and both $\V$ and $\W$ are obtained from \sreps,
it is well-known that $\V\twine W$ \iff\ $\W\twine V$ (\ie\ the
relation $\twine$ is symmetric).  In fact, when $A$ is a \csalg\ and
$T:\V\ra\W$ is an $A$-module map with polar decomposition $T=US$,
then $U$ is a non-trivial partially isometric module map, and so
establishes an isometric module isomorphism between certain submodules
of $\V$ and $\W$.  This fact provides the starting point for the
analysis of factor \vnas\ (as pursued in \cite{Schwartz}, for
example).

For more general
algebras is is not always true that $\twine$ is symmetric. For instance,
if $A$ is the matrix algebra of section~\ref{Section with triangular
  algebra} with elements of the form
$$
\TxT *, *, 0, *;
$$
acting on $\C^2$ then $\C e_1\twine \C^2$, but it is easy to check
that $\C^2 \not\twine\C e_1$. A less trivial example is afforded by the
\RA\ of section~\ref{rp does not imply crp}.  Using the notation of
that example, the characterisation of the submodule lattice
(\nref{Lattice is what it is}) shows that while $T$ is a module map
from $\H_1$ to $\H_2$, there are no non-zero module maps from $\H_2$
to $\H_1$.

\comment{
Another example which is perhaps less contrived is the disc algebra
$A(D)$. This algebra has a natural contractive \rep\ on the Hardy
space $H^2$, given by pointwise multiplication. On the other
hand, the one-dimensional \rep\ $f\mapsto f(0)$ makes $\C$ into an
$A(D)$-module. The map $T:H^2\ra \C$ given by $\xi\mapsto \int_0^{2\pi}
\xi(e^{it})\,dt$ is a non-zero module map, but there is no non-zero
module map from $\C$ into $\H^2$, since the image of $\C$ under such a
map must lie in $\ker z=\{0\}$. }

It is an important fact that 
the \crp\ rules out this behaviour: for \cra{}s the relation $\twine$
is symmetric.
\begin{proposition}
  \label{Intertwiner relation is symmetric}
  \label{Relation is symmetric}
Let $A\subseteq\BH$ be a \cra, and suppose $\v_1, \v_2$ are submodules
of $\H$ with
$\v_1\twine \v_2$. Then $\v_2\twine \v_1$.
\end{proposition}
\proof We may assume that $\V_1\cap \V_2=\{0\}$ by replacing $\H$ with
$\H\ot_2 \C^2$ if necessary and replacing $\V_1$ with $\V_1\ot_2 e_1$
and $\V_2$ with $\V_2\ot_2 e_2$. This device also ensures that
$\V_1+\V_2$ is closed.

Consider the subrepresentation
$\theta$ of $A$ on
$\v_1\oplus \v_2$. Supposing that $\v_2\not\twine \v_1$,
the elements of $\theta(A)'$  have the matrix form 
$$ \TbT b_{11}, 0, b_{21}, b_{22};,$$ where $b_{ij}$ intertwines
$\v_j\ra \v_i$.  Since $\V_1\twine\V_2$ there is a non-zero module map
$T:\v_1\ra \v_2$. For $\lambda\in\C$ consider the matrix
$$\TbT 1, 0, \lambda T, 0;\in \theta(A)'.$$ 
This is a projection
onto $\Gr {\lambda T}$ with kernel $\v_2$, and it is
the only projection onto $\Gr {\lambda T}$  
of the form given for elements of
$\theta(A)'$. However, as $\lambda\ra
\infty$ the angle between $\Gr \lambda T$ and $\v_2$ goes to
zero. Then \nref{Existence of \projconst{}s for crp} shows that $A$
cannot have the \crp.
\qed

The next few results give information about the algebra generated by a
\cra\ and its commutant. 

\begin{lemma}\label{Complements of joint submodules are uniquely complemented}
  Let $A\subseteq\BH$ be any \oa.
  If $\v\in \lat A$ is complemented by an element $\w\in \lat A \cap
  \lat A'$, then $\w$ is the {\em unique} complement to $\v$ in $\lat A$.
\end{lemma}
\proof
The decomposition $\H=\v\oplus \w$ induces the matrix form 
$$\TbT
b_{11},0,b_{21},b_{22};$$
 on the elements of $A'$. However, any projection in $A'$ onto $\v$ is
necessarily of the form 
$$\TbT 1,b_{12}, 0,0;.$$ 
The only matrix which
satisfies both conditions is the projection along $\w$.
\qed

\Nref{Complements of joint submodules are uniquely complemented} was
noted in \cite{Fong}. The next result also appeared in \cite{Fong} for
algebras with only the \rp, but under the strong assumption that all
graph transformations of $A$ are bounded. Here we are able to use
\nref{Intertwiner relation is symmetric} to replace this graph
transformation assumption with the \crp.

\begin{lemma}\label{Uniquely complemented submodules are central}
\label{Unique complement is central}
  Let $A\subseteq\BH$ be a \cra.
  If $\v\in \lat A$ is uniquely complemented in $\lat A$ by $\w$ then
  the projection onto $\v$ along $\w$ is central in $A''$.
\end{lemma}
\proof
Again, we let $\H=\v\oplus \w$ induce a $2\times 2$ matrix form on
$A'$. The projections onto  $\v$ in $A'$ are exactly those matrices 
$$\TbT
1,b_{12},0,0;$$ 
where $b_{12}$ intertwines $\w$ and $\v$. By the
uniqueness of $\w$, this implies that there are no non-zero
intertwiners from $\w$ to $\v$. Now \nref{Intertwiner relation is symmetric}
shows that there are no non-zero intertwiners from $\v$ to $\w$, and
hence elements of $A'$ are of the form 
$$\TbT b_{11}, 0, 0,
b_{22};.$$ 
Since the projection onto $\V$ along $\W$ commutes with such matrices,
this projection is central.
\qed

The example of section~\ref{rp does not imply crp} shows that
the \rp\ alone is insufficient for the conclusion of \nref{Uniquely
  complemented submodules are central}.

\begin{lemma}\label{Joint submodules are uniquely complemented}
  Let $A\subseteq\BH$ be a \cra.
  If $\v\in \lat A\cap \lat A'$, then $\v$ is uniquely complemented in
  $\lat A$, by $\w\in \lat A\cap \lat A'$.
\end{lemma}
\proof Let $\w\in \lat A$ be a complement to $\v$. By
\nref{Complements of joint submodules are uniquely complemented}, $\v$
is the unique complement to $\w$, and hence by \nref{Uniquely
  complemented submodules are central}, $\w\in \lat A\cap \lat A'$ and
$\w$ is the unique complement to
$\v$.  \qed
\begin{proposition}\label{Commutant has RA}
  Let $A\subseteq\BH$ be a \cra. Denote by $B$ the closed algebra generated
  by $A$ and $A'$. Then $\H$ has the \rp\ as a
  $B$-module, and invariant subspaces of $B$ correspond to central
  projections of $A''$.
\end{proposition}
\proof
We have $\lat B=\lat A\cap \lat A'$. \Nref{Joint submodules are
  uniquely complemented} now gives the result.
\qed

\begin{proposition}\label{Product has RP}
  If $A\subseteq\BH$ is a nondegenerate \cra, and $A\cdot A'$ denotes
  the algebra generated by products of elements of $A$ and $A'$, then
  $\H$ has the \rp\ as an $A\cdot A'$-module.
\end{proposition}
\proof
By 
\nref{Joint submodules are uniquely complemented} we need only
show that 
$\lat A\cdot A'=\lat A\cap \lat A'$. It is clear that
$\lat A\cap\lat A'\subseteq\lat A\cdot A'$. Conversely, since
$A'$ contains the identity operator, it is clear that $ \lat
A\cdot A'\subseteq\lat A$. To finish it suffices to show that 
$\lat A\cdot A'\subseteq\lat A'$. 

Taking $\V\in\lat A\cdot A'$, the \rp\ gives $\W\in\lat
A$ with $\V\oplus\W=\H$. The nondegeneracy assumption means
that $\H=\clos{A\H}\subseteq\clos{A\V}\oplus\clos{A\W}$, and so
$\clos{A\V}=\V$. Since $A'(A\V)\subseteq\V$, it follows that
$A'V=A'\clos{A\V}\subseteq\V$. 
\qed

\Nref{Commutant has RA} and \nref{Product has RP} are interesting
because the central projections of $A''$ are very well
behaved. In particular, \nref{Central projections are sa} shows that
we may assume that the central projections of $A''$ are \sa; in this
case the algebra generated by $A$ and $A'$ is a reductive algebra.

The results of section~\ref{Section with algebraic properties} demonstrate
the equivalence between central projections of $A''$ and \sw{}ly
closed ideals of $A''$. There is a similar connection with direct sum
decompositions of $\lat A$. 
\begin{proposition}\label{split lattice}
  Let $A\subseteq\BH$ be a \cra. Then there are abstract lattices $L_1,
  L_2$ such that $\lat A$ is lattice isomorphic to $L_1 \oplus
  L_2$ \iff\ there is a central projection $p\in A''$ with
  $\lats{p\H}A\cong L_1$ and $\lats{(1-p)\H}A\cong L_2$.
\end{proposition}
\proof The `if' direction is immediate.  For the `only if' direction,
suppose that $\lat A\cong L_1\oplus L_2$.  The submodule $\H$ has a
direct sum decomposition $\H\cong 1_{L_1}\oplus 1_{L_2}$. Then
$1_{L_1}$ is the unique lattice-complement to $1_{L_2}$ in $L_1\oplus
L_2$. Since a pair of complementary submodules in $\lat A$ are
certainly lattice-complements, this means that $1_{L_1}\in\lat A$ is
uniquely complemented (in the sense of submodules), and is hence the
range of a central projection in $A''$. 
\qed


\galley{This section ends rather suddenly. Check how it reads
  later. It doesn't have much text in it either.}

\chapter{Are \RA{}s similar to \csalgs?}
\label{Chapter with characterisations}

\section{Introduction}
\label{Section with reductive algebras}

The results of the previous chapter indicate that \oas\ with either
the \crp\ or the \trp\ share a certain amount of structure with the
class of \csalgs. This structure can be deployed in an attempt to
prove \nameref{cra conjecture}{}s~\ref{cra conjecture} and \ref{tra
  conjecture}.  

\comment{ Having established some of the behaviour of
  algebras with the various reduction properties, we now turn to the
  problem of finding more concrete \reps\ of these algebras. The
  general idea is that by utilising the \rp\ to `orthogonalise' the
  algebra it should be possible to establish some degree of \sa{}ness.
  The exact nature of this \sa{}ness will depend on which \rp\ is at
  hand: we know from the example of section~\ref{rp does not imply
    crp} that the \rp\ alone is generally too weak to produce any
  decent \sa{}ness; \nameref{cra conjecture}s~\ref{cra conjecture} and
  \ref{tra conjecture} indicate what we expect from the \crp\ and the
  \trp\ respectively.  }

Suppose that $A$ is an \oa\ with either the \crp\ or the \trp.
Any approach to the problem of `orthogonalising' $A$
can be broken into several stages.
Firstly, a suitable faithful \rep\ of $A$ as a concrete \oa\ needs to be
chosen, so that we may consider $A$ as a subalgebra of $\BH$.
Then there are two steps to complete.
\begin{enumerate}\label{Page with program}
\item Find a similarity $S\in\BH$ so that $\lat A^S$ is orthogonally
  complemented ($S$ `orthogonalises' $\lat A$).
\item Show that because of the appropriate choice of the \rep\ of $A$,
  the orthogonal complementation in $\lat A$ implies that $A$ is \sa.
\end{enumerate}
Both steps (i) and (ii) in this procedure raise interesting problems
in general.

The first step can be considered as a module analogue of the
complemented subspaces theorem. Recall that this theorem states: if
$X$ is a \bs\ such that every subspace of $X$ is complemented, then
$X$ is \isoc\ to a \hs. In other words, if the lattice $L$ of closed
subspaces of $X$ is topologically complemented then there is an
equivalent inner product norm on $X$ such that $L$ becomes
orthogonally complemented.
The strength of the analogy with step (i) is confirmed in
section~\ref{section with compact case}, where part of the proof of
the complemented subspaces theorem is adapted for exactly this
purpose.

The second step in the procedure is closely related to
the reductive algebra problem (\nref{reductive algebra problem}).
If step (i) in the above program is accomplished, the \sw\ closure of
the algebra $A^S$ will be a reductive algebra.  For $A^S$ to be \sa\ 
it is necessary for $\swclos{(A^S)}$ to be \sa, and so we need a
positive solution to the reductive algebra problem, at least for the
algebras which arise in this way from step (i).

This involvement of the reductive algebra problem would seem to
indicate that step (ii) is a significant step in the program outlined
above. However, in certain circumstances the entire problem is
contained in step (i). In particular, if step (i) can be completed for
every \cra, then \nref{cra conjecture} has a positive solution.

\begin{proposition}\label{RP -> reductive is desirable}
  Suppose that every \sw{}ly closed \cra\ is similar to a reductive algebra. Then
  every \sw{}ly closed \cra\ is similar to a \csalg.
\end{proposition}
\proof
Let $A\subseteq\BH$ be a \sw{}ly closed \cra. 
The amplification $\pma
A\infty$ has the \crp\ and so by assumption there is a similarity $S$
on $\pma\H \infty$ with $(\pma A\infty)^S$ reductive.
From \nref{Ideals split} $A$ is unital.
By \nref{Reflexive lemma} $\pma A\infty$ is reflexive, so 
$(\pma A\infty)^S$ is reflexive and reductive,
hence \sa. We have a completely
bounded faithful
\rep\ $\theta:(\pma A\infty)^S\ra A$ given by $(\pma a\infty)^S\mapsto
a$. This is similar to a \srep\ by 
\nref{cb reps are similar to sreps}. The image of this \srep\ is \sa\ and
similar to $A$.
\qed

Intriguingly, there are  situations where the only obstruction to
completing step (i) is the reductive algebra problem. The case of
abelian \oas\ discussed in the next section is an example of this.

\dealwithsectionbreaks
\section{Abelian algebras}
\label{Section abelian}

Abelian \oas\ are particularly interesting in terms of trying to prove
\nameref{tra conjecture}s~\ref{cra conjecture} and \ref{tra
  conjecture}. Unlike general \csalgs, the involution on an abelian
\csalg\ is uniquely determined by the (non-isometric) \ba\ structure
alone. Specifically, if $A$ is an abelian \csalg\ then $A$ is
isomorphic to $C_0(\p A)$, where $\p A$ is the maximal ideal space of
$A$. The only involution on $A$ compatible with a \cs~structure is the
natural involution on $C_0(\p A)$ given by complex conjugation of
functions.  This fact relieves the substantial technical burden of
making an arbitrary choice of involution when trying to erect a
\cs~structure on $A$ (\cf~the discussion on page~\pageref{Page with
  non-uniqueness}).

A theorem of Scheinberg 
provides the first step towards proving that abelian \tras\ are \isoc\
to \csalgs. 
This theorem states: every unital amenable uniform
subalgebra of $C(X)$ which separates $X$ is equal to $C(X)$
\cite{Helemskii-survey}.  This can be interpreted as a statement about
abelian \oas\ by recalling that every unital abelian \csalg\ is isometrically
\isoc\ to $C(X)$ for some $X$.

An examination of Scheinberg's proof shows that the
methods are `spatial' in character, and in fact apply also to 
\tra{}s.  While the statement of the next result appears rather
different to that of Scheinberg, the bulk of the proof is identical.
\begin{proposition}[{Proposition [Scheinberg]}]
  \label{Modified Scheinberg}
  \label{contained in abelian c* algebra is sufficient}
  Suppose $A\subseteq\BH$ is an abelian \tra\ contained in some
  abelian \cs-subalgebra of $\BH$. Then $A$ is \sa.
\end{proposition}
\proof By \nref{Unitisation lemma} the unital algebra $A_+$ generated
by $A$ and $1$ is a \tra. If $A_+$ is \sa\ then $A$ is also \sa, being
an ideal of $A_+$ \cite{DixmierC}. Thus we may assume that $A$ is unital.

Let $1\in A\subseteq B\subseteq\BH$ where $B$ is an abelian \csalg.
Since the \csalg\ $C^*(A)$ generated by $A$ is a subalgebra of $B$, we
may assume $B$ coincides with $C^*(A)$. This being the case, $B\cong C(X)$
for some compact Hausdorff $X$, and since $B$ is generated by $A$ as a
\sa\ algebra, $A$ separates the points of $X$.

We show that $A$ is all of $C(X)$. Let $\mu\in A^\perp\subseteq
M(X)$, $\lambda=|\mu|$ and $\xi=d\mu/d\lambda$. Then $|\xi|=1$
$\lambda$-a.e.\ and $\xi\in L^2(X,\lambda)$. The \hs\ $L^2(X,\lambda)$ is an
$A$-module, with the action given by pointwise multiplication of
functions. For $a\in A$ the equation $\<\mu,a>=\int a \xi\,d\lambda=0$ can
be written as $(a\cdot 1_X|\bar\xi)=0$, giving $\bar\xi\in\clos{A\cdot
  1_X}^\perp$. Since $A$ is a \tra, the submodule $\clos{A\cdot
  1_X}\subseteq L^2(X,\lambda)$
is the range of some module projection $p\in A'$. However, as all the
operators from $A$ are normal, Fuglede's theorem
\cite{RadjaviRosenthal} shows that $p$ commutes
with the \csalg\ generated by $A$. That is,
$p$ commutes with the action of $C(X)$ given by pointwise
multiplication. 

Since $\lambda$ is a regular measure, $C(X)$ is dense in
$L^2(X,\lambda)$ and there is a sequence $f_n\in C(X)$ with $f_n\ra
\bar\xi$ in $L^2(X,\lambda)$.  Since $A$ is unital we have
$1_X\in\clos{A\cdot1_X}$ and $p1_X=1_X$.  Then
$0=f_n\cdot(1-p)1_X=(1-p)f_n$, and so $f_n\in\clos{A\cdot1_X}$ for all
$n$. Thus $\bar\xi\in\clos{A\cdot1_X}\cap
\clos{A\cdot1_X}^\perp=\{0\}$. Since $|\xi|=1$ $\lambda$-a.e., it
follows that $\lambda=0$ and consequently $\mu=0$.  As $\mu$ was an arbitrary
element of $A^\perp$, the Hahn-Banach theorem implies that $A=C(X)$.
\qed

Since \reps\ of \sa\ \tras\ are similar to \sreps\ (\nref{Similarity
  property equivalent to trp}) we can generalise \nref{Modified
  Scheinberg} to give a characterisation dependent only on the \iso\ 
class of $A$.
\begin{proposition}
  Let $A\subseteq\BH$ be an abelian \tra. If $A$ is isomorphic to a
  closed subalgebra of an abelian \csalg, then $A$ is similar to a \csalg.
\end{proposition}
\proof If $\theta:A\ra C(X)$ is an isomorphism of $A$ onto a
subalgebra of an abelian \csalg, then \nref{Modified Scheinberg} shows
that $\theta(A)$ is an abelian \csalg. The inverse isomorphism
$\inv\theta:\theta(A)\ra A\subseteq\BH$ is a \rep\ of a \sa\ \tra, and
is thus similar to a \srep. The similarity takes $A$ to a \sa\ \oa.
\qed

These results show that in order to show that an abelian \tra\ 
$A\subseteq\BH$ is similar to a \csalg\ we need only embed $A$ in an
abelian \csalg.  In particular, a sufficient (and necessary) condition
is that $A''$ be similar to a \csalg. Since the spatial definition of
the \trp\  tends to give more direct information about $A'$ and
$A''$ than $A$ itself, this is a helpful simplification of the
problem.

It is possible to say something about the case of an abelian \cra\ as
well. It is not too surprising that the conclusion is in terms of the
\sw\ closure of $A$.
\begin{proposition}
  \label{cra Scheinberg}
  Suppose $A\subseteq\BH$ is a \cra\ contained in an abelian \csalg.
  Then the \sw\ closure of $A$ is \sa.
\end{proposition}
\proof
We prove first the case where $A$ acts nondegenerately.
Let $B\subseteq\BH$ be a \sa\ abelian \oa\ containing $A$.
Since the \sw\ closure of $B$ is also abelian and \sa, we may assume
that $B$ is \sw{}ly closed. Further, since the \sw\ closure of $A$ will
also have the \crp, we may assume that $A$ is \sw{}ly closed. As
$\infty$-fold amplification is a \sw\ homeomorphic isometric
isomorphism which preserves the involution on $\BH$, we
may replace $A$ with $\pma A\infty$, so that every \sw{}ly \cts\ functional
on $A$ is of the form $a\mapsto (a\xi|\eta)$ for some vectors $\xi,
\eta\in \H$. Finally, we may assume that
$C^*(A)$ is \sw{}ly dense in $B$.

Suppose that there is $a\in A$ with $a^*\not\in A$. Since $A$ is \sw{}ly 
closed this means that $a^*$ can be separated from $A$ by a \sw{}ly \cts\ 
functional, so there are vectors $\xi,\zeta\in\H$ with
$(A\xi|\zeta)=0$ and $(a^*\xi|\zeta)=1$. Writing $C=C^*(A)$, this
means that $\clos{A\xi}\neq\clos{C\xi}$. The $A$-module $\clos{C\xi}$
has the \rp, so there is an $A$-module projection $p$ from
$\clos{C\xi}$ onto $\clos{A\xi}$. Since $C$ is an abelian \csalg, its
restriction to $\clos{C\xi}$ is also an abelian \csalg, and consists
of normal operators on $\clos{C\xi}$.  Thus $A|_{\ssclos{C\xi}}$
consists of normal operators, and by Fuglede's theorem $p$ commutes
with $a^*|_{\clos{C\xi}}$. But by
\nref{Nondegeneracy lemma} and the assumption of nondegeneracy, we have
$\xi\in\clos{A\xi}=p(\clos{C\xi})$, so
$a^*\xi=a^*p\xi=pa^*\xi\in\clos{A\xi}$. As $\clos{A\xi}\perp \zeta$
and $(a^*\xi|\zeta)=1$, this is a contradiction. Hence $A$ is \sa.

To deal with the case when $A$ acts degenerately on $\H$, recall from
\nref{Degenerate case} that there is a projection $e\in A''$ with
$e\H=\clos{A\H}$. Since $A''\subseteq B''$ and $B''$ is an abelian
\csalg, $e$ must be \sa. Then $e\swclos{A}=\swclos{A|_{\clos{A\H}}}$,
and $\swclos{A}$ is \sa\ if $\swclos{A|_{\clos{A\H}}}$ is.  Since
$\clos{A\H}=e\H$ is a submodule for $B$, $B|_{\clos{A\H}}$ is an
abelian \csalg\ and the nondegenerate algebra
$A|_{\clos{A\H}}\subseteq B|_{\clos{A\H}}$ also satisfies the
hypotheses of the \nameref{cra Scheinberg}. The result follows.  \qed

An example illustrates that \sw\ closure is necessary in the statement
of this \nameref{cra Scheinberg}. Consider $\ell^\infty(\Zed)$ as an
\oa\ acting on $\ell^2(\Zed)$ by pointwise multiplication.  It is
readily verified that $\ell^\infty(\Zed)$ is a \sw{}ly closed \sa\ 
\oa, and so has the \crp.  Let $A\subseteq\ell^\infty(\Zed)$ be the
\oa\ generated by $c_0(\Zed)$ and the function $g:n\mapsto e^{in}$.
The ideal $c_0(\Zed)$ is \sw{}ly dense in $\ell^\infty(\Zed)$, so $A$
is also \sw{}ly dense in $\ell^\infty(\Zed)$. By \nref{cra and weak
  closure}, $A$ has the \crp.  However, $A$ is not \sa. To see this,
let $m$ be any invariant mean of $\Zed$, considered as a functional on
$\ell^\infty(\Zed)$ \cite{Paterson}. If $f\in c_0(\Zed)$ then
$\<m,f>=0$. We have $g^k(n)=e^{ikn}$, and denoting by $U$ the
right shift on $\ell^\infty(\Zed)$, also
$(U^lg^k)(n)=e^{ik(n-l)}=e^{-ikl}g^k(n)$. Then $\<m,
U^lg^k>=e^{-ikl}\<m,g^k>=\<m,g^k>$ by the translation invariance of
$m$. Since this holds for all $l, k\in \Zed$ it follows that
$\<m,g^k>=0$ for all $k>0$. The continuity of $m$ as a functional on
$\ell^\infty(\Zed)$ implies that $\<m,A>=0$.  On the other hand, the
\sa\ algebra generated by $A$ contains $\bar g$, and hence contains
$\bar gg=1_\Zed$. Since $\<m,1_\Zed>=1$, $A$ cannot be \sa.

\smallskip

If $A\subseteq\BH$ is a complete or total \RA\galley{Careful} then
\nref{Central idempotents are bounded} shows that the idempotents of
$A''$ are uniformly bounded, and a similarity may be applied to
transform these to \sa\ projections. Since the idempotents of abelian
\csalgs\ are always \sa, this is a necessary step towards
orthogonalising $A$.

Having applied such a similarity it follows from \nref{Commutant has
  RA} and the commutativity of $A$ that $A'$ is a reductive algebra. 
If it were known that every reductive algebra was
\sa, then this would imply that $A'$ and hence $A''$ are \sa, solving
both \nref{cra conjecture} and \nref{tra conjecture} for abelian \oas.

In the absence of a positive solution to the reductive algebra
problem, we can use direct integral theory to further reduce the
problem to a `transitive case'.
In \cite{AzoffFongGilfeather} a 
\label{Page with direct integral}%
disintegration theory is presented for non-\sa\ \oas. The idea is
that whenever $A$ is an \oa\ commuting with an abelian \vna\ $R$,  the
direct integral \rep\ of $R$ will yield some sort of direct integral
\rep\ of $A$. 
The result that we wish to use does not explicitly appear in
\cite{AzoffFongGilfeather}, and so we include a short proof. The
techniques used are the same as in \cite{AzoffFongGilfeather}, which
in turn are borrowed from the \sa\ direct integral theory
\cite{Schwartz}. We use the notation of \cite{Schwartz} and assume
that the reader is familiar with the basic ideas of direct integral theory.

A technical point arises here, which is that the direct integral
theory works properly only when the algebras and underlying \hs{}s
involved are separable. Therefore the following result will apply
directly only to separable algebras acting on separable \hs{}s.

\begin{proposition}
\label{Direct integral result}
  Let $\H$ be a separable \hs\ and $A\subseteq\BH$ a separable abelian
  \tra, such that the projections of $A''$ are \sa. Then there
  is a finite measure space $(\Lambda,\mu)$, a measurable family
  $\lambda\mapsto \H_\lambda$ of \hs{}s, and a measurable family
  $\lambda\mapsto \pi_\lambda$ of \reps\ of $A$ such that: 
  $\H$ is isometrically isomorphic to the direct integral space $\H=\dint
  \H_\lambda.d\lambda$, the
  projections of $A''$ are the diagonal projections with respect to
  this disintegration of $\H$, $a=\dint
  \pi_\lambda(a).d\lambda$ for all $a\in A$, and $\pi_\lambda(A)'$ is
  transitive for almost all $\lambda$.
\end{proposition}
\proof Let $R$ denote the abelian \vna\ generated by the central projections
of $A''$.  The usual direct integral theory for abelian \vnas\ shows
that there is a regular measure space $(\Lambda, \mu)$, and a
measurable family $\lambda\mapsto\H_\lambda$ of \hs{}s \st\ $\H=\dint
\H_\lambda .d\mu(\lambda)$, so that $R$ becomes the algebra of
diagonal operators with respect to this identification (see
\cite{Schwartz} for details).

Since $A$ is separable a standard construction allows us to choose a
family $\{\pi_\lambda:A\ra \B(\H_\lambda)\}$ of contractive \reps\ so
that $\id_A=\dint \pi_\lambda.d\lambda$.  We claim that $\pi_\lambda(A)'$
is transitive for almost all $\lambda$.

We may partition $\Lambda$ into sets where $\H_\lambda$ is of
constant dimension, and assume without loss of generality that on
these sets the map $\lambda\mapsto \H_\lambda$ is constant. Since
there are at most countably many of these sets it suffices to prove
the result in the case where there is only one such set. Then we may
replace $\H_\lambda$ with a constant \hs\ $\H'$.

Observe that since $\H$ is separable, there is a
strongly dense sequence $\{T_i\}$ in $A'$. 
In addition to the strong topology 
we will need to use the $*$-strong topology on $\BH$, which is defined
by the family of seminorms
$$
x\mapsto \norm{x\xi}+\norm{x^*\xi}, \hbox{ for $\xi\in\H$}.
$$
Since $T_i\in A'$, we have
$$
T_i=\dint T_i(\lambda).d\lambda
$$
where the function $\lambda\mapsto T_i(\lambda)\in\B(\H')$ is weakly
measurable. The Borel sets of the $*$-strong topology are also weak
Borel sets, and so the maps $\lambda\mapsto T_i(\lambda)$ are
$*$-strongly Borel. By Lusin's theorem, for every $m\in\N$ there is a
closed set $F_m\subseteq \Lambda$ with $\mu(\Lambda\setminus F_m)\leq
{1\over m}$ and $\lambda\mapsto T_i(\lambda)$ $*$-strongly \cts\ on
$F_m$ for
each $i$.

We put 
$$G_{n,m}=\{(\lambda,p): p \in\B(\H'), p^2=p, \lambda\in
F_m \mbox{ and } T_n(\lambda)p-p T_n(\lambda)=0\}.
$$
It follows that 
$G_{n,m}$ is a closed subset of $\Lambda\times \B(\H')$ 
for the given topology on $\Lambda$ and the $*$-strong topology on
$\B(\H')$. 
This implies that the set $G=\bigcup_m\bigcap_n G_{n,m}$ is a Borel
set. Define the Borel set $M=G\setminus \bigg(\{(\lambda, 0):
\lambda\in \Lambda\} \cup \{(\lambda, 1): \lambda\in \Lambda\}\bigg)$,
and let $N\subseteq\Lambda$ be the image of the projection of $M$ onto
$\Lambda$.

If $(\lambda,p)\in M$, then by the definition $p$ is a proper
idempotent in $\B(\H')$ commuting with $\{T_i(\lambda)\}$ for
all $i$. Conversely, if such an idempotent exists, then $\lambda\in
N$. 
By lemma 4.6 of \cite{AzoffFongGilfeather} the operators
$\{T_i(\lambda)\}$ generate $\pi_\lambda(A)'$ as a strongly closed
algebra for almost all $\lambda$.
Thus the idempotents $p$ appearing in $M$ 
are exactly the proper idempotents of
$\pi_\lambda(A)''$. Note that since all the $\pi_\lambda$ are
contractive, there is $K>0$ such that if
$p_\lambda$ is an idempotent of
$\pi_\lambda(A)''$ then $\norm{p_\lambda}\leq K$.

$N$ is an analytic set, so by the principle of measurable selection
there is a Borel set $E\subseteq N$ with $\mu(N\setminus E)=0$, and a
measurable choice function $\phi:E\ra M$ \cite{Schwartz}.  
This yields a measurable operator valued function $p:E\ra \B(\H')$
given by composing $\phi$ with projection onto the $\B(\H')$ component.
We may extend
the function $p$ to all of $\Lambda$ by setting $p(\lambda)=0$ for
$\lambda\not\in E$. Then $\dint p(\lambda).d\lambda$ is an idempotent
commuting with $\{T_i\}$ for all $i$. Hence $\dint
p(\lambda).d\lambda\in A''$, and $p(\lambda)\in\{0,1\}$ almost
everywhere. Since $M$ does not contain the operators $0$ or $1$ over
any point in $\Lambda$ we conclude that $\mu(E)=0$, and so $\mu(N)=0$.
This means that the set of $\lambda$ for which $\{T_i(\lambda)\}'$
contains a proper idempotent is null.  Thus the set of
$\pi_\lambda(A)''$ which contain proper idempotents is null.  Since
$\clos{\pi_\lambda(A)}$ is an abelian \tra\ this implies that
$\pi_\lambda(A)'$ is transitive for almost all $\lambda$.  \qed

This result allows us to reduce the problem of understanding abelian
\tra{}s to a special case of the transitive algebra problem.
Suppose that we were able to show that the transitive algebras
$\pi_\lambda(A)'$ are equal to $\B(\H')$. Then
$\pi_\lambda(A)\subseteq\pi_\lambda(A)''=\C1_{\H'}$, so we would get
an embedding of $A$ into the abelian
\csalg\ $\dint \C.d\lambda$. Then \nref{Modified Scheinberg} would show
that $A$ is \sa.

There are partial solutions to the transitive algebra problem which
allow us to address \nref{tra conjecture} in certain special cases. In
particular, we have the following result of Lomonosov, which solves
the transitive algebra problem for algebras containing compact operators.
\begin{lemma}[Lomonosov's Lemma]\label{Lomonosov}\index{Lomonosov's lemma}
  Let $A\subseteq\BH$ be a (not necessarily closed) algebra with $\lat
  A=\{0,\H\}$.  If $A$ contains any non-zero compact operator then $A$
  contains a compact operator with non-zero spectral radius.
  Consequently, any weakly closed transitive \oa\ containing a
  non-zero compact operator is all of $\BH$.
\end{lemma}
\proof See e.g., \cite[Lemma 8.22]{RadjaviRosenthal} for a proof.
\qed This `lemma' permits a complete analysis of abelian \tras\ which
commute with sufficiently many compact operators. The next result
shows that in this case the rather technical direct integral
methods of \nref{Direct integral result} can be avoided
completely.
\begin{lemma}
  \label{Discrete result}
  Suppose that $B\subseteq\K(\H)$ is an algebra of compact operators
  acting nondegenerately on $\H$. If $R$ is an abelian \vna\ commuting
  with $B$, then $R$ is generated as a \vna\ by its
  minimal projections.
\end{lemma}
\proof
Since $R$ is \sa, it commutes not only with $B$ but also with
$\setof{b^*:b\in B}$ and hence with the \csalg\ $C^*(B)$ generated by
$B$.  This is a \csalg\ of compact operators acting nondegenerately on
$\H$, so by the structure theory discussed on page~\pageref{Page with
  compact csalgs} there is a family $\H_\gamma$ of \hs{}s and
multiplicity function $\gamma\mapsto n_\gamma$, such that $C^*(B)$ is
unitarily \isoc\ to $\sum^{c_0}_\gamma\K(\H_\gamma)\ot 1_{n_\gamma}$.
The commutant $M$ of $C^*(B)$ may be read off from this as
$\sum_\gamma^{\ell^\infty} 1_{\H_\gamma}\ot
M_{n_\gamma}\cong\sum_\gamma^{\ell^\infty} M_{n_\gamma}$.  We write
$M_\gamma$ for the $\gamma^{th}$ summand. The projections from $M$
onto $M_\gamma$ are central projections of $M$.  Let $D$ be a maximal
abelian \sa\ subalgebra of $M$ containing $R$.  The maximality of $D$
implies that $D$ contains the centre of $M$, and so
$D=\sum_\gamma^{\ell^\infty} (D\cap M_{\gamma})$, where $D\cap
M_{\gamma}$ is a maximal abelian \sa\ subalgebra of $M_{\gamma}$.  The
maximal abelian \sa\ subalgebras of $M_\gamma$ are \isoc\ to
$\ell^\infty(n_\gamma)$, and so $D\cong \sum^{\ell^\infty}
\ell^\infty(n_\gamma)\cong \ell^\infty(\Omega)$ for some index set
$\Omega$.  This allows us to think of $R$ as a \sa\ subalgebra of
$\ell^\infty(\Omega)$.  Let $\Lambda$ be the set of equivalence
classes of $\Omega$ under the relation
$$
\omega_1 \equiv \omega_2 \Longleftrightarrow
x(\omega_1)=x(\omega_2) \hbox{ for all $x\in R.$}
$$
By this construction we may think of $R$ as a separating \sa\ 
subalgebra of $\ell^\infty(\Lambda)$.  Consider the set
$$P_\lambda=\{p\in R: 
\hbox{ $p$ is idempotent and $p(\lambda)=1$}\}.
$$ Since the projections in $R$ form a complete Boolean algebra,
$P_\lambda$ will possess an infimum, $p_\lambda$ say, which must be
given by a characteristic function of a subset of $\Lambda$ containing
$\lambda$. Because $R$ is generated as a norm closed algebra by its
projections and separates the points of $\Lambda$, we have
$p_\lambda=1_{\{\lambda\}}$, and so $R=\ell^\infty(\Lambda)$.  The
lattice of projections of $R$ is thus ${\cal P}(\Lambda)$, which is
atomic with atoms $\{1_{\{\lambda\}}\}$.  These are exactly the
minimal projections of $R$.  
\qed 
This means that if $A'$ contains a
nondegenerate subalgebra of compact operators, then the direct
integral of \nref{Direct integral result} is in fact a discrete direct
sum. Note that there are no separability requirements for $A$ or $\H$
in this case.

\begin{theorem}
\label{Abelian and compact}
  Let $A\subseteq\BH$ be an abelian \cra, and suppose that $A'\cap
  \Kom(\H)$ acts nondegenerately on $\H$. Then there exists a
  similarity such that $A$ is contained in an 
  abelian \csalg. Hence the \sw\ closure $\swclos{A}$ is similar to a
  \csalg. If $A$ also has the \trp, then $A$ is itself similar to a \csalg.
\end{theorem}
\proof Using \nref{Central projections are sa} we may assume that the
central projections of $A''$ are \sa.  Letting $R$ denote the \vna\ 
generated by these projections, \nref{Discrete result} implies that
$R$ is generated by its minimal projections.  This gives rise to an
orthogonal decomposition $\H=\sum^\oplus\H_\lambda$, where the
$\H_\lambda$ are the ranges of the minimal projections in $R$. Since
these projections are central, it follows that
$(A|_{\H_\lambda})''=A''|_{\H_\lambda}$, and by minimality there are
no proper idempotents in $(A|_{\H_\lambda})''$. Since
$\clos{A|_{\H_\lambda}}$ has the \crp, \nref{Commutant has RA} implies
that $(A|_{\H_\lambda})'$ is a transitive algebra. Since
$\clos{A'\cap\Kom(\H)\cdot\H}=\H$, we have
$\clos{(A|_{\H_\lambda})'\cap
  \Kom(\H_\lambda)\cdot\H_\lambda}=\H_\lambda$, and so
$(A|_{\H_\lambda})'$ contains non-zero compact operators. By
Lomonosov's lemma this implies that
$A|_{\H_\lambda}=\C1_{\H_\lambda}$, and so $A$ is contained in the
abelian \csalg\ $R$. The remaining assertions follow from \nref{cra
  Scheinberg} and \nref{Modified Scheinberg}.  \qed

\Nref{Abelian and compact} is related to a result of Rosenoer
\cite{Rosenoer1}, which shows that if $A$ is a \RA\ with
$A'\cap\Kom(\H)$ acting nondegenerately on $\H$, then under various
additional hypotheses (which amount to requiring that the algebra be
multiplicity free) $\lat A$ `admits spectral synthesis', meaning
that every $\V\in\lat A$ is the span of its one-dimensional
submodules. The techniques used in \cite{Rosenoer1} derive from the
theory of almost-periodic functions on groups rather than the purely
\oa{}ic approach we employ. \Nref{Abelian and compact} also contains
the result of Willis discussed on page~\pageref{Page with Willis}.

\subsection{A spectral approach to abelian \tras}

Another approach to proving \nref{tra conjecture} for abelian \tras\ 
is to capitalise on the observation of the beginning of this section,
that an abelian \ba\ $A$ is \isoc\ to an abelian \csalg\ \iff\ the
Gelfand transform is an \iso\ of $A$ onto $C_0(\p A)$.

It follows from \nref{Modified Scheinberg} that if $A$ is an abelian
\tra\ and $\Gamma$ is the Gelfand transform, then $\Gamma(A)$ is
uniformly dense in $C_0(\p A)$.  If $A$ has `enough' ideals, it is
possible to use the estimate of \nref{Bounded decomposition}
to show that in fact the Gelfand transform maps onto
$C_0(\p A)$. We can ensure the existence of a sufficiently rich
collection of ideals by assuming that $A$ is regular (\ie, the \ws\ 
topology on $\p A$ coincides with the hull-kernel topology).  

We will need a few preparatory results.
\begin{lemma}\label{Edwards lemma}
  Suppose $T:X\ra Y$ is a map between \bs{}s $X$ and $Y$, such that $\bd
  T(\bd X)$ is norm closed in $\bd Y$. Then $T(X)$ is norm closed in $Y$.
\end{lemma}
\proof
Let $f\in\B(X,Y)$ be a bounded linear map. Theorem 8.7.3 of
\cite{Edwards} states that if $f^*(Y^*)$ is norm closed in $X^*$ and
$f(X)$ is norm dense in $Y$ then $f(X)=Y$. If $f^*(Y^*)$ is norm
closed but $f(X)$ is not norm
dense, we may consider the corestricted map $\hat f:X\ra
\clos{f(X)}$. The dual map ${\hat f}^*:Y^*/\clos{f(X)}^*\ra X^*$ has
the same range as $f^*$, and so the closedness of $f^*(Y^*)$ implies
$f(X)=\clos{f(X)}$.

Applying this observation firstly to $T^*$ and then $T$ gives the
result.
\qed

The next result continues the theme established by \nameref{l-infty
  norms}s~\ref{Bounded decomposition} and \ref{l-infty norms}.
Note that we do not require $A$ to be abelian here.
\begin{proposition}
  \label{Bound via sufficient quotients}
  Let $A$ be a \tra\ and $\{K_i\}$ a finite family of two-sided ideals
  with $\bigcap_i K_i=\{0\}$. Let $M$ be the ideal constant of $A$.
  Then for all $a\in A$ we have $\norm{a}\leq
  (1+2M)^2M\sup_i\{\norm{a+K_i}\}$.
\end{proposition}
\proof
Consider the natural homomorphism $\theta:A\ra\sum^{\ell^\infty}
A/K_i$. The condition $\bigcap K_i=\{0\}$ implies that $\theta$ is
one-one. It is required to show that $\theta$ is an isomorphism onto a
closed subalgebra with $\norm{\inv\theta}\leq (1+2M)^2M$. 

Applying the bidual functor yields the homomorphism
$$
\bd \theta: \bd A\ra \sum{}\!^{\ell^\infty} \bd A/ \bd {K_i}.
$$
The ideals $\bd {K_i}\ideal \bd A$ correspond to central projections
$k_i$, and under the isomorphisms $\bd A/\bd {K_i}\cong (1-k_i)\bd A$
the map $\bd \theta$ becomes
$$
\theta': \bd A \ra \sum{}\!^{\ell^\infty} (1-k_i)\bd A,
$$
given by
$$
\theta'(a)=\sum{}\!^\oplus (1-k_i)a.
$$ We have the bounds $\norm{\bd\theta(a)}\leq \norm{\theta'(a)}\leq
M\norm{\bd\theta(a)}$, since the isomorphism from $\bd A/\bd K_i$ to
$(1-k_i)\bd A$ has norm at most $M$ and contractive inverse.  Embedding
$\bd A$ in $\BH$ via the GNS \rep, there is a similarity $S$ on $\H$
with $\norm{S}\norm{\inv S}\leq 1+2M$ such that the central
idempotents of $(\bd A)^S$ are \sa.  It is clear that the range of
$\theta'$ is norm closed, and has kernel $(\prod k_i)\bd A$.  Hence
$\bd \theta$ has norm closed range, and by \nref{Edwards lemma} it
follows that $\theta$ has norm closed range.  This means that $\theta$
is an isomorphism onto a closed subalgebra of
$\sum^{\ell^\infty}A/K_i$.

It remains to find a bound for $\norm{\inv\theta}$. Since $\theta$ is
an isomorphism so is $\bd \theta$ and also $\theta'$. This means
that $\prod k_i=0$ and $\bigvee (1-k_i)=1$. Let $\{j_i\}$ be a
finite set of central projections with $j_i\leq 1-k_i$, $j_ij_{i'}=0$
for $i\neq i'$ and
$\sum j_i=1$. If $a\in \bd A$ we have $a=\sum j_ia$,
and since $j_i\leq 1-k_i$ we have $\norm{(j_ia)^S}\leq
\norm{[(1-k_i)a]^S}$. This gives the estimate
\begin{eqnarray*}
\norm{a} & \leq &(1+2M)\norm{a^S} \\
         & = & (1+2M)\bnorm{\sum{(j_ia)}^S} \\
         & = & (1+2M)\supof{\norm{(j_1a)^S}} \\
         & \leq & (1+2M)^2\supof{\norm{(1-k_i)a}} \\
         & = & (1+2M)^2\norm{\theta'(a)} \\
         & \leq & (1+2M)^2M\norm{\bd \theta(a)}.
\end{eqnarray*}
Since $\inv{\bd\theta}$ extends $\inv\theta$ this gives the desired
bound.
\qed

The
following proof is inspired by the result of \cite{ElliotMultiplier},
where a similar conclusion is established in a rather different context.

\begin{theorem}
  Let $A$ be a unital, semisimple and regular abelian \tra. Then
  $A\cong C(\p A)$.
\end{theorem}
\proof
By the open mapping theorem it is sufficient to show that the range of
the Gelfand transform is all of $C(\p A)$. 
Let $f\in C(\p A)$ be a \cts\ function. We construct a Cauchy sequence
$\net{a_n}$ of elements from $A$ such that $\Gel{a}_n\ra f$. The
completeness of $A$ and the continuity of the Gelfand map then gives
the desired result.

It is sufficient to deal with the case where the range of $f$ is
contained in $[0,1]$. For each $n\in\N$ and $-1\leq i\leq 2^n$ let
$U_{n,i}=(2^{-n}i, 2^{-n}(i+2))$. The set $\{U_{n,i}\}_{-1\leq i\leq
  2^n}$ is a finite open cover for $[0,1]$ by intervals of length
$2^{1-n}$. Since $A$ is regular, the sets $\inv f (U_{n, i})$ are open
in the hull-kernel topology. We define the ideals $J_i^n=\ker
(\p A \setminus \inv f (U_{n, i}))$. 
If $\sum_i J_i^n\neq A$ then there is $\omega\in\p A$ with
$\ker\omega\supseteq\sum_i J_i^n$. Then
$\omega\in\bigcap\mathop{\hbox{hul ker}}(\p
A\setminus\inv f(U_{n,i}))=\emptyset$; this contradiction implies that
$\sum_i J_i^n=A$,
which means that we can use \nref{Bounded decomposition} to find
$x_i^n\in J^n_i$ with $\norm{x^n_i}\leq 2M$ and $1=\sum_i x_i^n$. We
set $a_n=\sum_i 2^{-n}i x_i^n$. It is straightforward to verify that
$\Gel a_n\ra f$ in the uniform norm on $C(\p A)$.

It remains to show that $\{a_n\}$ is a Cauchy sequence in $A$. 
For every integer $m\geq 1$ we define ideals $K_i^m$ by 
$$
K_i^m=J_{-1}^m + J_0^m + \cdots + J_{i-2}^m + J_{i+2}^m + \cdots +
J_{2^m}^m.
$$
From this definition we have $K^m_i\subseteq\ker {\inv f}(U_{m,i})$.
Thus 
$$\bigcap_i\K^m_i\subseteq\ker\bigcup_i {\inv f}(U_{i,m})=\{0\}.$$ 
Concentrating on fixed $m$ and $i$ we have the following equation in
the quotient $A/K_i^m$:
$$
a_m + K_i^m= 2^{-m}(i-1)x_{i-1}^m + 2^{-m}ix_{i}^m +
2^{-m}(i+1)x_{i+1}^m + K_i^m.
$$
For an integer $n< m$ it is possible to get a similar equation for
$a_n$ in the quotient $A/K^m_i$.  Note that at most three of the open
intervals $\{U_{n,j}: -1\leq j \leq 2^n\}$ overlap with $U_{m, i-1}
\cup U_{m, i} \cup U_{m, i+1}$.  For any integer $j$ with
$$
U_{n, j}\cap (U_{m, i-1} \cup U_{m, i} \cup
U_{m, i+1})=\emptyset
$$
we have $J_j^n\subseteq K_i^m$, since if $a\in
J_j^n$ then 
$$
a=a\cdot1=\sum_k ax_k^m=\sum_{k\not\in\{i-1,i,i+1\}} ax_k^m\in K_i^m.
$$
This means that $J_j^n\subseteq K_i^m$ except for at most
three consecutive values of $j$. Let $j$ be an integer
such that
$J_k^n\not\subseteq K_i^m$ implies $k\in\{j-1, j, j+1\}$. Then 
$$
a_n + K_i^m = 2^{-n}(j-1)x_{j-1}^n + 2^{-n}jx_j^n +
    2^{-n}(j+1)x_{j+1}^n + K^m_i,
$$
and $|2^{-n}j-2^{-m}i|\leq 2^{-n+1}$.
We also have 
$x_{i-1}^m+x_i^m+x_{i+1}^m+K_i^m = 1 + K_i^m$ and
$x_{j-1}^n+x_j^n+x_{j+1}^n+K_i^m = 1 + K_i^m$.
Thus
\begin{eqnarray*} 
a_n-a_m + K_i^m & = &
  a_n-m_m + 
        2^{-m}i\cdot 1 - 2^{-m}i\cdot 1 + K^m_i \\
     & = &
  (2^{-n}(j-1)x_{j-1}^n- 2^{-m}i x_{j-1}^n) \\
  & - & (2^{-m}(i-1)x_{i-1}^m- 2^{-m}i x_{i-1}^m) \\
  & + & 
  (2^{-n}jx_{j}^n- 2^{-m}i x_{j}^n) \\
  & - & (2^{-m}ix_{i}^m- 2^{-m}i x_{i}^m) \\
  & + & 
  (2^{-n}(j+1)x_{j+1}^n- 2^{-m}i x_{j+1}^n) \\
  & - & (2^{-m}(i+1)x_{i+1}^m- 2^{-m}i x_{i+1}^mx) + K_i^m.
\end{eqnarray*}
This gives the bound $\norm{a_n-a_m+K_i^m}\leq 12\cdot
M\cdot2^{-n+2}$. \Nref{Bound via sufficient quotients} 
then implies that $\norm{a_n-a_m}\leq 48(1+2M)^2M^2\cdot 2^{-n}$, so
$\{a_n\}$ is a Cauchy sequence in $A$.
\qed

It is unclear whether the hypothesis of regularity is needed here.  It
is only used to ensure that $\sum_i J_i^n=A$. It is interesting to
note that for regular uniform algebras this result gives an
alternative proof of Scheinberg's result, using totally different
methods.

It is surprising that this method does not refer to the either the
reductive or transitive algebra problem at all. The catch is that the
method also ignores the radical of $A$. Recall from \nref{CP is
  hereditary} that the radical of a \tra\ is also a \tra, and so it is
ultimately important to be able to deal with the
possibility of a radical \tra.\galley{Too defeatist?}

\comment{
It is shown in \cite{Ulger} that if $A$ is an amenable abelian \ba\
then the spectrum $\p A$ is discrete when considered as a subspace of
$A^*$ with the $\sigma(A^*,\bd A)$-topology (\ie\ the weak topology on
$A^*$). This is a generalisation
of a theorem from \cite{Gourdeau} which shows that $\p A$ is discrete
with the norm topology. 
We show here that if an abelian \oa\ $A$ has the \trp\ then $\p A$ is
discrete in the weak topology.

\begin{lemma}
  Let $A$ be an abelian \tra. Then $\p A$ is discrete in the weak topology.
\end{lemma}
\proof
Let $\omega\in\p A$. The kernel of $\omega$ is a codimension 1 ideal
$J_\omega\ideal A$. Taking the bidual of this inclusion gives $\bd
{J_\omega}\ideal \bd A$, where $\bd {J_\omega}$ is a codimension 1
\ws\ closed ideal. By \nref{Ideals split} there is a projection
$q_\omega\in \bd A$ with $\bd A=\bd {J_\omega}\oplus q_\omega \bd
A$. The kernel of $\omega$ is $\bd {J_\omega}$.

If $\omega'\in\p A$ is any other character with corresponding
projection $q_{\omega'}$, then since $q_\omega$ is minimal we have
$q_\omega q_{\omega'}=0$. Thus $q_{\omega'}\subseteq \bd
{J_\omega}=\ker \bd \omega$, and so $\bd
\omega(q_{\omega'})=\delta_{\omega,\omega'}$. This means that
$\omega'$ is clopen in the $\sigma(\p A, \bd A)$-topology.
\qed
}

\vfill\break

\dealwithsectionbreaks
\section{Algebras of compact operators}
\label{section with compact case}

We have seen that in certain cases the obstacle to proving \nref{tra
  conjecture} is the unknown status of the transitive algebra problem.
However, Lomonosov's lemma effectively solves this problem for
algebras containing non-zero compact operators. It is possible to use
the lever provided by Lomonosov's lemma to completely analyse \tras\ 
consisting of compact operators. In fact, we will prove that if
$A\subseteq\KH$ is a \cra, then $A$ is similar to a \csalg.

\smallskip
\noindent\hbox to \hsize{\hskip\parindent\hrulefill\hskip\parindent}

{\scshape Remark:} After the work for this section had been completed,
the paper \cite{RosenoerAnnoying} by Rosenoer came to my attention. In
that paper the author essentially studies \RA{}s which contain a
nondegenerate subalgebra of compact operators.  
Although Rosenoer studies the \rp\ in the context of algebras of
operators on a general \bs, his techniques are very similar to ours.
His use of the \rp\ rather than the \crp\ means that
he admits algebras like that of section~\ref{rp does not imply
  crp}. Consequently he does not find a result analogous to
\nref{Bounded similarity for compact}, and does not find a
characterisation of the subalgebras of $\KH$ which are similar to
\csalgs.

\noindent\hbox to \hsize{\hskip\parindent\hrulefill\hskip\parindent}
\smallskip

The description of \csalgs\ of compact operators given on
page~\pageref{Page with compact csalgs} will guide our approach.
Recall that if $B\subseteq\KH$ is a nondegenerate \csalg, then there
is a family $\{V_\lambda\}$ of irreducible submodules of $\H$ and
positive integers $\{n_\lambda\}$ so that $\H$ is isometrically \isoc\ 
with $\sum^{\oplus_2}\V_\lambda\ot_2\C^{n_\lambda}$ and $B$ is
unitarily equivalent to $\sum^{c_0} \K(\V_\lambda)\ot 1_{n_\lambda}$.
 
From this description it follows that if $A\subseteq\KH$ is similar to
a \csalg\ and $A''$ contains no proper central projections, then
$A\cong \K(\V)\ot 1_n$ for some irreducible submodule $\V\in\lat A$
and some positive integer $n$.  It is convenient to take this as a
starting point, by showing that if $A\subseteq\KH$ is a \cra\ and
$A''$ contains no proper central projections, then $\lat A$ contains
an irreducible submodule $\V$ and $A\cong \K(\V)\ot 1_n$ for some
$n\geq 1$.

The tools used to establish the existence of an irreducible submodule
are Lomonosov's lemma and an argument concerning quasinilpotent
compact operators due to Shul'man. To get to Shul'man's result it is
necessary to introduce the notion of a \cts\ nest of invariant
subspaces. 

If $N\subseteq\lat A$ is a totally ordered subset of $\lat
A$ and $\V\in N$, we write $\V_-$ for the subspace
$$
\V_-=\cspam{\W\in N: \W\subseteq\V}.
$$
When $\V_-=\V$ for all $\V\in N$ we say that $N$ is a \cts\ nest in
$\lat A$.  

\begin{lemma}\label{cts nest}
  Let $A\subseteq\KH$ have the \crp, and suppose that $\lat A$ contains no
  non-zero
  irreducible submodules. Then there is a \cts\ nest $N\subseteq\lat
  A$.
  Consequently, every element of $A$ is quasinilpotent.
\end{lemma}
\proof A simple
application of Zorn's lemma shows that maximal totally ordered subsets
of $\lat A$ always exist.
Let $N\subseteq\lat A$ be such a subset and take $\V\in
N$.  Since $A$ has the \crp, $\V$ has the \rp\ and there is $\w\in
\lat A$ with $\V=\V_-\oplus\w$. The maximality of $N$ implies that $\w$ is
irreducible.  Since $\lat A$ has no non-zero irreducible submodules, $\w=\{0\}$ and
$\V_-=\V$.

A well-known theorem of Ringrose \cite[Theorem 3.4]{Davidson} shows
that since every operator in $A$ is compact, the existence of a \cts\ 
nest in $\lat A$ implies the quasinilpotence of every element of $A$.
\qed

\begin{lemma}
  \label{qn lemma}
  Let $A\subseteq\KH$ be an \oa\ and suppose that $\lat A$ contains a
  \cts\ nest $N$. Then for any finite set $\{a_i\}_{i\leq n}\subseteq
  A$ and any $\epsilon>0$, there is a finite chain $0=\v_0\subset
  \v_1\subset \ldots \subset \v_m=\H$ of invariant subspaces from $N$
  with $\norm{p_{\v_{j+1}\ominus \v_j}a_i|_{\v_{j+1}\ominus
      \v_j}}<\epsilon$ for all $i\leq n$ and all $j<m$, where
  $p_{\v_{j+1}\ominus \v_j}$ denotes the orthogonal projection of $\H$
  onto $\v_{j+1}\ominus \v_j$.
\end{lemma}
\proof Let $N\subseteq\lat A$ be a \cts\ nest.  Then Theorem 3.5 of
\cite{Davidson} shows that if $n=1$ the result is true. For $n>1$ the
union of the chains obtained for each $a_i$ satisfies the requirements
of the lemma.  \qed

Shul'man's theorem now provides the critical information for dealing
with \cts\ nests.  We reproduce the proof here since the argument
is both elementary and surprising.
\begin{theorem}[{Theorem [Shul'man]}]
  \label{Shulman}
  Let $A\subseteq\KH$ be an \oa, and suppose that $\lat A$
  contains a \cts\ nest. Let $\{a_i\}_{i\leq
    n}\subseteq A$ and $\{b_i\}_{i\leq n}\subseteq A'$ be two finite
  sets. Then the operator $\sum a_ib_i$ is quasinilpotent.
\end{theorem}
\proof We need to show that $\lim_{k\ra \infty} \norm{(\sum
  a_ib_i)^k}^{1 \over k}\ra 0$. We may assume that the operators $a_i$
and $b_i$ are all contractions without loss of generality.  Write
$E=\{a_i\}_{i\leq n}$ and $F=\{b_i\}_{i\leq n}$.  The $k^{th}$ power
of $\sum a_ib_i$ expands to a sum of the form $\sum_{l\leq n^k} c_l
d_l$, where each $c_l$ is a product of $k$ operators from $E$ and
$d_l$ is a corresponding product of $k$ operators from $F$. Since
$\norm{d_l}^{1\over k}\leq 1$ for all $k$ and $l$, it will be
sufficient to show that for all $\epsilon>0$ we have
$$
\lim_{k\ra \infty} \left[\supof{\norm{a}^{1\over k}: \mbox{
    $a$ is a product of $k$ elements from $\{a_i\}$}}\right] < \epsilon.
$$
For, if this limit is established, then 
$$
\bnorm{\left(\sum a_ib_i\right)^k}^{1\over k}=
\bigg\|\sum_{l\leq n^k} c_l d_l\bigg\|^{1\over k}\leq 
\bigg(\sum_{l\leq n^k}\norm{c_l}\bigg)^{1\over k} \leq
n\epsilon
$$ for large $k$.

To obtain the estimate on norms of $k^{th}$ degree products from $a_i$,
take $\epsilon>0$ and 
let $0=\v_0\subset \v_1\subset \cdots \subset \v_m=\H$ be a chain as in 
\nref{qn lemma}. If $a_i\in E$ and $\xi\in \v_j$ has $\norm{\xi}\leq
1$, then $a_i \xi=\xi_0 \oplus \xi_1$, where $\xi_0\in \v_j\ominus
\v_{j-1}$, $\xi_1\in \v_{j-1}$, $\norm{\xi_0}\leq\epsilon$ and
$\norm{\xi_1}\leq 1+\epsilon$. Inductively applying this observation
shows that if $\xi\in \H$ has $\norm{\xi}\leq 1$ and $a$ is a product
of $k$ elements from $E$, then we have the decomposition
$$
a\xi=\xi_{0,0,\cdots,0}+\xi_{0,0,\cdots,1}+\cdots + \xi_{1,1,\cdots,1}
$$
where there are $k$ subscripts and $2^k$ terms in the sum. Each
subscript takes either the value $0$ or $1$, indicating whether or
not to `go down' one subspace. Since $\v_0=0$ and $\v_m=\H$, the terms
with $m$ or more $1$'s vanish. Also, the above remarks show that a
term with $i$ $1$'s and $k-i$ $0$'s has norm less than
$(1+\epsilon)^i\epsilon^{k-i}$. Thus we obtain the estimate
$$
\norm{a\xi}\leq \sum_{i<m} \choose k i (1+\epsilon)^i \epsilon^{(k-i)}.
$$
Assuming $\epsilon<1$ we can obtain a crude estimate for this which is
sufficient for our purposes. Note that 
\begin{eqnarray*}
  \norm{a\xi} & \leq  &
 \sum_{i<m} \choose k i (1+\epsilon)^i \epsilon^{(k-i)} \\
 & \leq &
 \left(\sum_{i<m} \choose k i\right) (1+\epsilon)^m \epsilon^{(k-m)}
\end{eqnarray*}
and if $k>2m$ then
$$
\sum_{i<m}\choose k i \leq m\choose k m = {k!\over (m-1)!(k-m)!}.
$$
We are interested in the value of $\norm{a\xi}^{1\over k}$. The above
shows that 
$$ \norm{a\xi}^{1\over k}\leq \left({k!\over
    (k-m)!(m-1)!}\right)^{1\over k}(1+\epsilon)^{m\over
  k}\epsilon^{1-{m\over k}} \leq k^{m\over k}(1+\epsilon)^{m\over
  k}\epsilon^{1-{m\over k}}\ra \epsilon.
$$
Since this estimate holds for any $\xi\in \H$ with $\norm{\xi}\leq 1$,
we have
$$
\lim_{k\ra \infty} \left[\supof{\norm{a}^{1\over k}: \mbox{
    $a$ is a product of $k$ elements from $\{a_i\}$}}\right] < \epsilon.
$$
This concludes the proof.
\qed
\begin{lemma}
  \label{Not radical}
  Let $A\subseteq\KH$ be a \cra\ such that $A''$ has no proper central
  projections.  Then $\lat A$ contains a non-zero irreducible submodule.
\end{lemma}
\proof If $\lat A$ contains no irreducible submodules then by
\nref{cts nest}, \nref{qn lemma} and \nref{Shulman} both $A$ and
$A\cdot A'$ consist of quasinilpotent compact operators. By
Lomonosov's lemma~\ref{Lomonosov} $A\cdot A'$ has a non-trivial
invariant subspace, and so by \nref{Commutant has RA} $A''$ contains a
proper central projection.  \qed

Having established the existence of an irreducible submodule $\V\in
\lat A$, the next step is to show that $A\cong \K(\V)\ot 1_n$ for some
integer $n$. The combination of Lomonosov's lemma and
section~\ref{Section with twine} makes this a fairly straightforward
task. To avoid circumlocution, for the remainder of this section we
will write `irreducible submodule' to imply `non-zero irreducible
submodule'.

\comment{
\Nref{Map from irreducible is zero} uses a technique adapted from \cite{Kissin}.
}
\begin{lemma}
\label{Map from irreducible is zero}
  Let  $A\subseteq \KH$ be a nondegenerate \cra, and suppose that $\v, \w\in \lat A$.
  If $\v$ is irreducible and $T:\v\ra \w$ is a non-zero module map, then
  the range of $T$ is closed and $T$ is an isomorphism onto its range.
\end{lemma}
\proof The (perhaps non-closed) restriction algebra $A|_\V$ is an
algebra of compact operators, such that $\V$ has no proper closed
submodules. Lomonosov's lemma implies that $A|_\V$ is weakly dense in
$\B(\V)$, and so $(A|_\V)'=\C1_\V$. Suppose $T:\V\ra \W$ is a non-zero
module map.  Replacing $\W$ with $\clos{T\V}$ we may assume that $T$
is dense-ranged.  By \nref{Relation is symmetric} there is a non-zero
module map $S:\W\ra \V$. Then $ST:\V\ra\V$ is a module map, non-zero
by the density of $T\V$. This means that $0\neq ST\in(A|_\V)'=\C1_\V$, and
so $T$ is bounded below. Thus, $T\V$ is closed and
$T$ is an isomorphism of $\V$ onto its range.
\qed

\comment{ Old Kissin proof, superseded...
As in 
\nref{Relation is symmetric}
we may assume that $\v\cap \w=\{0\}$ and $\v+\w$ is closed.

By the above \nameref{Compact and irreducible implies idempotent}
there is an idempotent $p\in A$ with non-zero restriction to $v$. 
Let $p_v$ and $p_w$ denote its
restriction to $\v$ and $\w$ respectively. Since $p_v$ is finite rank,
we may represent it as $p_v=\sum_{i\leq n} \xi_i\ot \zeta_i$.

Suppose that $T:\v\ra \w$ is a non-zero module map. Since $\v$ is
irreducible, $T$ must be one-to-one.
We will show that it is bounded below, and hence an isomorphism with
closed range. We
may assume that $T$ is dense-ranged, for if not we replace $\w$ with $\clos{T\v}$.

As $p_w$ is a restriction of a finite rank operator, it too has
finite rank, and so the 
equation $Tp_v=p_wT$ implies that $\ran p_w=\spam{T\xi_i}$.
Consider an arbitrary module map $S:\w\ra \v$. Since $Sp_w=p_vS$  
we have that $S$ maps $\spam{T\xi_i}$ into
$\spam{\xi_i}$, and so $ST$ maps $\spam{\xi_i}$ into itself. Since
$\spam{\xi_i}$ is finite-dimensional, there is an eigenvector
$\xi\in\spam{\xi_i}$ with eigenvector $\lambda$. Since $ST$ is a
module map on $\v$, we have $STa\xi=aST\xi=a\lambda\xi=\lambda a\xi$
for all $a\in A$, and hence $ST=\lambda\cdot I$ on $A\xi$. 
As $\v$ is irreducible, $A\xi$ is dense and so 
$ST=\lambda\cdot I$. Now if $T$
is not bounded below, then $\lambda=0$ and $ST=0$.  Since $T$ is
one-to-one and has dense range, this implies that $S=0$. Since $S$ was
an arbitrary module map from $\w$ to $\v$, \nref{Relation is symmetric}
implies that $T=0$. This contradiction shows that $T$ is bounded below.
\qed
}

\begin{lemma}
  \label{Map to irreducible is zero}
  Suppose $A\subseteq \KH$ is a nondegenerate \cra.  Let $\V\in\lat A$ be
  irreducible and $\W\in\lat A$ be arbitrary.  There is a non-zero module
  map $T:\w\ra \v$ \iff\ $\w$ contains a submodule isomorphic to $\v$.
\end{lemma}
\proof
Suppose that $T:\w\ra \v$ is a non-zero module map. 
\Nref{Relation is symmetric} implies that there is a non-zero module
map $S:\v\ra \w$, which is an \iso\ onto its range by 
\nref{Map from irreducible is zero}.

Conversely, if $\v'\subseteq \w$ is isomorphic to $\v$, there is a
module projection from $\w$ onto $\v'$. Composing this projection with the
isomorphism from $\v'$ to $\v$ yields the required non-zero module map.
\qed

\begin{lemma}\label{Sum of irreducible and other is closed}
  Let $A\subseteq\KH$ be a nondegenerate \cra.
  If $\v, \w\in \lat A$ and $\v$ is irreducible,
  then $\v+\w$ is closed.
\end{lemma}
\proof
If $\v\cap \w\neq \{0\}$ then $\v\cap \w=\v$ by irreducibility, and $\v+\w=\w$,
which is closed.

If $\v\cap \w=\{0\}$, we may find $\u\in \lat A$ such that $\u\oplus
\w=\H$. Let $p:\H\ra \u$ be the projection onto $\u$ along $\w$. Then
$T=p|_\V:\v\ra \u$ is a non-zero module map, and so, by \nref{Map from
  irreducible is zero}, $T$ is an isomorphism of $\v$ onto its range
contained in $\u$.  If $\v+\w$ is not closed there are $\xi\in \v,
\eta\in \w$ with $\norm{\xi}> 1$ and
$\norm{\xi+\eta}<\inv{\norm{\inv T}}\inv{\norm{p}}$. But then
$\norm{T(\xi)}=\norm{p(\xi)}=\norm{p(\xi+\eta)}<\inv{\norm{\inv T}}$,
whence $1<\norm{\xi}\leq \norm{\inv T}\norm{T(\xi)}<1$.  This
contradiction implies that $\v+\w$ is closed.  \qed

\begin{lemma}
\label{Compact and irreducible implies idempotent}
Let $A\subseteq\KH$ be nondegenerate \cra\ and suppose that $\v\in\lat
A$ is irreducible.  Then $A$ contains a non-zero projection which
restricts to a non-zero projection on $\v$.
\end{lemma}
\proof
Suppose that $\v\subseteq\H$ is irreducible, and
let $B$ denote the (possibly non-closed) algebra obtained by
restricting the operators of $A$ to $\v$. By Lomonosov's lemma, $B$
contains an operator with a non-zero eigenvalue, and hence so does
$A$. Since this operator is compact, we may use the holomorphic
functional calculus to obtain the desired idempotent \cite{BonsallDuncan}.
\qed

\comment{
We now want to reconstruct non-\sa\ versions of the results from
\csalg\ theory that say that faithful \sreps\ are Banach algebra
isomorphisms. Roughly speaking, the next two results show that this
also holds for
\reps\ of \cra{}s of compact operators involving irreducible modules.
}

\begin{lemma}
  \label{Only finite multiplicity}
  Let $A\subseteq \KH$ be a \oa\ with the \crp\
  and suppose  $\v\in \lat A$ is an irreducible submodule.
  Let $\W\subseteq\H$ be the closed span of a
  family of 
  submodules of $\H$ each isomorphic to $\v$. Then $\W$ is the direct sum of
  finitely many submodules \isoc\ to $\v$. 
\comment{
  The same conclusion holds
  for any submodule of $\W$.
}
\end{lemma}
\proof
Let $\V_1\subseteq\W$ be isomorphic to $\V$.
We inductively define  sequences 
$\{\v_i\}$ and $\{\w_i\}$ of submodules as follows:

If $\w_n=\cspam{\v_i}_{i\leq n}\neq \W$, then by the assumption on $\W$
there is a submodule $\v_{n+1}\subseteq \W$, isomorphic to $\v$ and not
contained in $\w_n$. Since $\v_{n+1}$ is irreducible, we have
$\v_{n+1}\cap \w_n=\{0\}$. 
By \nref{Sum of irreducible and other is closed} the
algebraic sum $\w_n+\v_{n+1}$ is closed, and so $\w_n+\v_{n+1}=\w_n\oplus
\v_{n+1}$. 

Using this to define inductively the sequences $\{\v_i\}$ and
$\{\w_i\}$, there are two possibilities: either $\w_n\varsubsetneq \W$
for all $n$ or $\w_n=\W$ for some $n$. However, the first case cannot
occur, for by \nref{Compact and irreducible implies idempotent} there
is an element $p\in A$ which restricts to a non-zero idempotent on
$\v$. Since each $\v_i$ is module isomorphic to $\v$, $p$ restricts to a
non-zero idempotent on $\v_i$.  As $p$ is compact and hence of
finite-rank, there can only be finitely many such $\v_i$ and so $\w_n=\W$
for some $n$.  By construction $\W=\w_n=\sum^\oplus_{i\leq n} \v_i$.
\qed

\comment{
\begin{lemma}
\label{At most one irreducible in simple algebra}
  Let $A\subseteq\KH$ has the \crp, and suppose that $A''$ has no central idempotents. 
  If $\v_1, \v_2\in\lat A$ are
  irreducible, then $\v_1$ and $\v_2$ are \isoc.
\end{lemma}
\proof
Suppose that $\v_1$ and $\v_2$ are not \isoc.
For $i\in \{1,2\}$ let $\X_i$ denote the closed span of all submodules of $\H$
\isoc\ to $v_i$. By \nref{Only finite multiplicity} $\X_i$ is the
direct sum of finitely many copies of $\v_i$.

Let $\w=\X_1\cap \X_2$. If $\w\neq \{0\}$,
then $\w$ is a finite direct sum of copies of
$\v_1$. Let $\v_1'\subseteq w$ be one such direct summand. Then 
$\v_1'$ is a direct sum of copies of $\v_2$
since $\v_1'\subseteq \w\subseteq \X_2$. However, $\v_1'$ is irreducible,
hence $\v_1'=\v_2'$ and $\v_1\cong \v_2$.

If $\w=\{0\}$, then 
\nref{Sum of irreducible and other is closed} implies
that $\X_1+\X_2$ is closed and so $\X_1+\X_2=\X_1\oplus \X_2$.
Let $\y$ be a module complement to
$\X_1\oplus \X_2$ in $\H$.
This induces a matrix representation of elements of
$A'$. If $\X_i$ is the direct sum of $n_i$ copies of $\v_i$,
then the commutant $A'$ has the matrix form
$$
\left[
  \begin{array}{ccc}
    M_{n_1} & 0 & 0 \\
    0 & M_{n_2} & 0 \\
    0 & 0 & (A|_y)' \\
  \end{array}\right],
$$
where the zeros are obtained from 
\nameref{Map from irreducible is zero}, 
\ref{Map from irreducible is zero} and
\ref{Map to irreducible is zero}.
From this it can be seen that
the matrix
$$
\left[
  \begin{array}{ccc}
    I & 0 & 0 \\
    0 & 0 & 0 \\
    0 & 0 & 0 \\
  \end{array}\right],
$$
is a central idempotent for $A'$, a contradiction.
\qed
}

\begin{proposition}
  \label{Compact algebra either continous or irreducible}
  Let $A\subseteq\KH$ be a \cra, and suppose $A''$ contains no proper
  central idempotents.  Then there exists an irreducible submodule
  $\V\in\lat A$, and $A$ is similar to $\Kom(\V)\ot1_n$ for some
  $n\in\N$.
\end{proposition}
\proof We have already seen that there is an irreducible submodule
$\V\in\lat A$.  Let $\W$ be the closed span of all \isoc\ copies of
$\v$ in $\lat A$.  From \nref{Only finite multiplicity}, $\W$ is a
finite direct sum of $n$ isomorphic copies of $\v$. If $\W\neq \H$
there is a non-zero $\U\in \lat A$ with $\H=\W\oplus \U$. By
definition $\U$ does not contain any isomorphic copies of $\v$.  In
this case the commutant $A'$ can be expressed in matrix form as
$$ \left[
  \begin{array}{cc}
    (A|_\W)' & 0 \\ 0 & (A|_\U)'
  \end{array}\right]
$$
where  the off-diagonal zeros are obtained from 
\nameref{Map from irreducible is zero}s~\ref{Map from irreducible is zero} and
\ref{Map to irreducible is zero}.
Thus the operator $$
\left[
  \begin{array}{cc}
    1 & 0 \\
    0 & 0
  \end{array}\right]
$$
is a central projection of $A''$. This contradiction implies that
$\W=\H$. 

This allows us to write $\H=\sum^\oplus_{i\leq n} \v_i$, where
$\v_i\cong \v$ via module isomorphisms $T_i:\V_i\ra V$.
We renorm $\H$ by $\norm{\sum^\oplus \xi_i}_{\rm
  new}^2=\sum\norm{T_i(\xi_i)}^2$. This renorming effects a similarity
on $\H$ under which $A$ is similar to $\K(\V)\ot 1_n$.
\qed

This establishes the desired result for those algebras $A\subseteq\KH$
with no proper central projections in $A''$. However, at this stage we
do not have a bound on the size of the similarity needed to implement
the \iso\ between $A$ and $\K(\V)\ot 1_n$ in terms of the \projconst\ 
for $A$. Clearly some such bound is necessary: if $\{A_i\}$ is a
sequence of \cra{}s with uniformly bounded \projconst{}s but which
require increasingly large similarities for the above \iso{}s, then
the direct sum $\sum^{c_0} A_i$ will be a \cra\ of compact operators
which is not similar to a \csalg. Happily, a modification of the proof
of \nref{Uniform bound on similarities} provides precisely the bound we
need. As for \nref{Uniform bound on similarities}, the proof is based on an analogous step in the proof of the
complemented subspaces theorem appearing in \cite{Day}.
\begin{lemma}\label{Bounded similarity for compact}
  Suppose that $A\subseteq\KH$ is a \cra\ such that $A''$ contains no
  proper central projections, and suppose $\pma\H \infty$ has
  \projconst\ $M$. Then there is a similarity $S$ on $\H$ with
  $\norm{S}\norm{\inv S}\leq 128 M^2$ such that $A^S$ is \sa.
\end{lemma}
\proof We know from \nref{Compact algebra either continous or
  irreducible} that there is an irreducible submodule $\V\in\lat A$
and $m\geq 1$ such that $\H$ is \isoc\ as an $A$-module to $\V\ot
\C^m$. Further,  for any $n\in\N$ the submodules of $\V\ot
\C^n$ are of the form $\V\ot \W$ where $\W\subseteq\C^n$ is an
arbitrary subspace. Consequently, any submodule of $\pma\V {2m}$ which
is module \isoc\ to $\pma\V m$ is actually isometrically \isoc\ to $\pma V
m$ as a module.

We define
$$
\alpha=\infof{\norm{S}\norm{\inv S}: \hbox{$S$ is a module \iso\ between
$\H$ and $\pma\V m$}}. 
$$
We may find a contractive module \iso\ 
$S:\H\ra\pma\V m$ with $\norm{\inv S}\leq 2\alpha$.
For $\mu\in\R^+$ we consider the graph subspace $\gr
\mu S\subseteq \H\oplus \pma\V m$. Observe that $\H\oplus\pma\V m$
is embedded isometrically in $\pma\H {m+1}$, and so $\H\oplus\pma V
m$ has the \rp\ with \projconst\ at most $M$. Let $p$ be a 
module projection from $\H\oplus\pma\V m$ onto $\gr \mu S$ with
$\norm{p}\leq M$.

Writing $p$ in components shows that
$$
p=\TxT 1+R \mu S, -R, \mu S(1+R \mu S), -\mu S R;
$$
for some module map $R:\pma\V m\ra \H$. The fact that $\norm{p}\leq M$
implies that $\norm{\mu S(1+R \mu S)}\leq M$ and $\norm{R}\leq M$.

We consider now the module map $T:\H\ra \pma\V m\oplus\pma\V m$ given by 
$$
T\xi={1\over 2}S \xi \oplus {1\over 2M} (\mu S(1+R \mu S)\xi).
$$
This is a contractive module isomorphism onto some closed submodule of
$\pma\V m\oplus\pma\V m$. By the above comments this image submodule
is isometrically \isoc\ to $\pma\V m$.  From the definition of
$\alpha$, this means that there is $\xi_0\in\H$ with $\norm{\xi_0}=1$
and $\norm{T\xi_0}\leq 2\inv \alpha$.

The remainder of the proof follows that of \nref{Uniform bound on
  similarities}. 
\comment{
Suppose that $\norm{S \xi_0}\leq \inv{(2M\mu)}$. From the second
term in the definition of $T$ we see that $\norm{T\xi_0}\geq
\mu/8M\alpha$. Since we know that $2\inv \alpha\geq \norm{T\xi_0}$ this
is impossible if we choose $\mu>16M$. 

Thus if we choose $\mu=16M+\epsilon$, we must have $\norm{S\xi_0}>
\inv{(2M\mu)}=\inv{(32M^2+2M\epsilon)}$. But the first term in the
definition of $T$ then gives us the inequality
$$
2\inv \alpha \geq \inv 2\norm{S\xi_0}>\inv{(64M^2 + 4M\epsilon)},
$$
and hence $\alpha\leq 128M^2$.}
\qed

We now have sufficient information to deal with general subalgebras
of $\K(\H)$. 
\begin{proposition}\label{A is c0 sum of simple}
  Suppose $A\subseteq\KH$ is a nondegenerate \cra.  Denote the set of
  minimal central projections of $A''$ by $P$.  For each $p\in P$ the
  algebra $A_p=pA$ is a closed two-sided ideal of $A$, and
  $A\cong\sum^{c_0}_{p\in P}A_p$. Moreover, considering $A_p$ as
  a subset of $\B(p\H)$, 
  the bicommutant $A_p{}''\subseteq\B(p\H)$ contains no proper central
  projections.
\end{proposition}
\proof Using \nref{Central projections are sa} we may assume that the
central projections of $A''$ are \sa. Let $R$ be the abelian \vna\ 
generated by these central projections.  \Nref{Discrete result} shows
that $R$ is generated by its minimal projections.  We denote the
minimal central projections of $A''$ by $P$ and the range of
$p\in P$ by $\H_p$. The nondegeneracy assumption implies
$\H=\sum^\oplus \H_p$ and $\sum_{p\in P} p=1$ (strong
convergence).

Let us write $A_p$ for the (possibly non-closed) algebra
$p A\subseteq\B(\H_p)$. 
Suppose that $q\in A_p{}''\subseteq\B(\H_p)$ 
is a central projection of $A_p{}''$ for some
$p\in P$. Then
$q\cdot p$ is central for $A''$, and $0\le
q\cdot p\leq
p$. Since $p$ is a minimal 
central projection of $A''$, either $q=0$ or $q=p$. In
either case, $q$ is not a proper central projection of $A_p{}''$.

We claim that $A_p$ is a closed ideal of $A$. To see
that $p A\subseteq A$, recall the duality $\bd
{\KH}=\TC{\H}^*=\BH$, where $\TC{\H}$ is the space of trace class
operators on $\H$ and the corresponding \ws\ topology on $\BH$ is the
\sw\ topology. \label{Page with HB}%
Suppose then that $p a\not\in A$ for
some $a\in A$. Since $p a$ is compact, by the Hahn-Banach theorem
there is $f\in \TC{\H}$ with $\<f, A>=0$ and $\<f,pa>=1$.
However, since $p a\in A''=\swclos{A}$, there is a net $\net{b_\mu}$
in $A$ which is \sw{}ly convergent to $p a$. 
Then we have $0=\<f,b_\mu>\ra \<f, p a>=1$. This
contradiction implies that $p A\subseteq A$, and that
$p A$ is norm closed (being the range of the projection
$p|_{A}$).  Since $p$ is central, $p A$ is a
two-sided ideal of $A$.

Let us write $\sum^{c_{00}}_{p\in P}A_p$ for the algebraic direct sum
(\ie\ considering elements with only finitely many non-zero terms). The
norm closure of $\sum^{c_{00}}A_p$ is $\sum^{c_0}A_p$.
Since the projections in $P$ are \sa\ and mutually
orthogonal, we have a natural isometric embedding $\sum^{c_{00}}
A_p\subseteq A$, and as $A$ is norm closed this implies that
$\sum^{c_0} A_p\subseteq A$.  On the other hand, as $\sum_{p\in P}
p=1$ the equality $a=\sum p a$ holds for all $a\in A$.
The compactness of $a$ then implies that $p\mapsto
\norm{p a}\in c_0(P)$ and so $A=\sum^{c_0}
A_p$.  \qed

\begin{theorem}\label{Compact algebra is csalg}
  Let $A\subseteq\KH$ be an \oa. Then $A$ has the \crp\ \iff\ $A$ is
  similar to a \csalg.
\end{theorem}
\proof Suppose that $A$ has the \crp.  If $A$ does not act
nondegenerately, there is a unique module $\V\in\lat A$ with
$\clos{A\H}\oplus\V=\H$ and $A\V=\{0\}$. Applying a similarity we may
arrange for $\clos{A\H}\perp\V$. Then $A$ will be \sa\ \iff\ 
$A|_{\clos{A\H}}$ is \sa, and so we may reduce to the case where $A$
acts nondegenerately.

By \nref{A is c0 sum of simple}, we may assume that there is an
orthogonal family $\H_\lambda$ of submodules of $\H$ such that 
$\H=\sum^\oplus\H_\lambda$ and $A=\sum^{c_0}A_\lambda$, where
$A_\lambda\subseteq\K(\H_\lambda)$ are \cra{}s and the bicommutants
$A_\lambda{}''$ have no proper central projections.  By \nref{Compact
  algebra either continous or irreducible}, for each $\lambda$ there is
a submodule $\V_\lambda\subseteq\H_\lambda$ and integer $n_\lambda$
such that $A_\lambda|_{\V_\lambda}=\K(\V_\lambda)$ and
$A_\lambda\cong\K(\V_\lambda)\ot 1_{n_\lambda}$.  Moreover, by
\nref{Bounded similarity for compact} the similarities
$S_\lambda\in\B(\H_\lambda)$ needed to make $A_\lambda$ \sa\ can be
chosen uniformly bounded.  The direct sum similarity $S=\sum^\oplus
S_\lambda$ will orthogonalise $A$.

Since all \sa\ \oas\ have the \crp\ the converse is immediate.
\qed

This result appears to be slightly better than might be expected from
the statements of \nameref{cra conjecture}s~\ref{cra conjecture} and
\ref{tra conjecture}, in that only the \crp\ is used to conclude that
a non-\sw{}ly closed \oa\ is similar to a \csalg. This apparent
anomaly is explained in terms of the biduality between $\K(\H)$ and
$\BH$. For suppose that $A\subseteq\KH$ is a \cra. Then by
\nref{Double commutant theorem} $A''$ is a
\sw{}ly closed \cra, and so according to \nref{cra conjecture} we
expect that $A''$ is similar to a \csalg. Consequently we expect that
the ideal $A''\cap \K(\H)$ is also similar to a \csalg. However, since
$A''=\swclos{A}$, it follows from the biduality $\bd {\K(\H)}=\BH$
that $A''\cap \K(\H)=A$. Thus \nref{Compact algebra is csalg} 
is exactly what we expect from  \nref{cra conjecture}.

\Nref{Compact algebra is csalg} can be extended to treat \cra{}s which
simply contain enough compact operators.
\begin{theorem}
  \label{sufficiently many compacts}
  If $A\subseteq\BH$ is a \sw{}ly closed \cra\ such that $A\cap
  \K(\H)$ acts nondegenerately on $\H$, then $A$ is similar to a
  \csalg. In fact, there is a family $\{\V_\lambda:
  \lambda\in\Lambda\}$ of \hs{}s and positive integers $\{n_\lambda:
  \lambda\in\Lambda\}$ such that $A\cong
  \sum^{\ell^\infty}\B(\V_\lambda)\ot 1_{n_\lambda}$.
\end{theorem}
\proof Let us denote $A\cap \K(\H)$ by $B$. Then the \sw\ closure of
$B$ is a \sw{}ly closed two-sided ideal of $A$, and by \nref{Ideals
  split} there is a central projection $p\in A\cap A'$ such that
$\swclos{B}=pA$. 
The nondegeneracy of $B$ implies that $p=1_\H$, and so
$\swclos{B}=A$. Thus by \nref{cra and weak closure} $B$ also has the \crp.
Then there is a family $\{\V_\lambda: \lambda\in\Lambda\}$ of \hs{}s
and positive integers $\{n_\lambda: \lambda\in\Lambda\}$ such that
$\H=\sum^{\oplus_2}\V_\lambda\ot \C^{n_\lambda}$ and $B\cong
\sum^{c_0}\K(\V_\lambda)\ot 1_{n_\lambda}$. Taking \sw\ closures shows
that $A$ contains the algebra $\sum^{\ell^\infty}\B(\V_\lambda)\ot
1_{n_\lambda}$.  Suppose that $T\in\BH$ does not lie in
$\sum^{\ell^\infty}\B(\V_\lambda)\ot 1_{n_\lambda}$. This implies that
either there is $\lambda\in \Lambda$ such that
$T\V_\lambda\not\subseteq\V_\lambda$ or there is $\lambda\in \Lambda$
and $i, j\leq n_\lambda$ with $T|_{\V_\lambda\ot e_i}\neq
T|_{\V_\lambda\ot e_j}$. In either case $T$ may be multiplied on the
right by a suitable operator from $B$ to obtain a compact operator
which does not belong to $B$. Thus $T\not\in A$ and so $A\cong
\sum^{\ell^\infty}\B(\V_\lambda)\ot 1_{n_\lambda}$.  \qed

\dealwithsectionbreaks
\section{Dealing with general \cra{}s}

If $A\subseteq\BH$ is a \sw{}ly closed \cra, then in order for $A$ to be
\sa\ it is necessary that the central projections of $A''$ be \sa. We
have seen that it is always possible to arrange for this by
application of a suitable similarity (\nref{Central projections are
  sa}). From here, the direct integral theory for non-\sa\ \oas\ given
in \cite{AzoffFongGilfeather} and discussed on page~\pageref{Page with
  direct integral} can be used to decompose $A$ into a
direct integral of \oas, each containing no proper central
projections. Thus the case where $A$ contains no proper central
projections is a test case for \nref{cra conjecture}. As we have seen,
for abelian $A$ this translates into the transitivity of $A'$.

A \vna\ which possesses no proper central projections is called a
factor. These are important in the theory of \vnas\ for the same
reason that they are important for \sw{}ly closed \cra{}s: every
\vna\ can be canonically and uniquely decomposed into a direct
integral of factors. Extending the \vna\ terminology, we will refer to
\sw{}ly closed \cra{}s which possess no proper central projections
as factors.

It does not immediately follow that a factor \cra\ has trivial centre.
\Nref{Commutant has RA} shows that if $A$ is a factor \cra, then the
algebra generated by $A$ and $A'$ is transitive. However, the centre
of $A$ will be trivial exactly when the algebra generated by $A$ and
$A'$ is $\BH$.

In the absence of a solution to \nref{cra conjecture} for factors, it
is appropriate to explore alternative means of understanding their structure.

The usual approach to understanding a \vna\ factor $R$ is  by examining its
(\sa) projections, or, equivalently, the invariant subspaces of the
factor $R'$.
The basic result for
factor \vnas\ can be stated in terms of $\lat R'$ as follows: if
$\V_1, \V_2\in \lat R'$ then there are non-zero submodules
$\W_1\subseteq\V_1$ and $\W_2\subseteq\V_2$ with $\W_1\cong\W_2$ as
$R'$-modules. This provides the basis for the comparability of
projections in $R$, and gives rise to the type classification of
factors \cite{Schwartz}.

It is interesting that to a certain extent a similar result holds for
\cra{}s. 
Recall from section~\ref{Section with twine} that the relation $\twine$ on
$\lat A$ is defined by 
$$\V\twine\W \Longleftrightarrow \mbox{ there is a non-zero
module map from $\V$ to $\W$.}
$$ 
It is convenient to say that $\V$ and $\W$ are disjoint if
$\V\not\twine\W$.
In \nref{Relation is symmetric} it was shown that $\twine$ is a
symmetric relation, in agreement with the \vna\ case. In further
accord with the \vna\ situation, if $A$ is a factor algebra then
every pair $\V, \W\in \lat A$ has $\V\twine\W$. The next few lemmas
lead up to the proof of this in \nref{All submodules are non
  disjoint}.

\begin{lemma}
  \label{Maximal disjoint}
  Suppose $A\subseteq\BH$ is a \cra.  Let $\V, \W\in\lat A$ be
  non-zero disjoint submodules. Then there are disjoint submodules
  $\V'\supseteq\V$ and $\W'\supseteq\W$, each maximal with respect to
  being disjoint from the other.
\end{lemma}
\proof Suppose that $\setof{\V_\alpha}$ is a chain of submodules each
disjoint from $\W$.  If $T:\clos{\bigcup\V_\alpha}\ra \W$ is a module
map, then the disjointness of $\V_\alpha$ and $\W$ implies
$T|_{\V_\alpha}=0$ for each $\alpha$. The boundedness of $T$ then
gives $T=0$, so $\clos{\bigcup\V_\alpha}$ is disjoint from $\W$.
Zorn's lemma now implies that there is a submodule $\V'\supseteq\V$
which is maximal with respect to being disjoint from $\W$.

Repeating this argument for $\V'$ and $\W$, there is $\W'\supseteq\W$
which is maximal with respect to being disjoint from $\V'$.
If $\V''\supset\V'$ 
is an arbitrary module strictly larger than $\V'$, 
then $\V''\twine\W$ by the maximality of
$\V'$. Hence $\V''\twine\W'$, so $\V'$ is maximally disjoint from $\W'$.
\qed

\begin{lemma}
\label{Factor implies summands are not disjoint}
Let $A\subseteq\BH$ be a \cra, and suppose $\H=\V\oplus\W$ is a module
direct sum decomposition of $\H$. If $A''$ has no proper central
idempotents then $\V\twine\W$.
\end{lemma}
\proof
If $\V\not\twine\W$ then by 
\nref{Relation is symmetric} we have $\W\not\twine\V$. This implies
that elements of $A'$ have the matrix form 
$$
\TxT *, 0, 0, *;
$$ with respect to the decomposition $\H=\V\oplus\W$. Then 
$$
\TxT 1, 0, 0, 0
;\in A''\cap A'.
$$
This contradicts the lack of proper central idempotents in $A''$.
\qed

\begin{lemma}
  \label{Disjoint implies closed sum}
Let $A\subseteq\BH$ be a \cra, and suppose $\V,\W\in\lat A$ are two disjoint
submodules. Then $\V+\W$ is closed.
\end{lemma}
\proof
Since $\V$ and $\W$ are disjoint we have $\V\cap\W=\{0\}$. Let $p$ be
a module projection from $\H$ onto $\V$. 
Then $p|_\W:\W\ra \V$ is a bounded module map, hence zero by
disjointness. This means that $\W\subseteq\ker p$, and so $\V+\W$
is closed.
\qed

\begin{proposition}
  \label{All submodules are non disjoint}
Let $A\subseteq\BH$ be a \cra. If $A''$ contains no proper central
idempotents then $\V\twine \W$ for any non-zero $\V, \W\in \lat A$.
\end{proposition}
\proof Suppose $\V, \W\in \lat A$ are non-zero submodules with
$\V\not\twine\W$. 
Using \nref{Maximal disjoint} we find $\V'\supseteq\V$ and
$\W'\supseteq\W$ maximal mutually disjoint covers of $\V$ and $\W$
respectively.  From \nref{Disjoint implies closed sum} the sum
$\V'+\W'=\V'\oplus\W'$ is a closed submodule of $\H$. If
$\V'\oplus\W'\neq \H$ then let $\U'$ be a module complement to
$\V'\oplus\W'$.  For any non-zero submodule $\U\subseteq\U'$, if
$\V'\not\twine\U$ then $\V'\not\twine\U\oplus\W'$, contradicting the
maximality of $\W'$. Thus $\V'\twine\U$ and by symmetry $\W'\twine\U$
for any non-zero submodule $\U\subseteq\U'$.  Let $\xi\in\U'$ be a
non-zero vector; then there is a non-zero module map $T:\V'\ra
\clos{A\xi}$, and consequently there is $0\neq\eta\in\V'$ with
$T\eta\neq 0$. Since $\clos{AT\eta}\twine\W'$,  there is a non-zero
module map $S:\clos{AT\eta}\ra \W'$, with $0\neq S(aT\eta)=aS(T\eta)$
for some $a\in A$. Hence $ST\eta\neq0$ and $ST:\V'\ra\W'$ is a
non-zero module map. This contradiction shows that $\V'\oplus\W'=\H$.
Since $\V'$ and $\W'$ are disjoint \nref{Factor implies summands are
  not disjoint} implies that $A''$ has a proper central projection.
\qed

\Nref{All submodules are non disjoint} suggests that it might be
possible to make a type analysis of factor algebras. Accordingly, if
$A$ is a factor algebra, we say that $A$ is type $I$ if $\lat A$
contains an irreducible submodule (\ie\ if $A'$ contains a minimal
projection).  In this case, if $\H$ is the direct sum of $n$
irreducible modules we say that $A$ is of type $I_n$; otherwise we
say $A$ is of type $I_\infty$.  

It is interesting to consider the
barriers to showing that $A$ is similar to a \vna\ 
in terms of its type.
We expect that type $I$ algebras should be the least difficult to
analyse. This is particularly so since if a transitive \oa\ 
$A\subseteq\BH$ is known to be similar to a \sa\ algebra, then $A$ is
already \sa\ (it is equal to $\BH$, after all). This means that the
technical problem of selecting a particular involution on $A$ is
avoided. 
We have already seen an example of this phenomenon in \nref{Compact
  algebra either continous or irreducible}. 

The problem of dealing with general type $I$ algebras reduces to solving a
certain case of the transitive algebra problem, as follows:
\begin{question}
  \label{special transitive algebra problem}
  Let $A\subseteq\BH$ be a transitive algebra with the \crp. Does this
  imply $A=\BH$?
\end{question}

If $A$ is of type $I_n$, then $\H$ is the direct sum of $n$
irreducible submodules.  In this situation a necessary and sufficient
condition for $A$ to be similar to a \sa\ algebra is that the answer
to \nref{special transitive algebra problem} is `yes'.

Moreover, a positive answer to \nref{special transitive algebra
  problem} would also deal with the type $I_\infty$ case, where $\H$
is the span of infinitely many copies of the irreducible submodule. In
this case, the argument of \nref{Bounded similarity for compact} could
be recycled to establish a uniform bound in the size of the
similarities needed to orthogonalise an arbitrary finitely generated
submodule of $\H$, and a standard ultrafilter argument used to piece
the finitely generated parts together (\cf~\nref{ultrafilter
  argument}). 

The analysis of non-type~$I$ factors appears to constitute a much greater
challenge. The lack of irreducible submodules means that a
`non-classical' generalisation of the complemented subspaces theorem
seems to be necessary in order to establish step (i) of the program
outlined on page~\pageref{Page with program}. I suspect that if there
is hope for a proof of \nref{cra conjecture}, it lies in finding a
suitable generalisation of the complemented subspaces theorem, and
using \nref{RP -> reductive is desirable} to sidestep directly
addressing the reductive algebra problem.


\vfill\eject

\bibliographystyle{plain}

\bibliography{thesis}

\begin{thebibliography}{10}

\bibitem{ArvesonBook}
W.~Arveson.
\newblock {\em An Invitation to \cs\ Algebra}.
\newblock Springer-Verlag, New York, Heidelberg, Berlin, 1976.

\bibitem{AzoffFongGilfeather}
E.~A. Azoff, C.~K. Fong, and F.~Gilfeather.
\newblock A reduction theory for non-self-adjoint operator algebras.
\newblock {\em Trans. Amer. Math. Soc.}, 224(no. 2):351--377, 1976.

\bibitem{BonsallDuncan}
F.~F. Bonsall and J.~Duncan.
\newblock {\em Complete Normed Algebras}.
\newblock Springer-Verlag, Berlin, Heidelberg, New York, 1973.

\bibitem{cb-operators}
E.~Christensen and A.~M. Sinclair.
\newblock A survey of completely bounded operators.
\newblock {\em Bull. London Math. Soc.}, 21(5):417--448, 1989.

\bibitem{CurtisLoy}
P.~C. Curtis~Jr and R.~J. Loy.
\newblock The structure of amenable {B}anach algebras.
\newblock {\em J. London Math. Soc.}, 40:89--104, 1989.

\bibitem{Davidson}
K.~R. Davidson.
\newblock {\em Nest Algebras}.
\newblock {P}itman Research Notes in Mathematics Series. Longman Scientific \&
  Technical, 1988.

\bibitem{DavisDeanSinger}
W.~J. Davis, D.~W. Singer, and I.~Singer.
\newblock Complemented subspaces and {$\Lambda$} systems in {B}anach spaces.
\newblock {\em Israel J. Math.}, 6:303--309, 1968.

\bibitem{Day}
M.~M. Day.
\newblock {\em Normed Linear Spaces}.
\newblock Springer Verlag, Berlin, Heidelberg, New York, 3rd edition, 1973.

\bibitem{DixmierC}
J.~Dixmier.
\newblock {\em \csalgs}.
\newblock North Holland Publishing Company, Amsterdam, New York, Oxford, 1977.

\bibitem{DixmierV}
J.~Dixmier.
\newblock {\em Von {Neumann} Algebras}.
\newblock North Holland, Amsterdam, 1981.

\bibitem{DixonTrick}
P.~G. Dixon.
\newblock Approximate identities in normed algebras.
\newblock {\em Proc. London Math. Soc.}, 26(3):485--496, 1973.

\bibitem{Edwards}
R.~E. Edwards.
\newblock {\em Functional Analysis}.
\newblock Holt, Rinehart and Winston, New York, Chicago, San Francisco,
  Toronto, London, 1965.

\bibitem{Effros-quantize}
E.~G. Effros.
\newblock Advances in quantized functional analysis.
\newblock In {\em Proceedings of the International Congress of Mathematicians},
  pages 906--916, 1986.

\bibitem{EffrosRuan}
E.~G. Effros and Z.-J. Ruan.
\newblock On non-self-adjoint operator algebras.
\newblock {\em Proc. Amer. Math. Soc.}, 110(4):915--922, December 1990.

\bibitem{ElliotMultiplier}
G.~A. Elliot.
\newblock An abstract {D}auns--{H}offman--{K}aplansky multiplier theorem.
\newblock {\em Canad. J. Math.}, {XXVII}(4):827--836, 1975.

\bibitem{ElliotOleson}
G.~A. Elliot and D.~Oleson.
\newblock A simple proof of the {Dauns}--{H}ofmann theorem.
\newblock {\em Math. Scand.}, 34:231--234, 1974.

\bibitem{Fong}
C.~K. Fong.
\newblock Operator algebras with complemented invariant subspace lattices.
\newblock {\em Indiana Univ. Math J.}, 26(6):1045--1056, 1977.

\bibitem{Gardner}
L.~T. Gardner.
\newblock On isomorphisms of {\csalgs}.
\newblock {\em Amer. J. Math.}, 87:384--396, 1965.

\bibitem{Gourdeau}
F.~Gourdeau.
\newblock Amenability of {L}ipschitz algebras.
\newblock {\em Math. Proc. Camb. Phil. Soc.}, 112:581--588, 1992.

\bibitem{GronbaekWillis}
N.~Gr{\o}nb{\ae}k and G.~A. Willis.
\newblock Embedding nilpotent finite-dimensional {B}anach algebras into
  amenable {B}anach algebras.
\newblock {\em J. Funct. Anal.}, 145(1):99--107, 1997.

\bibitem{HaagerupCyclic}
U.~Haagerup.
\newblock Solution of the similarity problem for cyclic representations of
  {\csalgs}.
\newblock {\em Ann. Math}, 118:215--240, 1983.

\bibitem{Helemskii}
A.~Ya. Helemski.
\newblock {\em The Homology of {B}anach and Topological Algebras}.
\newblock Kluwer Academic Publishers, Dordrecht, Boston, 1989.

\bibitem{Helemskii-survey}
A.~Ya. Helemskii.
\newblock Flat {B}anach modules and amenable algebras.
\newblock {\em Trans. Moscow Math. Soc.}, 47:179--218, 1984.
\newblock Amer. Math. Soc. Translations, (1985) 199--224.

\bibitem{HewittRoss}
E.\ Hewitt and K.~A. Ross.
\newblock {\em Abstract Harmonic Analysis}, volume~{II}.
\newblock Springer-Verlag, Berlin, Heidelberg, New York, 1970.

\bibitem{Johnson1}
B.~E. Johnson.
\newblock {\em Cohomology in {B}anach Algebras}, volume 127 of {\em Mem. Amer.
  Math. Soc.}
\newblock Amer. Math. Soc., Providence, R.I., 1972.

\bibitem{JohnsonPerturbation}
B.~E. Johnson.
\newblock Perturbations of {B}anach algebras.
\newblock {\em Proc. London Math. Soc.}, 34(3):439--458, 1977.

\bibitem{Kadison}
R.~V. Kadison.
\newblock On the orthogonalization of operator representations.
\newblock {\em Amer. J. Math.}, 77:600--620, 1955.

\bibitem{LauLoyWillis}
A.~T.-M. Lau, R.~J. Loy, and G.~A. Willis.
\newblock Amenability of {B}anach and {$C^*$}-algebras on locally compact
  groups.
\newblock {\em Studia Math.}, 119(2):161--178, 1996.

\bibitem{LindenstraussTzafriri}
J.~Lindenstrauss and L.~Tzafriri.
\newblock On the complemented subspaces problem.
\newblock {\em Israel J. Math.}, 9:263--269, 1971.

\bibitem{Palmer}
T.~Palmer.
\newblock {\em {B}anach Algebras and the General Theory of *-Algebras}.
\newblock Cambridge University Press, 1994.

\bibitem{Paterson}
A.~L. Paterson.
\newblock {\em Amenability}.
\newblock Number~29 in Mathematical Surveys and Monographs. AMS, Providence,
  Rhode Island, 1988.

\bibitem{Paulsen}
V.~Paulsen.
\newblock {\em Completely Bounded Maps and Dilations}.
\newblock Pitman Research Notes in Mathematics. Longman Scientific and
  Technical, Essex, 1986.

\bibitem{PisierSimilarity}
G.~Pisier.
\newblock {\em Similarity Problems and Completely Bounded Maps}, volume 1618 of
  {\em Lecture Notes in Mathematics}.
\newblock Springer, 1996.

\bibitem{PisierDegree}
G.~Pisier.
\newblock Le degr\'e de similarit\'e d'une alg\`ebre d'op\'erateurs.
\newblock {\em C. R. Acad. Sci. Paris S\'er. I Math.}, 324(3):287--292, 1997.

\bibitem{RadjaviRosenthal}
H.~Radjavi and P.~Rosenthal.
\newblock {\em Invariant Subspaces}.
\newblock Springer-Verlag, Berlin, Heidelberg, New York, 1973.

\bibitem{RaeburnPerturbation}
I.~Raeburn and J.~L. Taylor.
\newblock {H}ochschild cohomology and perturbations of {B}anach algebras.
\newblock {\em J. Funct. Anal.}, 25(3):258--266, 1977.

\bibitem{StormerIdeals}
E.~St\o rmer.
\newblock Two-sided ideals in {\csalgs}.
\newblock {\em Bull. Amer. Math. Soc.}, 73:254--257, 1967.

\bibitem{Rosenoer1}
S.~Rosenoer.
\newblock Completely reducible operator algebras and spectral synthesis.
\newblock {\em Canad. J. Math.}, XXXIV(5):1025--1035, 1982.

\bibitem{Rosenoer2}
S.~Rosenoer.
\newblock Completely reducible operators that commute with compact operators.
\newblock {\em Trans. Amer. Math. Soc.}, 299(1):33--40, 1987.

\bibitem{RosenoerAnnoying}
S.~Rosenoer.
\newblock Completely reducible algebras containing compact operators.
\newblock {\em J. Operator Theory}, 29:269--285, 1993.

\bibitem{Schwartz}
J.~T. Schwartz.
\newblock {\em {$W^*$} Algebras}.
\newblock Thomas Nelson and Sons, 1968.

\bibitem{Willis-operator-normal}
G.~A. Willis.
\newblock When the algebra generated by an operator is amenable.
\newblock {\em J. Operator Theory}, 34:239--249, 1995.

\end{thebibliography}

\printindex

\end{document}